\documentclass[oneside]{elsarticle}
\usepackage{etoolbox}
\usepackage{graphicx}
\usepackage{caption}
\usepackage{subcaption}
\usepackage{amssymb}
\usepackage{amsmath}
\usepackage{booktabs}
\usepackage{xcolor} 
\usepackage{multirow}
\usepackage{makecell}
\usepackage[colorlinks, citecolor=blue]{hyperref}
\graphicspath{ {./} }
\usepackage{scrextend}
\usepackage{fancyhdr}
\usepackage{graphicx}
\usepackage{xcolor,colortbl}
\definecolor{Red}{RGB}{230,4,0}
\definecolor{Green}{RGB}{0,175,80}
\definecolor{Yellow}{RGB}{255,192,0}
\usepackage{subcaption}
\usepackage{algorithm2e}
\usepackage{lineno}
\usepackage[normalem]{ulem}
\usepackage{setspace}
\newtheorem{remark}{Remark}[section]

\SetKwInput{KwInput}{Input}
\SetKwInput{KwOutput}{Output}

\biboptions{numbers,sort&compress}

\author[LMPS,ENS]{Thomas Melkior}
\ead{thomas.melkior@ens-paris-saclay.fr}

\author[ENS]{Harsha S. Bhat}
\ead{harsha.bhat@ens.fr}

\author[LMPS]{Faisal Amlani\corref{cor1}}
\ead{faisal.amlani@ens-paris-saclay.fr}
\cortext[cor1]{Corresponding author: faisal.amlani@ens-paris-saclay.fr}
\address[LMPS]{Université Paris-Saclay, CentraleSupélec, ENS Paris-Saclay, CNRS, LMPS - Laboratoire de Mécanique Paris-Saclay, Gif-sur-Yvette, France}

\address[ENS]{Laboratoire de G\'eologie,  \'{E}cole Normale Sup\'erieure, CNRS-UMR 8538, PSL Research University, 75006 Paris, France}

\title{Tsunami modeling with dynamic seafloors: a high-order solver validated with shallow water benchmarks
}
\date{Draft updated: \today}

\def\equationautorefname#1#2\null{%
  Eq.#1(#2\null)%
}

\begin{document}

% \pagestyle{headings}	
% \newpage
% \setcounter{page}{1}
% \renewcommand{\thepage}{\arabic{page}}
	
% \captionsetup[figure]{labelfont={bf},labelformat={default},labelsep=period,name={Figure }}	\captionsetup[table]{labelfont={bf},labelformat={default},labelsep=period,name={Table }}
% \setlength{\parskip}{0.5em}
	
% \maketitle

	\begin{abstract}

      Recent scientific studies have suggested that, in certain physical configurations, the time-dependent behavior of earthquake rupture and seafloor (bathymetry) motion can leave observable near-field signatures in tsunami wave generation and propagation. 
    However, dynamic ground movement is often neglected in conventional tsunami models, which commonly assume instantaneous ground displacement (sourcing). 
      This work introduces a pseudo-spectral algorithm for the solution of the nonlinear shallow water equations with time-dependent seafloor displacement and velocity. Based on a Fourier continuation (FC) methodology for the accurate trigonometric interpolation of a non-periodic function, the solver provides  high-order convergence in space and time; mild (linear) Courant-Friedrichs-Lewy (CFL) constraints for explicit time integration; and results that are effectively free of numerical dispersion (or "pollution") errors. Such properties enable the efficient and robust resolution of the different space-time scales involved modeling tsunamis generated by dynamic earthquake ground motion (including over long distances).  Numerical experiments attesting to accuracy and computational performance are presented with direct comparisons to high-order finite difference methodologies. The solver is physically validated by a number of classical and semi-classical benchmark cases based on simulated or experimental data. Additionally, a seismologically realistic, first-of-its-kind parametric study based on earthquake speed is introduced, whose results---easily facilitated by the FC-based approach proposed herein, with minimal numerical tuning---further demonstrate the potential importance of (and the motivation herein for) incorporating time-dependent seafloor behavior in quantitative tsunami hazard assessment.
       
    %   This work introduces a pseudo-spectral algorithm for the solution of the nonlinear shallow water equations with time-dependent seafloor displacement and velocity (which is often neglected in conventional

% 		\let\thefootnote\relax\footnotetext{
% 			\small $^{*}$\textbf{Corresponding author.}}
% 		\textbf{\textit{Keywords}}: \textit{Tsunami; Numerical simulation; Shallow water equations; Benchmark; Spectral methods}
	\end{abstract}
	

\maketitle

\section{Introduction}
The human and economic cost of tsunamis has been increasing in recent decades \cite{Imamura2019}. In order to better understand the underlying mechanisms of tsunami generation and propagation, numerical modeling has been indispensable, including for real-world events~\cite{kowalik2007, sugawara2021, dynamicsource1}. In particular, using physics-based simulations, it has been only recently shown~\cite{amlanibhat2022,Elbanna2021} that certain strike-slip earthquakes (which have minimal vertical displacement) can create significant waves, e.g., the deadly 2018 Palu bay, Indonesia tsunami~\cite{amlanibhat2022} that was mysteriously and unexpectedly generated from such a fault (as opposed to from a classical tsunami-generating subduction zone, where the final, static, non-negligible vertical ground displacement is singularly important). This unusual tsunamigenesis was explained in a proof-of-concept work~\cite{amlanibhat2022} by incorporating the time-dependent effects of so-called "supershear" earthquake rupture displacement and velocity into a simplified fluid-structure model. Such a finding has uncovered possible new risks (where there was thought to be little before) in a host of similar fault locations and configurations elsewhere~\cite{amlanibhat2022}, ultimately implying that there may be signatures yet to be discovered in fully dynamic considerations of seismogenic tsunami generation (including those by subduction zone/thrust fault ruptures). While others have also suggested the importance of considering dynamics in the near-field~\cite{dutykh2016, dynamicsource1}, further study of this hypothesis requires fast, accurate, and robust numerical tools capable of more full-order, physically-faithful seismogenic tsunami modeling on realistic large-scale computational domains subjected to dynamic ground motion.

% or the unusual 2018 Palu tsunami \cite{amlanibhat2022}

Conventionally, tsunami propagation is calculated at different scales: near-fields where the depth is low and where many strong nonlinearities can take place (e.g. wave breaking), and far-fields where the depth is high and some of these nonlinearities can be neglected \cite{TUNAMI}. Comprehensive tsunami solvers often use nested meshes where the choice of numerical method for each mesh depend on the scale \cite{TUNAMI}. At the near-field, methods able to treat can easily sharp gradients are better suited, such as the finite volume method \cite{LeVeque2011, Delestre2015}, the lattice Boltzmann method \cite{Zhou2002}, or the smooth particle hydrodynamics (SPH) method \cite{Roselli2019, Wei2015}. However, many of these methods can become costly for larger-scale configurations \cite{TUNAMI}. For long-distance and long-time tsunami wave propagation sourced by earthquakes, more computationally-efficient methods are often used, such as finite elements (FE) \cite{FEMtsunami} or finite differences (FD) \cite{TUNAMI, NAMIDANCE, HEINRICH2021}. However, numerical dispersion or diffusion ("pollution") errors are well-known to be problematic for FD-based \cite{Durran1999} or FE-based \cite{Deraemaeker1999} solvers. For tsunami problems, this numerical pollution manifests as an accumulation of errors that grows as waves propagate further through space and time, eventually leading to considerable loss of accuracy. Consequently, such methods often require significant refinement that can quickly become prohibitively expensive in terms of computational cost for large-scale problems, potentially undermining their efficacy. Additionally, there are very few tsunami solvers in general that can consider time-dependent ground displacement and velocity induced by a rupture; those that can are often low-order accurate or too computationally expensive and sometimes only consider simple mathematical approximations to the full-order dynamic ground behavior~\cite{dynamic1,dynamic2} (e.g., for black-box implementations). In order to fully incorporate dynamics, existing and commonly-used solvers require sometimes difficult or impractical modifications to the fundamental governing equations as well as their associated numerical methodologies: configurations involving such time-dependent sources (which can be of very high-frequency compared to the tsunami~\cite{amlanibhat2022}) bring additional numerical and computational challenges in resolving all scales.
% (e.g., the additional high-frequency and time-dependent complexities can be more sensitive to the method of discretization or solution).} 

Hence the objective of this work is to introduce a new nonlinear shallow water equation (tsunami) solver based on a Fourier continuation (FC) approach for the trigonometric interpolation of non-periodic functions \cite{amlani2016,bruno2010}. Such an approach extends the applicability of Fourier series representations to general functions and boundary conditions while minimizing the well-known Gibbs "ringing effect" and enabling the use of the fast Fourier transform (FFT). FC has been employed successfully as a spectral method for solving various hyperbolic or parabolic partial differential equation (PDE) systems~\cite{bruno2010,AmlaniPahlevan2020,bruno2025,gaggioli2019,Fontana2020,Bruno2022,albin2011,amlani2016,Amlani2023,amlaniPOF,arian2025}, achieving high accuracy by means of relatively coarse discretizations, a faithful preservation of the dispersion or diffusion characteristics of the underlying continuous operators (minimal numerical pollution errors), and mild (linear) CFL constraints on explicit time integration (properties that have facilitated applications to challenging scientific and engineering problems in both solid and fluid dynamics~\cite{amlani2016, Amlani2023, amlaniPOF, arian2025, bilgi2023,amlanibhat2022}). These features are well suited for resolving shallow water equations since tsunami waves can propagate over long distances and long times. Additionally, such high-order capabilities are essential for the resolution of the different temporal scales that might be involved in modeling earthquake-induced tsunamis. Indeed, a particular novelty of the proposed solver is that it treats a more generalized or enriched version \cite{dutykh2016} of the classical shallow water equations in order to account for dynamic ground motion. Most conventional solvers employ static seafloor displacement, which is better suited for long wavelength far-field tsunami models \cite{wang2006analysis,abrahams2023comparison}; however, as described above, the dynamics of ground displacement and velocity may represent an important mechanism in some tsunamigenesis scenarios, particularly in the near-field~\cite{amlanibhat2022, Elbanna2021, dynamicsource1}.

Following an encouraging 1D proof-of-concept work on Palu~\cite{amlanibhat2022}, which demonstrates the potential of FC in addressing such challenges (the authors built a preliminary 1D solver for that specific case only), this contribution presents a general, full-fledged, high-order 1D and 2D tsunami solver that treats moving seafloors with high accuracy, FFT-speed, and minimal numerical pollution. The ultimate goal is to enable efficient simulations towards further parametric studies and explorations of the dynamic effects of earthquake ground motion on tsunamigenesis. After introducing the governing equations in \autoref{Governing equations}, the solver based on FC is detailed in \autoref{Methodology}. \autoref{Numerical study} presents a comprehensive analysis of the numerical performance of the solver (attesting to high accuracy and efficiency), including comparisons to fourth- and sixth-order finite difference methods. \autoref{Benchmarks} provide assessments of the solver on a number of realistic physical configurations that have been proposed extensively as benchmarks in literature (including comparisons to both experimental and simulated data). \autoref{sec:Application} presents the application of the FC-based solver to a newly-proposed parametric problem configuration with highly-dynamic seismic sources of different surface speeds, providing insight on the relevance of the new solver on scenarios that highlight the importance of time-dependent seafloor displacement. Concluding remarks are discussed in \autoref{Conclusions}.

% \red{The dispersion errors of finite differences will be highlighted in \autoref{Numerical study}, it will furthermore be observed that these errors are not resolved by increasing the order of the scheme.}

\section{Governing equations} \label{Governing equations}

The high-order solver introduced in this work treats a physical model representing the behavior of tsunamis that is based on the shallow water equations (SWE) of Saint-Venant \cite{SaintVenant1871}, which are a depth-averaged version of the Euler equations. The validity of these shallow water (tsunami) equations to model fluid flows rely on various principal hypotheses \cite{Kundu2007,dutykh2007mathematical}:
\begin{itemize}
    \item The water depth $h(x,y,t)$ is small with respect to the wavelength of water waves.
    \item The fluid is incompressible.
    \item The flow is irrotational.
    \item The horizontal components of the velocity fields $u(x,y,t)$ (along $x$) and $v(x,y,t)$ (along $y$) are constant along the vertical component $z$.
\end{itemize}
In the present contribution, the system is treated in Cartesian coordinates---a valid approximation when propagation distances remain reasonable compared to the earth curvature~\cite{Kundu2007,dutykh2007mathematical}. For larger-scale propagation, the system can be considered in spherical coordinates, for which the methodology introduced in this work can be straightforwardly extended. 

\begin{figure}[!t]
	\centering
	\includegraphics[width=0.5\textwidth]{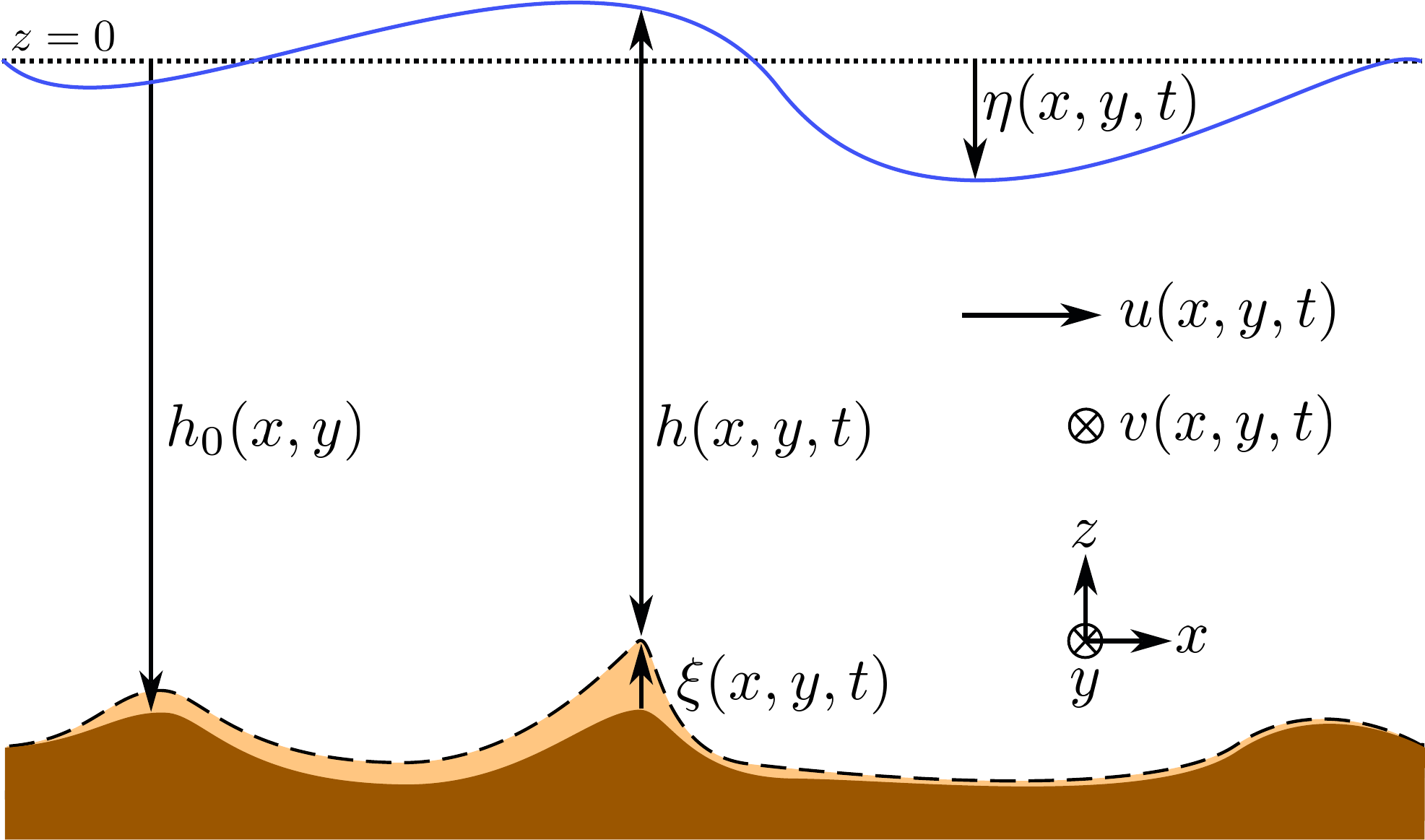}
	\caption{A representative illustration of the mathematical notations and physical definitions for the governing shallow water systems treated in this work.}
	\label{fig:SWE}
\end{figure}

The fundamental physical configuration is illustrated in \autoref{fig:SWE}, where the total water depth $h(x,y,t)$ is decomposed as $h(x,y,t)=h_0(x,y)+\eta(x,y,t)-\xi(x,y,t)$ for given still water depth $h_0(x,y)$ corresponding to a given unperturbed sea floor geometry (bathymetry) $h_0(x,y)$, a given ground motion vertical displacement $\xi(x,y,t)$ (the source), and unknown wave height $\eta(x,y,t)$ (i.e., displacement with respect to the at-rest free surface). The vertical displacement of the sea floor can be considered by different contributions \cite{dutykh2016}, here, a first order approximation that neglects horizontal displacements is implemented. Commonly, $\xi$ is enforced by an instantaneous water uplift \cite{Setiyowati2019} (i.e., an instantaneous free surface perturbation corresponding to $\eta(x,y,t=0)=\max_t \xi(x,y,t)$). However, as discussed previously, the interest of this work is to develop a solver that can account for ground motion time history in order to model and study dynamic earthquake source effects (represented by displacement $\xi(x,y,t)$ and velocity $\xi_t(x,y,t)$\footnote{For compactness of notation, partial derivatives in $x,y$, and $t$ are presented in subscript form throughout this work, e.g., $\partial \xi/\partial t (x,y,t) = \xi_t(x,y,t)$.}) on tsunamigenesis. For the present work, the unperturbed sea depth $h_0(x,y)$ and the perturbation $\xi(x,y,t)$ (with its corresponding velocity  $\xi_t(x,y,t)$) are given \emph{a priori} (the latter possibly by separate rupture simulations or geophysical data~\cite{amlanibhat2022,Elbanna2021}).

Although the focus of the solver is 2D, the 1D equations are first introduced in what follows for facilitating the presentation and for completeness with respect to the 1D numerical experiments and benchmark cases that are discussed in \autoref{Numerical study} and \autoref{Benchmarks}. This work considers all systems in a classical nonlinear and nondispersive form, although the methodology can be easily extended to include additional terms such as those modeling friction, viscosity, and Coriolis effects.

\subsection{1D shallow water PDE system}
In the 1D formulation of the SWE configuration shown in \autoref{fig:SWE}, the transverse dimension $y$ is considered invariant: the quantities do not depend on $y$ and hence velocity in that direction is also neglected, i.e. $v=0$. The corresponding 1D shallow water PDE system is given by \cite{amlanibhat2022}
\begin{equation}
    \begin{cases}
        h_t+(hu)_x=0, \\ (hu)_t+(hu^2)_x+gh \eta_x=0,
    \end{cases}
\end{equation}
where $g$ is acceleration due to gravity (assumed to be 9.81 m/s$^2$ in this work unless indicated otherwise), where $h=h(x,t)=h_0(x)+\eta(x,t)-\xi(x,t)$, and where the unknowns are given by velocity $u=u(x,t)$ and free-surface displacement $\eta = \eta(x,t)$ for $x\in \Omega$ of spatial domain $\Omega \subset \mathbb{R}$, and a time interval $t\in[0, T]$ for $T\in \mathbb{R}^+$. The corresponding wave speed, $c = \sqrt{gh}$, can be approximated via a linearization of such equations~\cite{Kundu2007,dutykh2007mathematical,SaintVenant1871}. Expanding, rearranging, and assuming $h>0$~\cite{dutykh2013use} yields the system given by
\begin{equation} \label{1D SWE}
    \begin{cases}
        \eta_t - \xi_t +(hu)_x=0, \\ u_t+ u u_x+g \eta_x=0, \\ h(x,t)=h_0(x)+\eta(x,t)-\xi(x,t),
    \end{cases}
\end{equation}
where the (given) perturbation of the sea floor $\xi(x,y,t)$ appears in both the $h$ and $h_t$ terms (again, such time-dependent ground motion is often neglected for tsunami modeling~\cite{amlanibhat2022,dutykh2007mathematical}). In vectorial form, defining $\mathbf{U}=(\eta~u)^\text{T}$, this can be written as
\begin{equation} \label{Vectorial 1D SWE}
    \mathbf{U}_t+\mathbf{F}(\mathbf{U},\mathbf{U}_x)=\mathbf{S},
\end{equation}
where
% $$\mathbf{F}(\mathbf{U},\mathbf{U}_x)=\begin{pmatrix}h_x u+h u_x\\u u_x +g \eta_x \end{pmatrix}=\begin{pmatrix}(\eta_x+h_{0x}-\xi_x)u+(\eta+h_0-\xi)u_x\\u u_x +g \eta_x \end{pmatrix} \quad \text{and} \quad \mathbf{S}=\begin{pmatrix}\xi_t\\0\end{pmatrix}.$$
$$\mathbf{F}(\mathbf{U},\mathbf{U}_x)=\begin{pmatrix}(\eta_x+h_{0x}-\xi_x)u+(\eta+h_0-\xi)u_x\\u u_x +g \eta_x \end{pmatrix} \quad \text{and} \quad \mathbf{S}=\begin{pmatrix}\xi_t\\0\end{pmatrix}.$$
The well-posedness of this system is completed by various Dirichlet, Neumann, or non-reflecting boundary conditions. The two classical boundary conditions for the shallow water system considered in this work are the so-called \textit{wall boundary conditions} given by \cite{bristeau2001boundary}
\begin{equation} \label{eq:1Dwallbc}
    \begin{cases}
        \eta_x(x,t)=0, \\
        u(x,t)=0,
    \end{cases} \quad x\in{\partial \Omega},~~t\in [0,T],
\end{equation}
and the \textit{radiation} (non-reflecting) boundary conditions given by \cite{Durran1999}
\begin{equation} \label{eq:1Dradiationbc}
    \begin{cases}
        \eta_t(x,t)+c_n \eta_x(x,t)=0, \\
        u_t(x,t)+c_n u_x(x,t)=0,
    \end{cases} \quad x\in{\partial \Omega},~~t\in [0,T],
\end{equation}
where $c_n=\pm \sqrt{gh}$ is the propagation speed directed outward to the boundary at the right endpoint (i.e., $x = x_\text{max}$) and left endpoint (i.e., $x=x_\text{min}$), respectively.
\subsection{2D shallow water PDE system}
The corresponding 2D SWE system (depth-averaged over the third dimension $z$) is given by \cite{DELIS2005,wei2006well,dias2007dynamics,dutykh2007mathematical}
\begin{equation}
    % \begin{cases}
        % h_t +(hu)_x+(hv)_y=0, \\ u_t + uu_x+vu_y+g\eta_x=0, \\ v_t + uv_x+vv_y+g\eta_y=0, \\ h=h_0+\eta-\xi,
    % \end{cases} \text{ or, }  
    \begin{cases}
        \eta_t -\xi_t +(hu)_x+(hv)_y=0, \\ u_t + uu_x+vu_y+g\eta_x=0, \\ v_t + uv_x+vv_y+g\eta_y=0, \\ h=h_0+\eta-\xi,
    \end{cases}
\end{equation}
where $u=u(x,y,t), h_0 = h_0(x,y), \eta = \eta(x,y,t)$ and $\xi = \xi(x,y,t)$ for $(x,y)\in \Omega$, $\Omega \subset \mathbb{R}^2$ being the considered spatial domain, and $t\in[0, T]$. 
% \begin{equation}
%     \begin{cases}
%         \eta_t -\xi_t +(hu)_x+(hv)_y=0 \\ u_t + uu_x+vu_y+g\eta_x=0 \\ v_t + uv_x+vv_y+g\eta_y=0 \\ h=h_0+\eta-\xi
%     \end{cases}
% \end{equation}
In vectorial form, defining $\mathbf{U}=( \eta~u~v)^\text{T}$, this can be written, similarly to the 1D case above, as
\begin{equation} \label{2D SWE}
    \mathbf{U}_t+\mathbf{F}(\mathbf{U},\mathbf{\nabla U})=\mathbf{S},
\end{equation}
where
\begin{equation*}
  \begin{aligned}
    \mathbf{F}(\mathbf{U},\mathbf{\nabla U}) 
    % & = \begin{pmatrix}h_x u+h u_x + h_y v+h v_y\\u u_x + v u_y +g \eta_x \\ u v_x + v v_y + g \eta_y \end{pmatrix}\\
      & = \begin{pmatrix} (\eta_x+h_{0,x}-\xi_x)u+(\eta+h_0-\xi) u_x + (\eta_y+h_{0,y}-\xi_y) v+(\eta+h_0-\xi) v_y\\u u_x + v u_y +g \eta_x \\ u v_x + v v_y + g \eta_y \end{pmatrix}
  \end{aligned}\text{and}~~\mathbf{S}=\begin{pmatrix}\xi_t\\0\\0\end{pmatrix}.
\end{equation*}
% and
% \begin{equation*}
%     \mathbf{S}=\begin{pmatrix}\xi_t\\0\\0\end{pmatrix}.
% \end{equation*}
% Similarly to the 1D system, the well-posedness of the 2D system is completed by boundary conditions.
The corresponding 2D wall boundary conditions are given by \cite{bristeau2001boundary}
\begin{equation} \label{eq:2Dwallbc}
    \begin{cases}
        \left( \mathbf{\nabla} \eta \right). \mathbf{n}=0, \\
        \mathbf{u}. \mathbf{n} =0,
    \end{cases} \quad (x,y)\in{\partial \Omega}, t\in [0,T],
\end{equation}
where $\mathbf{u}=(u \ v)^\text{T}$ and $\mathbf{n}$ is the normal outward to the boundary. The corresponding non-reflective (radiation) boundary conditions, assuming a propagation normal to the boundaries, are given by \cite{Durran1999}
\begin{equation} \label{eq:2Dradiationbc}
    \begin{cases}
        \eta_t+c_n \left(\mathbf{\nabla} \eta \right) . \mathbf{n}=0, \\
        u_t + c_n \left(\mathbf{\nabla} u \right) . \mathbf{n}=0, \\
        v_t + c_n \left(\mathbf{\nabla} v \right) . \mathbf{n}=0, \\
    \end{cases} \quad (x,y)\in{\partial \Omega},~~t\in [0,T],
\end{equation}
where $c_n=\pm \sqrt{gh}$ is the propagation speed directed outward to the boundary (positive sign for right end, i.e. $x=x_\text{max}$ or $y=y_\text{max}$ boundaries).

\section{A high-order SWE solver based on Fourier continuation (FC)} \label{Methodology}
This section introduces a new, high-order, pseudo-spectral solver for treating both the 1D and 2D SWE systems described above. The methodology consists of employing a high-order Fourier continuation (FC) technique \cite{amlani2016,bruno2010} for trigonometric interpolation in space (and calculation of the subsequent derivatives in Equations~\eqref{1D SWE} and~\eqref{2D SWE}), together with an explicit high-order time integration scheme. The general principles of FC are presented in \autoref{FC}, and an accelerated version of this method is briefly described in \autoref{Accelerated FC}. The complete solver is summarized in \autoref{SWE solver}.

\subsection{Spatial discretization}  \label{FC}

\begin{figure}[!t]
	\centering
	\includegraphics[width=0.4\textwidth]{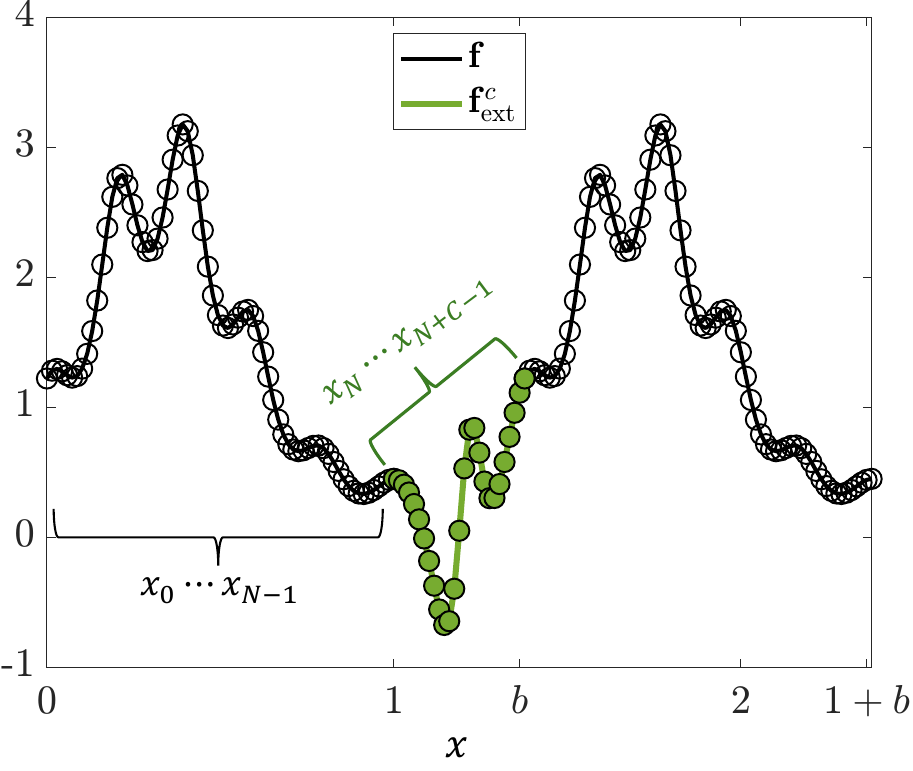}
	\caption{A representative illustration of a periodic extension produced by Fourier continuation.}
	\label{fig:FC}
\end{figure}

Let $N \in \mathbb{N}^*$ represent the number of points discretizing a smooth (possibly non-periodic) function $f:[0,1] \rightarrow \mathbb{R}$. Here, the unit interval is considered without loss of generality (but can be applied to any interval by simple affine transformation, as is the case for all results of this paper). Defining $\{x_j\}$ as the grid points corresponding to a uniform discretization of the domain $[0,1]$, i.e., $x_j=j \Delta x,~0 \leq j \leq N-1$ where $\Delta x =1/(N-1)$, the FC method consists of extending $f$ into a function $f^c$ that is periodic on a slightly larger  interval $[0,b],~b>1$ and that is given by a band-limited trigonometric (Fourier) series as
\begin{equation}\label{eq:FCcontinuous}
    f^c:
    \begin{cases}
        [0,b] \rightarrow \mathbb{C}\\x \mapsto \displaystyle\sum\limits_{k=-M}^{M}a_k \exp\left( \frac{2\pi i k x}{b}\right),
    \end{cases}
\end{equation}
where $M=(N+C)/2$ (for a number of $C$ extension points that define the interval $b$, i.e., $b=(N+C-1)\Delta x$). The coefficients $\{a_k\},~k=-M,\dots,M$ are determined by minimizing the discrete error between $f^c$ and the original function $f$ evaluated at the grid points $\mathbf{x} = \left(x_0~\dots~x_{N-1}\right)^\text{T}$, i.e., 
\begin{equation}\label{eq:L2}
    \{a_k\} =\arg \min \left\lVert f^c(\mathbf{x})-f(\mathbf{x}) \right\rVert_2.
\end{equation}
In other words, $f^c$ approximates a trigonometric interpolation of $f$ on its discrete points.  \autoref{fig:FC} illustrates the corresponding extension that is achieved for an original discretized function $\mathbf{f} = \left(f\left(x_0\right)~\dots~f\left(x_{N-1}\right)\right)^\text{T}$, where $\mathbf{f}^c_\text{ext}$ represents the $C$ interpolated values in $[1,b]$. The complete discrete continuation function $\mathbf{f}^c$ is given discretely by 
$$\mathbf{f}^c =\left(\mathbf{f}~\mathbf{f}^c_{\text{ext}}\right)^\text{T}=  \left(f\left(x_0\right),\dots,f\left(x_{N-1}\right),f^c_\text{ext}\left(x_N\right),\dots,f^c_\text{ext}\left(x_{N+C-1}\right)\right)^\text{T}.$$
% $$\mathbf{f}^c =\begin{pmatrix}\mathbf{f}\\ \mathbf{f}^c_{\text{ext}}\end{pmatrix}= \begin{pmatrix}f\left(x_0\right) \\ \vdots \\ f\left(x_{N-1}\right) \\ f^c_\text{ext}\left(x_N\right)\\ \vdots \\ f^c_\text{ext}\left(x_{N+C-1}\right)\end{pmatrix}.$$
Once the coefficients of the Fourier series expansion $f^c$ in \autoref{eq:FCcontinuous} have been determined, its corresponding spatial derivatives can be obtained analytically by term-wise differentiation, i.e., 
\begin{equation}\label{eq:termwise}
    (f^c)'(x)=\sum\limits_{k=-M}^{M}\frac{2\pi i k}{b}a_k \exp \left( \frac{2\pi i k x}{b} \right).
\end{equation}
Thus, one obtains an approximation of the derivative of $f$ by restricting $f^c$ to the original interval, i.e.,
\begin{equation}
    f'(x) \approx \left(f^c \right)'(x) \quad \forall x \in [0,1].
\end{equation}
% The accuracy of spectral methods are dependent on the smoothness of the function $f$. 
As in other spectral approaches, if $f$ is infinitely differentiable ($\mathrm{C}^\infty$), the truncation error of such an approximation converges to zero faster than any finite power of $\Delta x$ \cite{Durran1999}.

In practice, for a time-domain PDE solver, such extensions needs to be constructed at each timestep ($\{a_k\}$ changes at each timestep). That is, \autoref{eq:L2} must be solved, for example using a Singular Value Decomposition (SVD), at each timestep. The corresponding computational complexity of such minimizations is well-known to be $\mathcal{O} \left(N^3 \right)$, which is too computationally expensive for dynamic problems where the SVD must be applied at each timestep. For this, an accelerated method known as FC(Gram) is employed in the solver proposed in this work (detailed in the following section).

% It can be slightly accelerated by calculating $(f^c)'(x_k)$ using the Fast Fourier Transform (FFT) \cite{FFT}.

\subsection{Accelerated Fourier continuation: FC(Gram)}   \label{Accelerated FC}
\begin{figure*}[!t]
    \centering
    \begin{subfigure}[t]{0.4\textwidth}
        \centering
        \includegraphics[width=\textwidth]{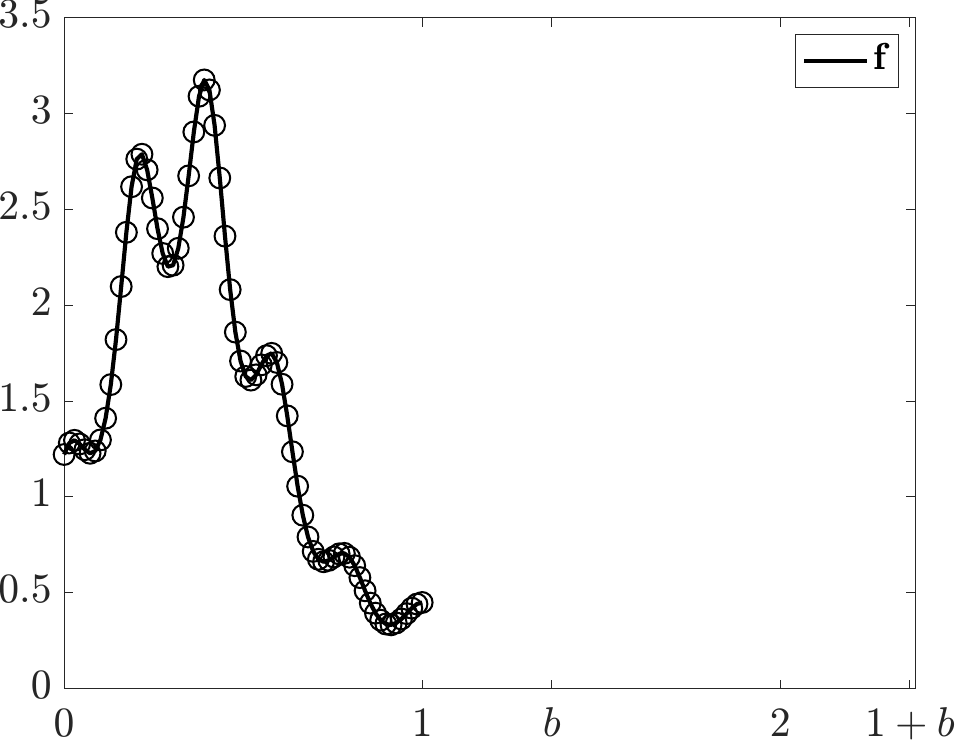}
        \caption[]%
        {{\footnotesize Discretization of a non-periodic function $f$. \\}}
        \label{function f}
    \end{subfigure}\quad\quad
    \begin{subfigure}[t]{0.4\textwidth}
        \centering
        \includegraphics[width=\textwidth]{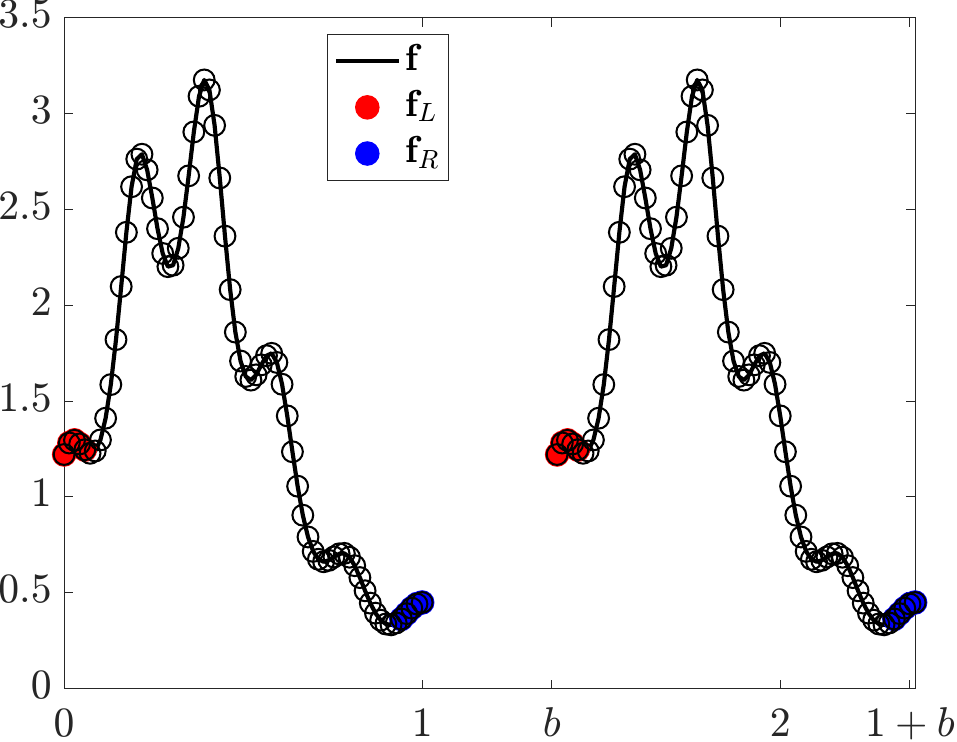}
        \caption[]%
        {{\footnotesize Translation by a distance $1+b$, where red and blue denote the left and right matching points, respectively. }}
        % left (red) and right (blue) matching points. $f$ is replicated for the continuation.}}
        \label{cont points}
    \end{subfigure}
    
    \medskip
    
    % \begin{subfigure}[t]{0.45\textwidth}
    %     \centering
    %     \includegraphics[width=\textwidth]{FC2.pdf}
    %     \caption[]%
    %     {{\footnotesize Continuation of the left matching points $\{x_0,\dots,x_{L-1} \}$ to zero.}}
    %     \label{left cont}
    % \end{subfigure}\hfill
    \begin{subfigure}[t]{0.4\textwidth}
        \centering 
        \includegraphics[width=\textwidth]{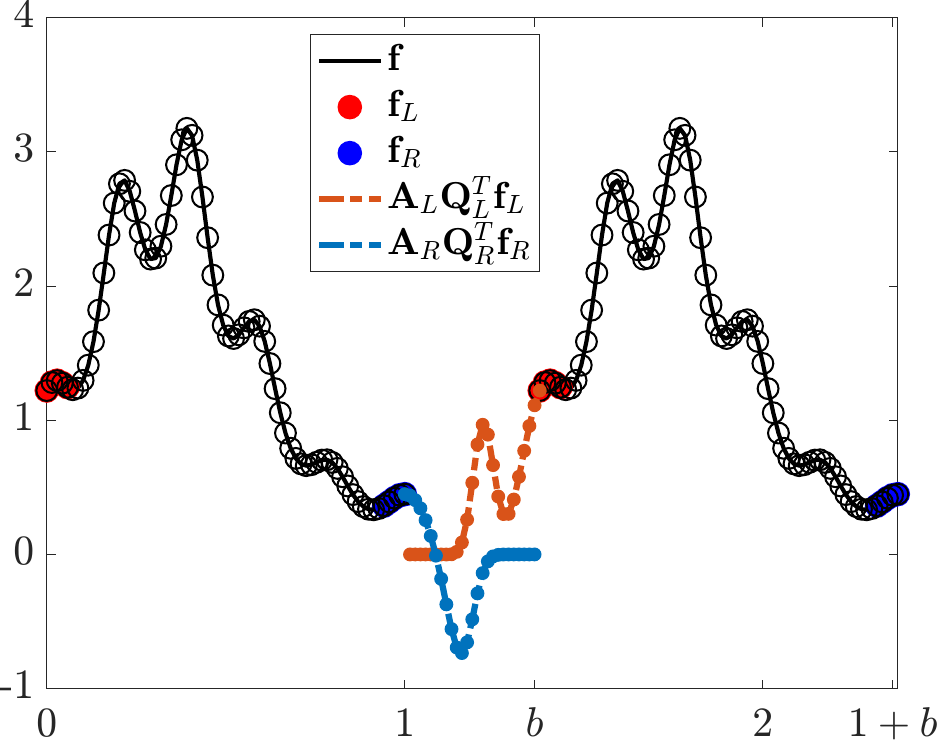}
        \caption[]%
        {{\footnotesize Continuation to zero of the left matching points $\{x_0,\dots,x_{L-1} \}$ (red) and right matching points $\{x_{N-R},\dots,x_{N-1} \}$ (blue).}}
        \label{right cont}
    \end{subfigure}\vspace{0.03\textwidth}\quad\quad
    \begin{subfigure}[t]{0.4\textwidth}
        \centering 
        \includegraphics[width=\textwidth]{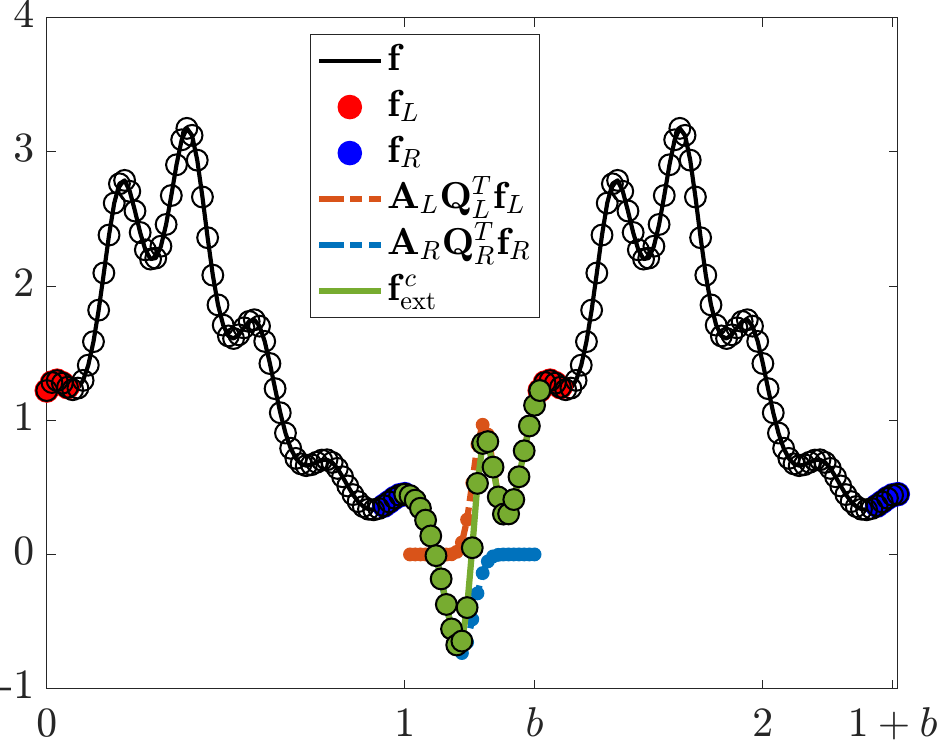}
        \caption[]%
        {{\footnotesize Summation of the two (leftward and rightward) continuations to produce a periodic function $f^c$, represented by the values $
    \mathbf{f}^c=\left(\mathbf{f}~\mathbf{f}^c_\text{ext}\right)^\text{T}$.}}
        \label{global cont}
    \end{subfigure}
    
    \vspace*{-5mm}
    
    \caption{A summary of the accelerated FC(Gram) continuation methodology.}
    \label{fig:Accelerated continuation}
\end{figure*}

An accelerated version of FC has been developed in order to pre-calculate extensions to enable its use with reasonable computational cost while maintaining its interesting properties (high order of accuracy, minimal numerical dispersion) \cite{amlani2016,bruno2010,AmlaniPahlevan2020,albin2011}. This accelerated version is referred as FC(Gram) because it involves Gram polynomials on which a continuation \emph{basis} is precomputed and on which a handful of left and right endpoint function values are projected to form the complete continuation. A visualization of the accelerated continuation procedure, which is detailed in what follows, is presented in \autoref{fig:Accelerated continuation}.

Define two subsets of points on the left and right ends of $\{x_0,\dots,x_{N-1}\}$ (see \autoref{cont points}): $\{x_0,\dots,x_{L-1}\}$ (where $L \in \mathbb{N}^*$ is the number of left points, $L \leq N$) and $\{x_{N-R},\dots,x_{N-1}\}$ (where $R \in \mathbb{N}^*$ is the number of right points). These points are referred to as \textit{matching points} in this work because each of these two subsets of points will be interpolated by a polynomial and used to continue the original discrete function. The continuation is done separately for each sides, enabling more flexibility in the choice of number of matching points as well as particular boundary conditions. 

On the left side, there exists a polynomial of degree $L-1$ that considers discrete function values $\{f \left(x_0 \right),\dots,f \left(x_{L-1} \right)\}$ at the corresponding matching points $\{x_0,\dots,x_{L-1}\}$ (additionally, this polynomial is unique). Finding the corresponding coefficients $c_0, \dots, c_{L-1}$ of an interpolating polynomial is given by the system
\begin{equation} \label{matching system}
    \begin{cases}
    c_0+c_1x_0+c_2 \left(x_0 \right)^2\dots+c_{L-1} \left(x_0 \right)^{L-1}=f \left(x_0 \right), \\ \vdots \\ c_0+c_1x_{L-1}+c_2 \left(x_{L-1} \right)^2+\dots+c_{L-1} \left(x_{L-1} \right)^{L-1}=f \left(x_{L-1} \right),
    \end{cases}
\end{equation}
which can equivalently be written in matrix form as
\begin{equation*}
    \mathbf{P}_L \mathbf{c}_L=\mathbf{f}_L,
\end{equation*}
where
\begin{equation}
    \mathbf{P}_L = \begin{pmatrix}
    1 & x_0 & \left(x_0 \right)^2 & \cdots & \left(x_0 \right)^{L-1} \\ \vdots & \vdots & \vdots & & \vdots \\ 1 & x_{L-1} & \left(x_{L-1} \right)^2 & \cdots & \left(x_{L-1} \right)^{L-1}
    \end{pmatrix}, \quad \mathbf{c}_L = 
    \begin{pmatrix}
    c_0 \\ \vdots \\ c_{L-1}
    \end{pmatrix}, \quad \text{and} \quad \mathbf{f}_L =    
    \begin{pmatrix}
    f \left(x_0 \right) \\  \vdots \\ f \left(x_{L-1} \right)
    \end{pmatrix}.
\end{equation}
% and
% \begin{equation*}
%     \mathbf{f}_L =    
%     \begin{pmatrix}
%     f \left(x_0 \right) \\  \vdots \\ f \left(x_{L-1} \right)
%     \end{pmatrix}.
% \end{equation*}
By employing a modified Gram-Schmidt orthonormalization, the real-valued square matrix $\mathbf{P}_L$ (also known as a \emph{Vandermonde matrix}) can be decomposed as $\mathbf{P}_L=\mathbf{Q}_L \mathbf{R}_L$, where $\mathbf{Q}_L$ is an orthogonal matrix ($\mathbf{Q}_L \mathbf{Q}_L^\text{T}=\mathbf{I}_{L}$, the identity matrix of order $L$) and $\mathbf{R}_L$ is an upper triangular matrix. The columns of $\mathbf{Q}_L$ form an orthonormal basis of polynomials. The strategy of FC(Gram) consists in performing the continuation on this basis of polynomials (each column of the corresponding $\mathbf{Q}_L$), before projecting the matching points onto such basis. This enables one to invoke the SVD-based procedure that form continuations independently of the specific values of the matching points $\mathbf{f}_L$  (which may evolve in time). That is, the extensions are pre-computed as a basis \emph{a priori}. A similar procedure can then be performed for the rightward matching points to construct another basis of polynomials forming the columns of a matrix $\mathbf{Q}_R$.

With the orthonormalized Gram basis in hand, the first step of the FC(Gram) procedure is to build operators $\mathbf{A}_L$ and $\mathbf{A}_R$ that continue $\mathbf{Q}_L$ and $\mathbf{Q}_R$ to zero (i.e., leftwardly and rightwardly), so that their sum gives a global continuation that smoothly connects the left and right sides of $f$. For example, to continue $\mathbf{Q}_L$ to zero, one takes its column vectors $\mathbf{q}^j_L$, which contain the $L$ values $q_{j,L} \left(x_0 \right), q_{j,L} \left(x_1 \right), \dots, q_{j,L} \left(x_{L-1} \right)$, and searches for an auxiliary interpolating trigonometric polynomial of the form given by~\cite{amlani2016}
\begin{equation}
    q_{j,L}^c:
    \begin{cases}
        \left[0,\left(L+2C+Z \right)\Delta x \right) \rightarrow \mathbb{C}\\x \mapsto \displaystyle\sum\limits_{k=-M}^{M}a_{k,j} \exp \left(\frac{2\pi i k x}{\left(L+2C+Z-1 \right)\Delta x} \right),
    \end{cases}
\end{equation}
%where $q_{j,L}^c$ is $\left(\left( L+2C+Z \right)\Delta x\right)$-periodic and approximates $q_{j,L}$ on $\left[0,(L-1)\Delta x \right]$, while smoothly blending to zero in $\left[L\Delta x, \left( L+C \right)\Delta x \right)$. 
where $Z$ is a number of points equal to zero to match for blending and where $\{a_{k,j} \}$ are the unknown coefficients. The role of $q_{j,L}^c$ is to:
\begin{itemize}
    \item approximate the $q_{j,L}$ values in $\left[0, \left( L-1 \right)\Delta x \right]$;
    \item smoothly continue (or ``blend") over $C$ points the $q_{j,L}$ values to zero in $\left[L\Delta x, \left( L+C-1 \right)\Delta x \right]$;
    \item smoothly approximate the $Z$ zero-valued points in $\left[\left(L+C \right)\Delta x, \left(L+C+Z-1 \right)\Delta x \right]$ that force the continuation to zero;
    \item smoothly continue the zero values over a certain number of extra points (here, taken to be $C$) in $\left[\left(L+C+Z \right)\Delta x, \left( L+2C+Z-1 \right)\Delta x \right]$ to obtain periodicity over the auxiliary interval.
\end{itemize}
%Discretely, for the points $x_m = m\Delta x, $
%It is possible to write the values of $q_{j,L}^c$ on the discrete points as
%\begin{equation}
%    q_{j,L}^c \left(x_m \right)=\sum\limits_{k=-M}^{M}a_{k,j} \exp \left(\frac{2\pi i k m}{L+2C+Z-1} \right)
%\end{equation}
The coefficients $a_{k,j}$ satisfying the above can be found by solving an $L^2$ minimization problem, i.e., 
\begin{equation}\label{eq:L2aux}
    (a_{k,j})_k=\arg \min \left\lVert q_{j,L}^c(\mathbf{x}_L')-\begin{pmatrix}q_{j,L}(\mathbf{x}_L) \\ \mathbf{0}_{Z\times 1} \end{pmatrix} \right\rVert_2,
\end{equation}
where $\mathbf{x}_L=\left( x_0~\dots~x_{L-1} \right)^\text{T}$, $\mathbf{x}_L'=\left( x_{0}~\dots~x_{L+Z-1} \right)^\text{T}$ and $\mathbf{0}_{Z\times 1}$ is a zero vector of length $Z$.
%\begin{equation}
%    \left( a_{k,j} \right)_k=\arg \min \left[\sum\limits_{m=0}^{L-1} \left(q_{j,L}^c \left(x_m \right)-q_{j,L} \left(x_m \right) \right)^2 + \sum\limits_{m=L}^{L+C-1} q_{j,L}^c \left( x_m \right)^2 \right]
%\end{equation}
The solution defined by \autoref{eq:L2aux} can be obtained for each of the $j$ base polynomials $q_{j,L}$ via a truncated SVD. The resulting  periodic continuation operator for $\mathbf{Q}_L$ formed by the $\left( q_{j,L}^c \right)_j$ is hence given by
\begin{equation}
    \left( \mathbf{A}_L \right)_{m,j}= q_{j,L}^c \left(x_{m-1} \right), \quad m \in \{1, \dots, L+C\}, \quad j \in \{1, \dots, L\}.
\end{equation}
Such a continuation operator is clearly pre-computed since the minimization problem does not depend on the specific function values of $f$ (merely the interpolating basis). In a similar manner to all of the above, one can straightforwardly construct an operator $\mathbf{A}_R$ that continues $\mathbf{f}_R$ to zero rightwardly with $R$ matching points.

\begin{remark}
For equal left and right matching points (i.e., $L=R$), the corresponding FC operators $\mathbf{Q_L}, \mathbf{Q_R}$ (respectively $\mathbf{A}_L, \mathbf{A}_R$) differ only by column-wise ordering (respectively row-wise). Hence the precomputation of such operators needs only to be conducted on one side.
\end{remark}

The second step of the FC(Gram) procedure is to project the left and right endpoints $\mathbf{f}_L$ and $\mathbf{f}_R$ by
\begin{equation}
\begin{cases}
    \text{proj} \left(\mathbf{f}_L \right)=\mathbf{Q}_L^\text{T} \mathbf{f}_L, \\
    \text{proj} \left(\mathbf{f}_R \right)=\mathbf{Q}_R^\text{T} \mathbf{f}_R,
\end{cases}
\end{equation}
to get the corresponding coefficients in the continuation basis $\mathbf{A}_L$ and $\mathbf{A}_R$, respectively. Hence, the quantity $\mathbf{A}_L \mathbf{Q}_L^\text{T} \mathbf{f}_L$ represents the leftward continuation of the left values $\mathbf{f}_L$ towards zero and, similarly, the quantity $\mathbf{A}_R \mathbf{Q}_R^\text{T} \mathbf{f}_R$ represents the rightward continuation of the right values $\mathbf{f}_R$ towards zero (see \autoref{right cont}). In the end, the sum $\mathbf{A}_L \mathbf{Q}_L^\text{T} \mathbf{f}_L + \mathbf{A}_R \mathbf{Q}_R^\text{T} \mathbf{f}_R$ represents the complete continuation that matches the left end of the original discrete function $\mathbf{f}$ with the right over an interval of $C$ points (\autoref{global cont}). Hence the final discrete continued function $\mathbf{f}^c$ is given by
\begin{equation}
    \mathbf{f}^c=
    \begin{pmatrix}
    \mathbf{f} \\ \mathbf{A}_L \mathbf{Q}_L^\text{T} \mathbf{f}_L + \mathbf{A}_R \mathbf{Q}_R^\text{T} \mathbf{f}_R
    \end{pmatrix},
\end{equation}
which matches $f$ discretely on the original interval $[0,1]$, but is periodic on a slightly larger interval $[0,b]$ (discretely, $b=(N+C-1)\Delta x$).
Hence the vector $\mathbf{f}^c$ represents discrete values of a corresponding continuous periodic function $f^c$.

Thanks to this acquired periodicity, one can now approximate $f^c$ by the truncated Fourier series
\begin{equation*}
    x \mapsto \sum\limits_{k=-M}^{M}\alpha_k \exp \left(\frac{2\pi i k x}{b} \right),
\end{equation*}
where the coefficients $\{ \alpha_k\}$ can be obtained directly using the Fast Fourier Transform (FFT) \cite{FFT} of the function's corresponding discretization $\mathbf{f}^c$, i.e., $\mathbf{\hat{f}}^c=\text{fft} \left(\mathbf{f}^c \right) = \left(\alpha_{-M/2}~\dots~\alpha_{M/2} \right)^\text{T}$. 
The discrete derivative can finally be obtained by forming the coefficients of the term-wise differentiation of the series (see \autoref{eq:termwise}) and performing the inverse (Fast) Fourier transform given by
\begin{equation}
    \mathbf{f}' \approx \left( \mathbf{f}^c \right)' \approx \text{ifft} \left( \mathbf{K}\mathbf{\hat{f}}^c \right),
\end{equation}
where $\mathbf{K}$ is a matrix containing identical column vectors whose $k$-th component is ${2\pi i \alpha_k k}/{b}$. The overall algorithm described above incurs significantly less computational cost than the simplest FC described in \autoref{FC} which requires an SVD each time an extension needs to be produced. Indeed, the $\mathcal{O} \left(N^3 \right)$ continuation is only done once for the Gram polynomials in the accelerated version, after which it is stored. Once these pre-computed operators have been loaded, the resulting PDE resolution in space is dominated by the FFT, whose complexity is $\mathcal{O}(N \log N)$ \cite{FFT}.

\begin{remark}
The complete FC(Gram) algorithm employed in the present solver uses a number of extra zero points for the continuation, as well as an oversampling of values in order to mitigate the numerical ill-conditioning well-known for Vandermonde matrices. Additionally, the overall precomputation procedure uses symbolic algebra for the QR decompositions and high-precision arithmetic for the SVDs (256 digits instead of 16 for a double). For simplicity and clarify of the above presentations, these details are left to the reader  \cite{amlani2016}. 
\end{remark}

\subsubsection{Extension to higher spatial dimensions}
\begin{figure}[!t]
	\centering
	\includegraphics[width=0.49\textwidth]{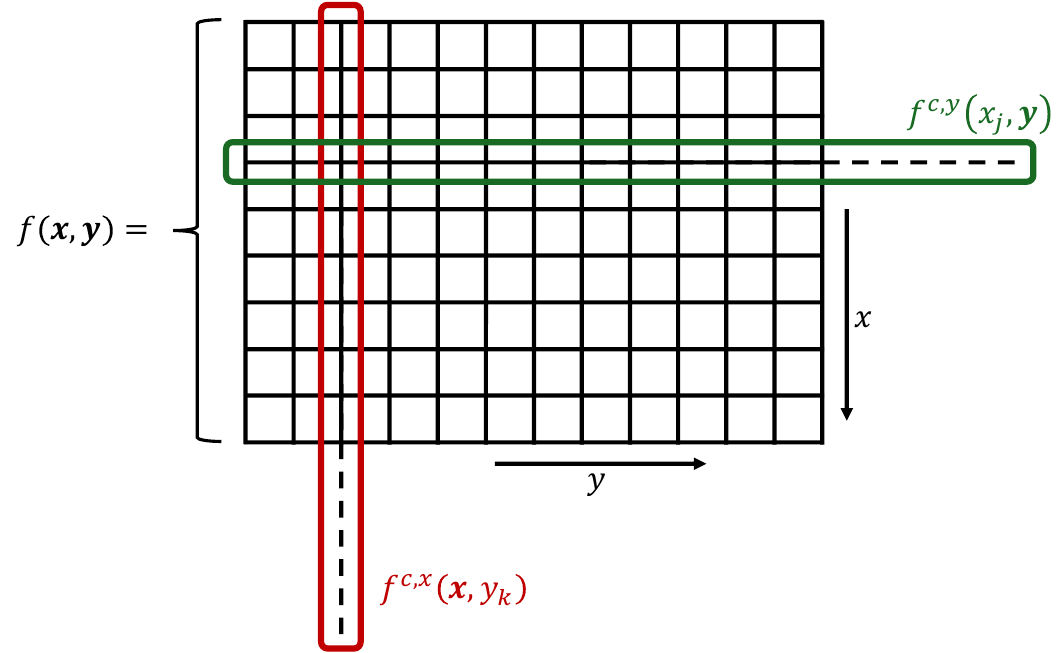}
	\caption{A representative illustration of a 2D "slice-wise" Fourier continuation.}
	\label{fig:2D FC}
\end{figure}

Consider a discretization of $\left(N_x, N_y \right) \in \left(\mathbb{N}^* \right)^2$ points for a smooth function $f:[0,1]^2 \rightarrow \mathbb{R}$, i.e.,  $x_j=j \Delta x, \ 0 \leq j \leq \left(N_x-1 \right)$ for $\Delta x =1/ \left(N_x-1 \right)$ and $y_k=k \Delta y, \ 0 \leq k \leq \left(N_y-1 \right)$ for $\Delta y =1/ \left(N_y-1 \right)$.
The accelerated FC(Gram) can be extended to 2D problem line-wise. The discrete values of $f$, now represented by a matrix $\textbf{f} \in \mathbb{R}^{N_x \times N_y}$, can be considered as a set of 1D slices on which to apply the FC procedure. Thus, two extensions are obtained: $f^{c,x}$, representing the continued function with respect to $x$ for each $y$-slice, and $f^{c,y}$, representing the continued function with respect to $y$ for each $x$-slice. These functions can be discretely represented at the points $(x_j,y_k)$ by the matrices $\mathbf{f}^{c,x}$ and $\mathbf{f}^{c,y}$, respectively. This procedure is not a proper 2D continuation, but a "2D slice-wise" continuation. The \autoref{fig:2D FC} illustrates this method of application of the 1D FC(Gram) procedure over the 2D grid.
Once the continuations are properly constructed, one can differentiate slice by slice the $x$ and $y$ continuations, just as before:
\begin{equation}
    \frac{\partial f}{\partial x} \approx \frac{\partial f^{c,x}}{\partial x} \quad \text{and} \quad   \frac{\partial f}{\partial y} \approx \frac{\partial f^{c,y}}{\partial y}.
\end{equation}

% A possible lead to accelerate this 2D slice-wise FC could be to work in the context of two-dimensional periodicity. In this regard, recent work has been done to implement a proper 2D Fourier Continuation \cite{Bruno2022}.
\begin{remark}
In the rest of this manuscript, and particularly in the results, `FC' refers to the accelerated FC(Gram) methodology.
\end{remark}

\subsection{Time integration}   \label{SWE solver}
The full shallow water PDE solver developed herein exploits the precision of FC in its accelerated version. The FC(Gram) method that enables high-order computation of the spatial derivatives, is completed with a time integration scheme in order to build a complete PDE solver. Classical integration methods are applicable, including the forward Euler method. However, to remain consistent with the high-order and dispersionless properties of FC, which are needed to model earthquake-generated tsunamis, a suitably high-order time scheme is employed in the solver introduced in this work. Although the well-known fourth-order Runge-Kutta (RK4) scheme is a common choice, in the interest of computational speed, a more convenient choice is fourth-order Adams-Bashforth (AB4), whose regions of stability are still adequate \cite{Demailly2016}, including for FC-based methods~\cite{amlani2016,amlaniPOF}. An obvious advantage of using the AB4 of \autoref{eq:AB4} instead of the RK4 of \autoref{eq:RK4} is that the latter can be slow in the current multi-dimensional PDE context since it relies on calculating three intermediate time stages for each timestep update, leading to four evaluations per timestep of the right-hand side (whose cost is dominated by the FC procedure). On the other hand, AB4 requires only one calculation of the right-hand side per timestep (at the lower cost of stocking the previous three right-hand sides). Additionally, at the same timestep, AB4 and RK4 have been shown to have equivalent accuracy in other FC-based solvers~\cite{Amlani2023}.

% Even if the CFL condition for RK4 is milder than AB4 (due to its larger stability domain), the trade-off between the two methods is in favor of AB4 to ensure an order five of spatial convergence.
After spatial discretization and differentiation, the corresponding ordinary differential equations (ODEs) in time are of the form
\begin{equation}
    \mathbf{U}'(t)=\textbf{RHS}\left(t,\mathbf{U}(t)\right),
\end{equation}
where $\mathbf{U}(t)= (\boldsymbol{\eta}(t)~\mathbf{u}(t)~
\mathbf{v}(t))^\text{T}$, $\boldsymbol{\eta}(t)=\{\eta(x_j,y_k,t)\}$, and the corresponding right-hand side (\autoref{2D SWE}) given by
\begin{equation} \label{RHS}
    \textbf{RHS}\left(t,\mathbf{U}(t)\right)=-\mathbf{F}(\mathbf{U}(t),\mathbf{\nabla U}(t))+\mathbf{S}(t).
\end{equation}
Employing AB4 \cite{Demailly2016}, the resulting discrete formulation for evolving to a timestep $n+1$ (corresponding to time $t^{n+1}$) is given by 
\begin{equation}\label{eq:AB4}
    \mathbf{U}^{n+1}=\mathbf{U}^n+\frac{\Delta t}{24} \left[55\mathbf{RHS} \left(t^n,\mathbf{U}^n \right)-59\mathbf{RHS} \left(t^{n-1},\mathbf{U}^{n-1} \right)+37\mathbf{RHS} \left(t^{n-2},\mathbf{U}^{n-2} \right)-9\mathbf{RHS} \left(t^{n-3},\mathbf{U}^{n-3} \right) \right].
\end{equation}
The computation of $\textbf{RHS}$ involves determining the spatial derivatives, which is accomplished  using the 2D-slicewise FC(Gram) method described in \autoref{Accelerated FC}. The corresponding procedure for computing $\textbf{RHS}$ is summarized in Algorithm \ref{alg:RHS}.

\RestyleAlgo{ruled}
\begin{algorithm}[!h]
\small
\DontPrintSemicolon
\caption{Computation of the right-hand side for time integration.}\label{alg:RHS}
\SetKwFunction{FMain}{RHS}
\KwInput{vector of discrete unknowns $\mathbf{U}$, the current timestep $t$, and a given dynamic sea floor motion $\boldsymbol{\xi}$ with derivatives $\boldsymbol{\xi}_x, \boldsymbol{\xi}_y, \boldsymbol{\xi}_t$}
\KwOutput{right-hand side $\mathbf{RHS}$}
\SetKwProg{Pn}{Function}{:}{\KwRet}
  \Pn{\FMain{$t,\mathbf{U}, \boldsymbol{\xi}, \boldsymbol{\xi}_x, \boldsymbol{\xi}_y, \boldsymbol{\xi}_t$}}{
        $\mathbf{U}_x \gets \text{FC\_deriv\_x}(\mathbf{U})$ \tcp*{$\mathbf{U}_x= (\boldsymbol{\eta_x}~\mathbf{u}_x~\mathbf{v}_x )^\text{T}$ via \autoref{Accelerated FC}}
        $\mathbf{U}_y \gets \text{FC\_deriv\_y}(\mathbf{U})$ \tcp*{$\mathbf{U}_y= (\boldsymbol{\eta_y}~\mathbf{u}_y~\mathbf{v}_y )^\text{T}$ via \autoref{Accelerated FC}}
        $\mathbf{h} \gets \mathbf{h}_{0}+\boldsymbol{\eta}-\boldsymbol{\xi}$\;
        $\mathbf{h}_x \gets \mathbf{h}_{0x}+\boldsymbol{\eta}_x-\boldsymbol{\xi}_x$\;
        $\mathbf{h}_y \gets \mathbf{h}_{0y}+\boldsymbol{\eta}_y-\boldsymbol{\xi}_y$\;
        $\mathbf{RHS}\gets -\begin{pmatrix} \mathbf{h}_x \mathbf{u}+\mathbf{h} \mathbf{u}_x+\mathbf{h}_y \mathbf{v}+\mathbf{h} \mathbf{v}_y \\ \mathbf{u} \mathbf{u}_x+\mathbf{v} \mathbf{u}_y+g \boldsymbol{\eta}_x \\ \mathbf{u} \mathbf{v}_x+\mathbf{v} \mathbf{v}_y+g \boldsymbol{\eta}_y \end{pmatrix} + \begin{pmatrix} \boldsymbol{\xi}_t \\ \mathbf{0} \\ \mathbf{0} \end{pmatrix}$\;
  }
\end{algorithm}

\begin{remark} Since the use of AB4 requires knowledge of three previous timesteps, the proposed SWE solver is initialized by computing the first three iterations (timesteps) via RK4 via the time evolution expression given by
\begin{equation}\label{eq:RK4}
    \mathbf{U}^{n+1}=\mathbf{U}^n+\frac{\Delta t}{6} \left(\mathbf{k_1}+2\mathbf{k_2}+2\mathbf{k_3}+\mathbf{k_4} \right),
\end{equation}
where
\begin{equation*}
    \begin{cases}
        \mathbf{k_1}=\mathbf{RHS} \left(t^n,\mathbf{U}^n \right), \\
        \mathbf{k_2}=\mathbf{RHS} \left(t^n+\dfrac{\Delta t}{2},\mathbf{U}^n+ \dfrac{\Delta t}{2} \mathbf{k_1}\right), \\
        \mathbf{k_3}=\mathbf{RHS} \left(t^n+\dfrac{\Delta t}{2},\mathbf{U}^n+ \dfrac{\Delta t}{2} \mathbf{k_2}\right), \\
        \mathbf{k_4}=\mathbf{RHS} \left(t^n+ \Delta t,\mathbf{U}^n+ \Delta t \mathbf{k_3}\right).
    \end{cases}
\end{equation*}
Alternatively, one could simply start the simulation after three `zero-state' iterations, which is acceptable for a number of realistic problems that may begin at rest.
\end{remark}

\subsection{Implementation details \& solver summary}\label{remark:FCparams}

The resulting overall solver is summarized in \autoref{alg:SWE solver}, where the selection and application of the boundary conditions are indicated by the function \texttt{ApplyBC}, which considers walls (\autoref{eq:2Dwallbc}), radiation boundary conditions (\autoref{eq:2Dradiationbc}), or prescribed (Dirichlet) values. FC-based solutions in this work are obtained in MATLAB 2023b Update 6 employing a self-developed code. Similarly, the FD-based solvers employed for comparison utilize fully self-implemented code in MATLAB. All results herein are obtained from computations performed on the same MacBook Air M2 2022 with 16GB of RAM.
% \end{remark}

For the simulations and studies that follow, the number of matching points is taken to be $L=R=5$ and the number of continuation points is taken to be $C=25$. Additionally, as is the case for classical spectral methods \cite{hesthaven2008} as well as FC-based methods~\cite{amlani2016,Amlani2023, amlaniPOF}, a filter $\sigma$ is used in the context of time-domain PDEs such that the $n$-th component of the discrete Fourier coefficients is appropriately scaled for stability, i.e.,  ${2\pi \sigma \left(\alpha_n \right) i n}/{b}$, where $\sigma(n)=\exp(-\alpha (n/M)^{2p})$ for $\alpha=16c\Delta t \log \left(10^{-2} \right)/\Delta x$, $p=4$, approximate wave speed $c \approx \sqrt{gH}$ (where $H$ is the mean or still water depth), $M$ as defined in \autoref{FC}, and the frequency index $n$ (such that $-M\leq n\leq M$). Unless otherwise stated, the timestep chosen for all simulations in this work is determined from the CFL condition given by $\Delta t=\text{CFL} \ \Delta x/c$ in 1D and $\Delta t=\text{CFL} \ \text{min}(\Delta x,\Delta y)/c$ in 2D, with corresponding CFL constants given in 1D for the FC solver as well as reference fourth- and sixth-order finite differences (FD4 and FD6, respectively) by $\text{CFL}_{1D}^\text{FC}=0.17$, $\text{CFL}_{1D}^\text{FD4}=0.21$, $\text{CFL}_{1D}^\text{FD6}=0.18$, respectively, and for 2D by  $\text{CFL}_{2D}^\text{FC}=0.1$, $\text{CFL}_{2D}^\text{FD4}=0.17$, $\text{CFL}_{2D}^\text{FD6}=0.14$, respectively. The values of these constants have been empirically determined to be the largest such values that ensure stability, where FC carries slightly smaller (but similar) CFL constants to the FD-based schemes considered for comparison in  \autoref{Numerical study}.

\begin{algorithm}[!h]
\small
\DontPrintSemicolon
\caption{A summary of the FC-based solver for shallow water systems.}\label{alg:SWE solver}
\KwInput{$\boldsymbol{\eta}^1=\eta(\mathbf{x},\mathbf{y},t=0)$, $\mathbf{u}^1=u(\mathbf{x},\mathbf{y},t=0)$, $\mathbf{v}^1=v(\mathbf{x},\mathbf{y},t=0)$, $\Delta t$, $N_t$, $\boldsymbol{\xi}$, $\boldsymbol{\xi}_x$, $\boldsymbol{\xi}_y$, $\boldsymbol{\xi}_t$}
\KwOutput{$\boldsymbol{\eta}^n, \mathbf{u}^n, \mathbf{v}^n, \forall  n \in \{1, \dots, N_t \}$ (complete space-time solutions)} \;
$\mathbf{U}^1 \gets \left( \boldsymbol{\eta}^1, \mathbf{u}^1, \mathbf{v}^1\right)^\text{T}$ \tcp*{Initial values}
\;

\For{$n = 1,\dots,3$\tcp*{Initialize first three timesteps via RK4}}{
$\mathbf{k_1} \gets \texttt{RHS} \left(t^n,\mathbf{U}^n, \boldsymbol{\xi}, \boldsymbol{\xi}_x, \boldsymbol{\xi}_y, \boldsymbol{\xi}_t \right)$\;
$\mathbf{k_2} \gets \texttt{RHS} \left(t^n+\dfrac{\Delta t}{2},\mathbf{U}^n+\dfrac{\Delta t}{2}\mathbf{k_1}, \boldsymbol{\xi}, \boldsymbol{\xi}_x, \boldsymbol{\xi}_y, \boldsymbol{\xi}_t \right)$\;
$\mathbf{k_3} \gets \texttt{RHS} \left(t^n+\dfrac{\Delta t}{2},\mathbf{U}^n+\dfrac{\Delta t}{2}\mathbf{k_2}, \boldsymbol{\xi}, \boldsymbol{\xi}_x, \boldsymbol{\xi}_y, \boldsymbol{\xi}_t \right)$\;
$\mathbf{k_4} \gets \texttt{RHS} \left(t^n+\Delta t,\mathbf{U}^n+\Delta t \mathbf{k_3}, \boldsymbol{\xi}, \boldsymbol{\xi}_x, \boldsymbol{\xi}_y, \boldsymbol{\xi}_t \right)$\;
$\mathbf{U}^{n+1}=\mathbf{U}^n+\dfrac{\Delta t}{6} \left(\mathbf{k_1}+2 \mathbf{k_2}+2 \mathbf{k_3}+\mathbf{k_4}, \boldsymbol{\xi}, \boldsymbol{\xi}_x, \boldsymbol{\xi}_y, \boldsymbol{\xi}_t \right)$\;
$\mathbf{U}^{n+1} \gets \texttt{ApplyBC} \left(t^{n+1},\mathbf{U}^{n+1} \right)$ \tcp*{Boundary conditions}
}
$\mathbf{RHS}^{n-2} \gets \texttt{RHS} \left(t^1,\mathbf{U}^1, \boldsymbol{\xi}, \boldsymbol{\xi}_x, \boldsymbol{\xi}_y, \boldsymbol{\xi}_t \right)$\;
$\mathbf{RHS}^{n-1} \gets \texttt{RHS} \left(t^2,\mathbf{U}^2, \boldsymbol{\xi}, \boldsymbol{\xi}_x, \boldsymbol{\xi}_y, \boldsymbol{\xi}_t \right)$\;
$\mathbf{RHS}^n \gets \texttt{RHS} \left(t^3,\mathbf{U}^3, \boldsymbol{\xi}, \boldsymbol{\xi}_x, \boldsymbol{\xi}_y, \boldsymbol{\xi}_t \right)$\;
\;

\For{$n = 4,\dots,N_t$\tcp*{Main loop (AB4)}}{
$\mathbf{RHS}^{n-3} \gets \mathbf{RHS}^{n-2}$\;
$\mathbf{RHS}^{n-2} \gets \mathbf{RHS}^{n-1}$\;
$\mathbf{RHS}^{n-1} \gets \mathbf{RHS}^n$\;
$\mathbf{RHS}^n \gets \texttt{RHS} \left(t^n,\mathbf{U}^n, \boldsymbol{\xi}, \boldsymbol{\xi}_x, \boldsymbol{\xi}_y, \boldsymbol{\xi}_t \right)$\;
$\mathbf{U}^{n+1}=\mathbf{U}^n+\dfrac{\Delta t}{24} \left(55 \mathbf{RHS}^n-59\mathbf{RHS}^{n-1}+37\mathbf{RHS}^{n-2}-9\mathbf{RHS}^{n-3} \right)$\\
$\mathbf{U}^{n+1} \gets \texttt{ApplyBC} \left(t^{n+1},\mathbf{U}^{n+1} \right)$\;
}
\end{algorithm}

\begin{remark}
The complexity of the 1D solver scales as $\mathcal{O} \left(N_t N \log (N) \right)$ (since the additional $C$ continuation points are fixed), where $N_t$ is the number of timesteps, and $N$ the number of 1D spatial discretization points. In 2D, the corresponding complexity is $\mathcal{O} \left(N_t N_x N_y \log \left(N_{\text{max}} \right) \right)$, where $N_x$, $N_y$ are the number of spatial discretization points in the two directions, and $N_{\text{max}}=\text{max} \left(N_x,N_y \right)$. An important remark on this result is that the ratio of complexity between FC and finite differences is only 
$\log \left(N_{\text{max}} \right)$, i.e., FC is only a factor $\log \left(N_{\text{max}} \right)$ more costly than FD. This is an acceptable compromise for the high-order convergence and numerically dispersionless behavior of an FC solver that is demonstrated in \autoref{Numerical study}, which ultimately results in simulations that are cheaper than FD (to achieve the same error) for realistic settings, especially those including long-distance or long-time propagation.
\end{remark}

\section{Numerical performance of the solver}   \label{Numerical study}
The objective of this section is to verify accuracy, convergence, and code implementation of the FC SWE solver introduced in this work (including via comparisons with standard high-order finite difference methods). The classical method of manufactured solutions \cite{Roache2002,amlani2016,Amlani2023} is employed throughout this section, allowing one to impose any given analytical function as a solution to the PDE of interest by prescribing a non-trivial right-hand-side forcing/source term (determined via substitution of the postulated solution into the PDE) and choosing appropriate initial values and boundary conditions. Such an approach enables the error of a methodology to be determined accurately and its solver implementation to be verified numerically. For the results that follow in this section, the solver detailed in \autoref{Methodology} is assessed against self-implemented finite difference solvers employing exactly the same temporal integration techniques (RK4+AB4), i.e., differing only by spatial differentiation. The FD solvers considered here employ central difference schemes of order four and six, where boundary derivatives are calculated via ghost points~\cite{fedkiw1999}.

% \begin{remark}\label{remark:computer}
% The FC-based solver has been developed in-house using MATLAB 2023b Update 6. Similarly, the FD-based solvers employed for comparison have been fully self-implemented. All results throughout this work are obtained from simulations performed on the same MacBook Air M2 2022 with 16GB of RAM.
% \end{remark}

\subsection{Convergence experiments}\label{sec:1Dconvergence}
\begin{figure*}[!t]
    \centering
    \captionsetup[subfigure]{oneside,margin={0cm,0cm}}
    \begin{subfigure}[b]{0.3\textwidth}
        \centering 
        \includegraphics[width=.8\textwidth]{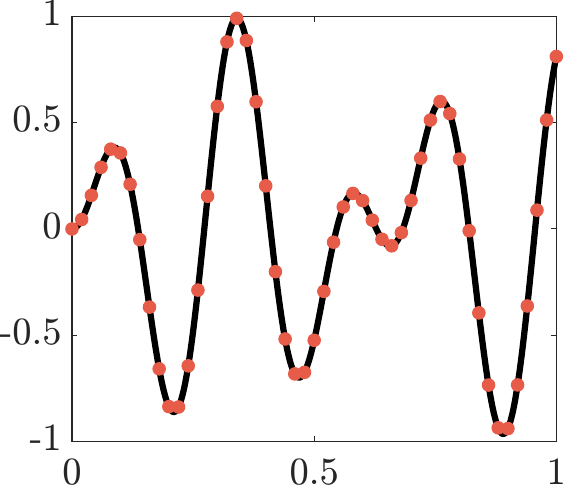}
        \caption[]%
        {{\footnotesize $\eta(x,t=0)$}}
        \label{eta frame 1}
    \end{subfigure}
    \begin{subfigure}[b]{0.3\textwidth}
        \centering
        \includegraphics[width=.8\textwidth]{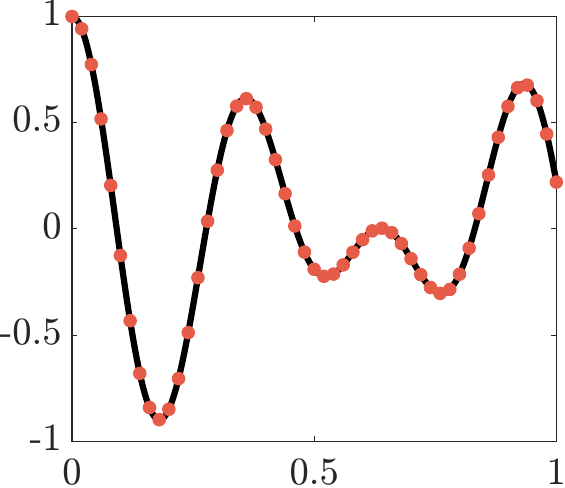}
        \caption[]%
        {{\footnotesize $u(x,t=0)$}}
        \label{u frame 1}
    \end{subfigure}
    \begin{subfigure}[b]{0.3\textwidth}
        \centering 
        \includegraphics[width=.8\textwidth]{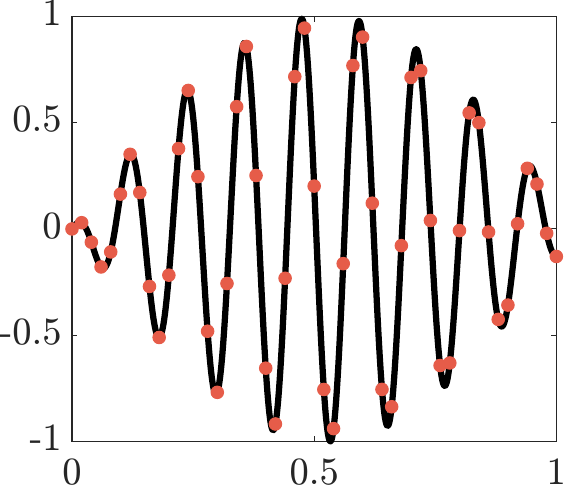}
        \caption[]%
        {{\footnotesize $\xi(x,t=0)$}}
        \label{ksi frame 1}
    \end{subfigure}
    
 \medskip
    \begin{subfigure}[b]{0.3\textwidth}
        \centering 
        \includegraphics[width=.8\textwidth]{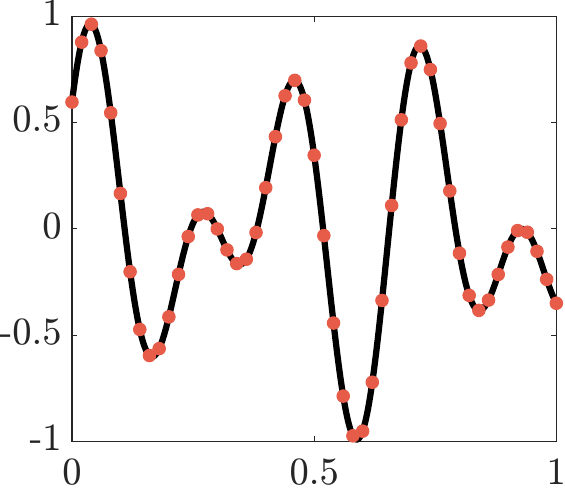}
        \caption[]%
        {{\footnotesize $\eta(x,t=1/2)$}}
        \label{eta frame 2}
    \end{subfigure}
    \begin{subfigure}[b]{0.3\textwidth}
        \centering
        \includegraphics[width=.8\textwidth]{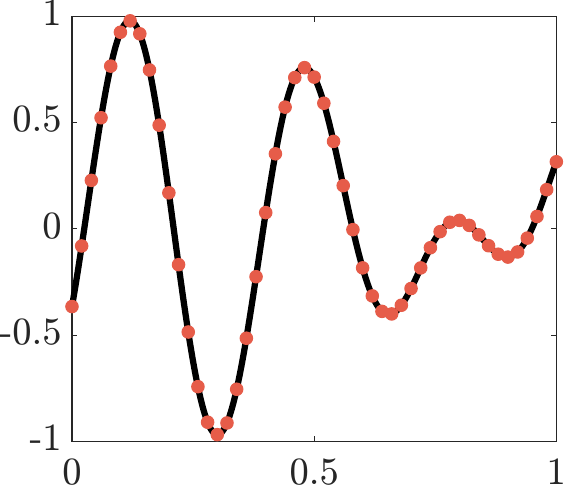}
        \caption[]%
        {{\footnotesize $u(x,t=1/2)$}}
        \label{u frame 2}
    \end{subfigure}
    \begin{subfigure}[b]{0.3\textwidth}
        \centering 
        \includegraphics[width=.8\textwidth]{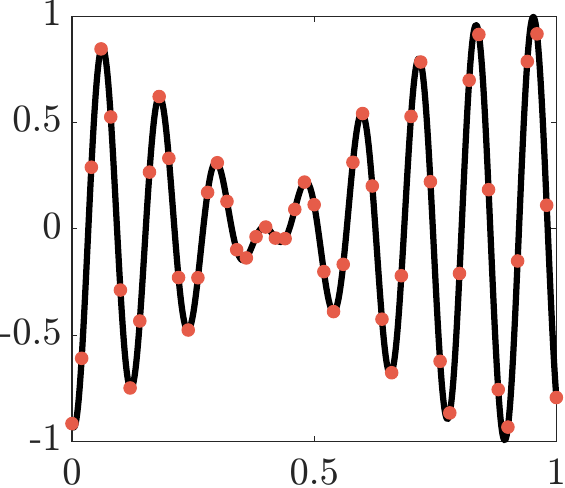}
        \caption[]%
        {{\footnotesize $\xi(x,t=1/2)$}}
        \label{ksi frame 2}
    \end{subfigure}
    
 \medskip
    \begin{subfigure}[b]{0.3\textwidth}
        \centering 
        \includegraphics[width=.8\textwidth]{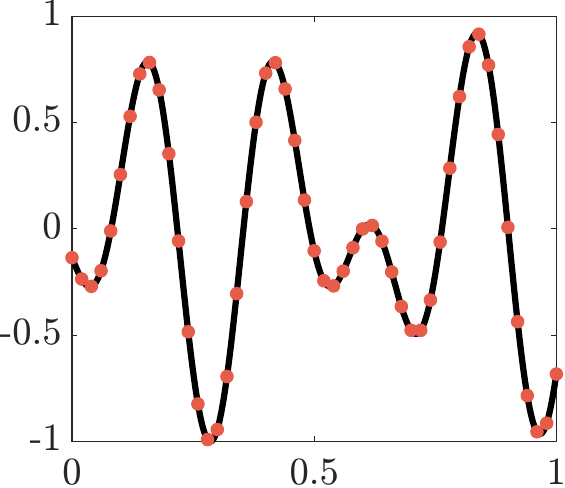}
        \caption[]%
        {{\footnotesize $\eta(x,t=1)$}}
        \label{eta frame 3}
    \end{subfigure}
    \begin{subfigure}[b]{0.3\textwidth}
        \centering
        \includegraphics[width=.8\textwidth]{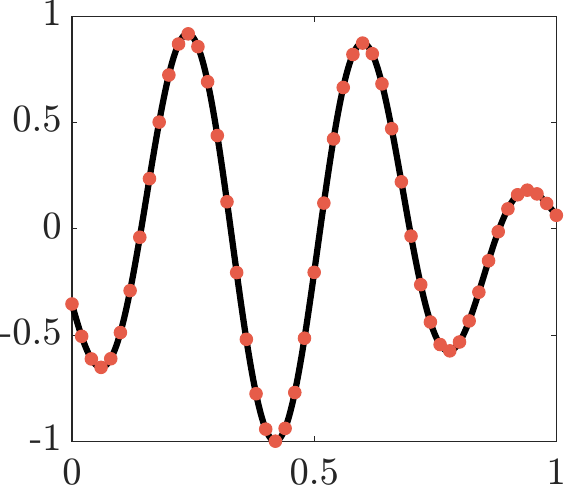}
        \caption[]%
        {{\footnotesize $u(x,t=1)$}}
        \label{u frame 3}
    \end{subfigure}
    \begin{subfigure}[b]{0.3\textwidth}
        \centering 
        \includegraphics[width=.8\textwidth]{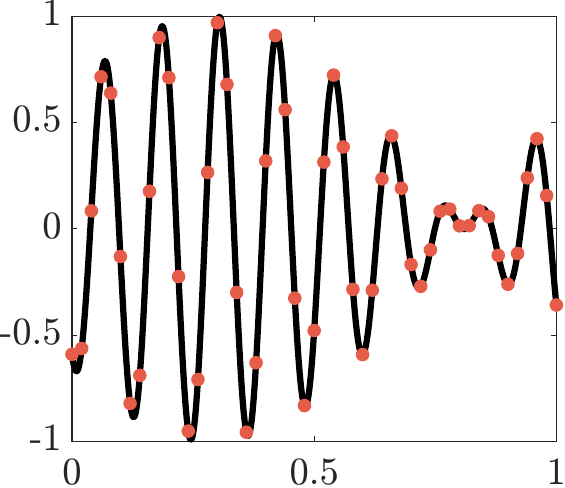}
        \caption[]%
        {{\footnotesize $\xi(x,t=1)$}}
        \label{ksi frame 3}
    \end{subfigure}
    \caption{\emph{Convergence studies.} Snapshots at various times of the manufactured solution corresponding to \autoref{MMS 1D} and \autoref{MMS 1D ksi}. Red dots denote the discretization corresponding to the coarsest mesh size employed in the convergence analysis.}
    \label{fig:1D sinusFast}
\end{figure*}

In order to assess convergence, the following manufactured reference solution over a 1D spatial domain $[0,1]$ is considered, which contains complex variations in time and space:
\begin{equation} \label{MMS 1D}
    \begin{cases}
        \eta(x,t)=\sin(5x-3t)\sin(23x-5t), \\ u(x,t)=\cos(2.5x-t)\cos(17x-4t).
    \end{cases}
\end{equation}
Since a principal objective of the SWE methodology introduced in this work is to treat highly dynamic sea-floor movement, a corresponding time-dependent (coseismic) displacement $\xi$ is additionally considered:
\begin{equation} \label{MMS 1D ksi}
    \xi(x,t)=\sin(53x-13t)\sin(3x-15t).
\end{equation}
The coefficients employed for $\xi$ have been intentionally chosen to be of higher frequency in both space and time than either $\eta$ or $u$, inspired by seismogenic configurations where earthquake dynamics are significantly faster than the corresponding water waves. With an additional consideration of a flat bathymetric (unperturbed) profile given by $h_0(x)=5$, the total $h(x,t)=\eta(x,t)+h_0(x)-\xi(x,t)$ remains positive at all times $t$. Snapshots of the manufactured $\eta, u$ and $\xi$ given above are presented in \autoref{fig:1D sinusFast}, overlaid with the corresponding discrete representation that is produced by the coarsest discretization considered.

\begin{figure}[!t]
	\centering
	\includegraphics[width=.8\textwidth]{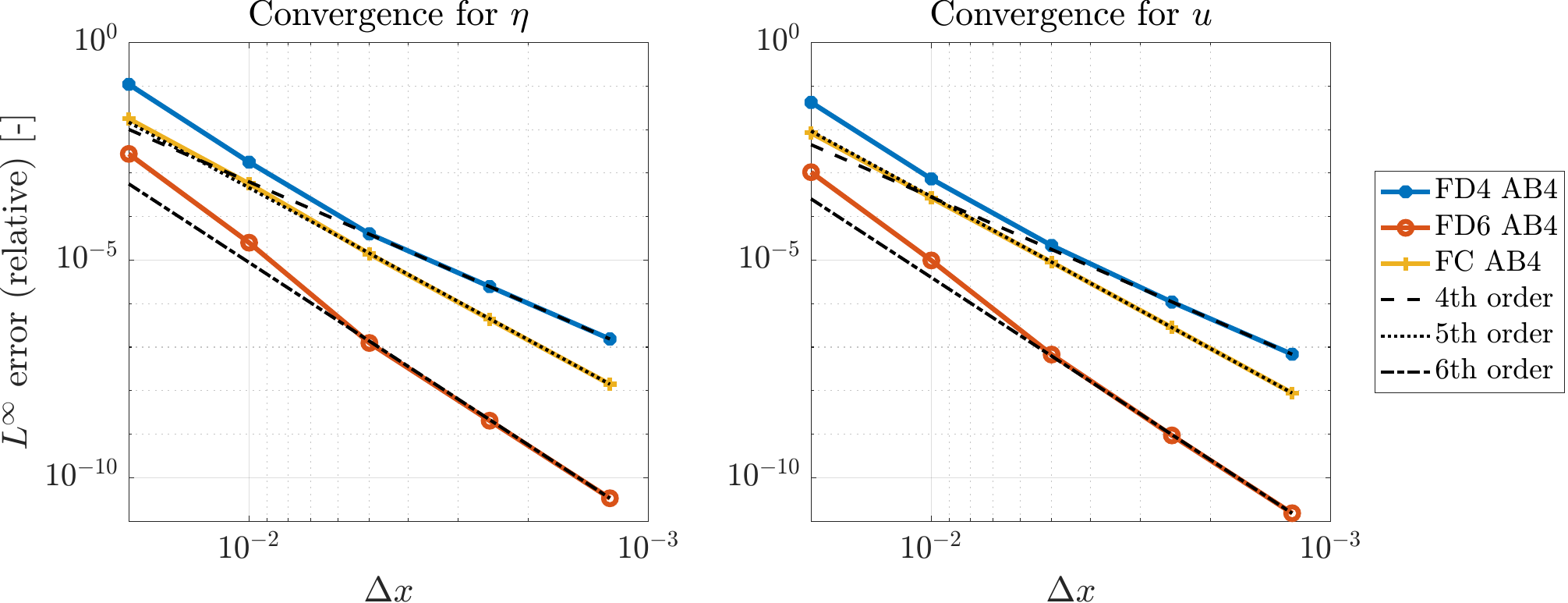}
	\caption{\emph{Convergence studies.} Errors as a function of increasing spatial resolution for FC, FD4, and FD6 applied to the 1D problem corresponding to the exact solution given by \autoref{MMS 1D} and \autoref{MMS 1D ksi}.}
	\label{fig:1D convergence}
\end{figure}

Simulations are performed employing $g=1$ in the shallow water equations for a series of increasing spatial resolutions of $\Delta x = 0.02, 0.01, 0.005, 0.00025, 0.000125$. All simulations are advanced to a final time of $t_\text{max}=1$ employing a timestep of $\Delta t = 0.2 \times \text{CFL} \times \Delta x_\text{min}/c$ (a timestep fine enough so that the errors are dominated by their respective spatial discretization sizes). \autoref{fig:1D convergence} presents the $L^\infty$ errors (maximum errors over all space and time) for both the new FC solver as well as fourth- and sixth-order finite differences (all three employing the same time integration method). Such errors are calculated (e.g. for $\eta$) by the expression given by
% The size of the domain, $n$, evolves between 1 and 40, for a total of 40 iterations.
\begin{equation}\label{eq:error_dispersion}
    L^{\infty} \ \text{error}=\frac{\max_{x,t} \left(|\eta(x,t)-\eta_{\text{ref}}(x,t)| \right)}{\max_{x,t} \left(|\eta_{\text{ref}}(x,t)| \right)}.
\end{equation}
As expected by construction (using $R=L=5$ matching points for the Gram polynomial basis), the tests here indicate that the FC method converges at a fifth-order rate to $10^{-8}$. After this point, the rate of convergence is observed to slow down as errors begin to be dominated by the timestep choice at a maximum CFL constant. This limit can be pushed to $10^{-10}$ or below by further refining the timestep well below the CFL condition such that errors are fully dominated by the spatial discretization (as have been observed for other FC-based methods~\cite{albin2011}). Future work entails improving the current limit, possibly by investigating the construction of the interpolating polynomial basis as well as investigating the choice of spectral filter~\cite{amlani2016}. Although sixth-order FD converges at an (expectedly) higher rate, the following section demonstrates that it still incurs numerical pollution errors that accumulate and make it less performant for long-time or long-distance propagation.
%future work entails improving the current implementation, possibly by investigating the choices of filter and FC parameters~\cite{amlani2016}. 

\begin{figure*}[!t]
    \centering
    \captionsetup[subfigure]{oneside,margin={0cm,0cm}}
    \begin{subfigure}[b]{0.24\textwidth}
        \centering 
        \includegraphics[width=\textwidth]{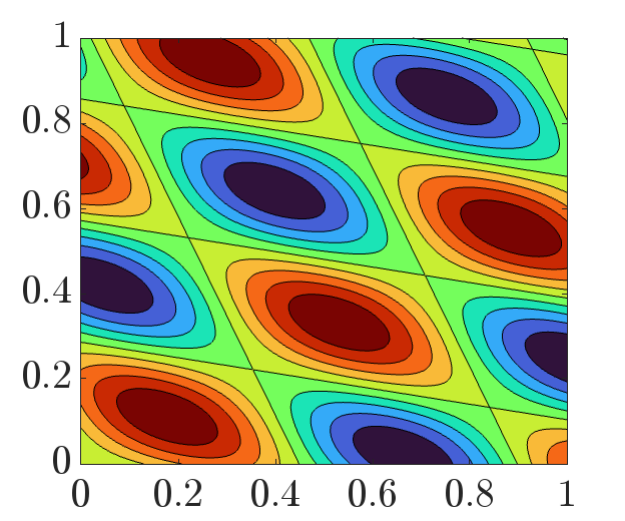}
        \caption[]%
        {{\footnotesize $\eta(x,y,t=0)$}}
        \label{2D eta frame 1}
    \end{subfigure}
    \begin{subfigure}[b]{0.24\textwidth}
        \centering
        \includegraphics[width=\textwidth]{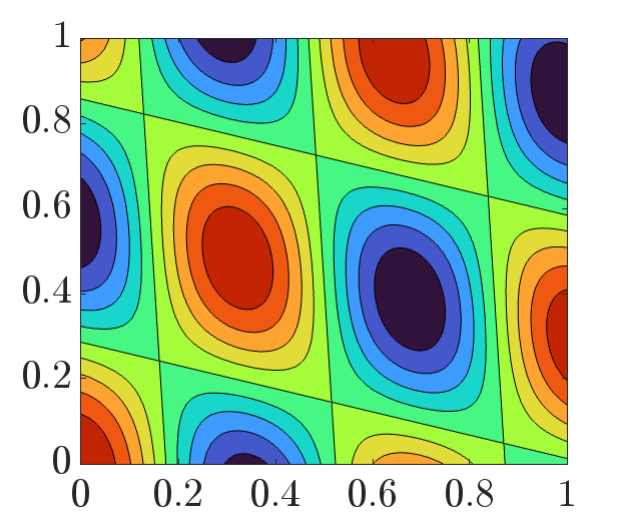}
        \caption[]%
        {{\footnotesize $u(x,y,t=0)$}}
        \label{2D u frame 1}
    \end{subfigure}
    \begin{subfigure}[b]{0.24\textwidth}
        \centering
        \includegraphics[width=\textwidth]{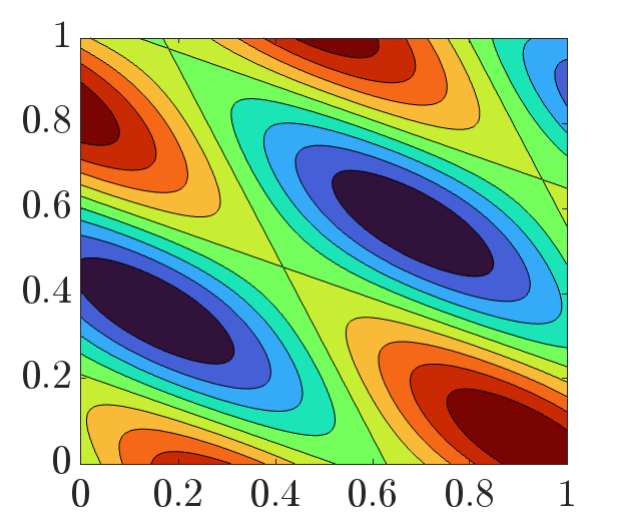}
        \caption[]%
        {{\footnotesize $v(x,y,t=0)$}}
        \label{2D v frame 1}
    \end{subfigure}
    \begin{subfigure}[b]{0.24\textwidth}
        \centering 
        \includegraphics[width=\textwidth]{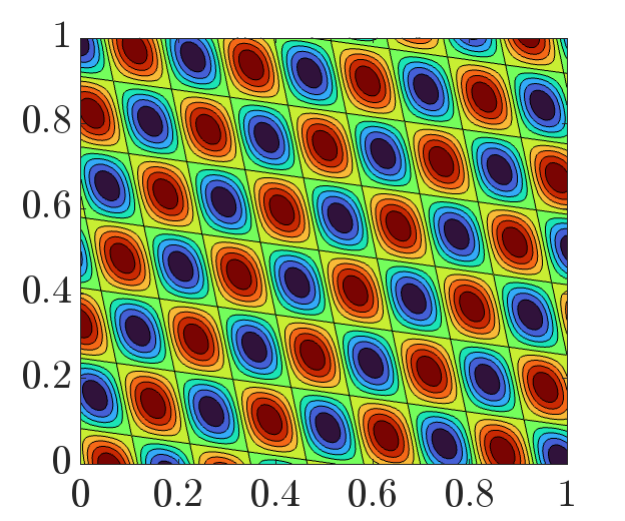}
        \caption[]%
        {{\footnotesize $\xi(x,y,t=0)$}}
        \label{2D ksi frame 1}
    \end{subfigure}
    \vskip\baselineskip
    \begin{subfigure}[b]{0.24\textwidth}
        \centering 
        \includegraphics[width=\textwidth]{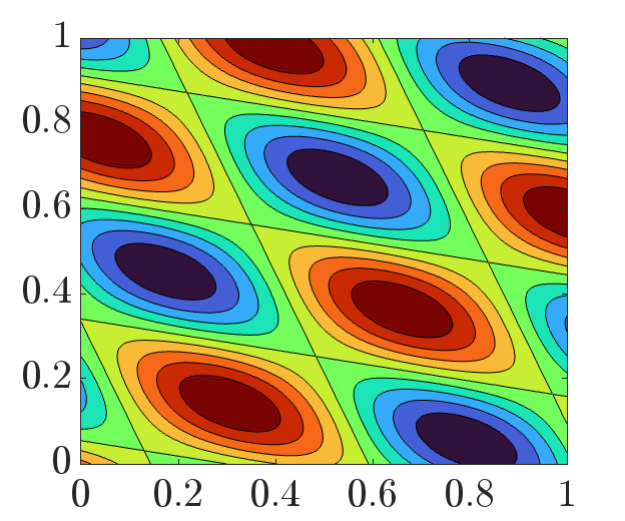}
        \caption[]%
        {{\footnotesize $\eta(x,y,t=1/2)$}}
        \label{2D eta frame 2}
    \end{subfigure}
    \begin{subfigure}[b]{0.24\textwidth}
        \centering
        \includegraphics[width=\textwidth]{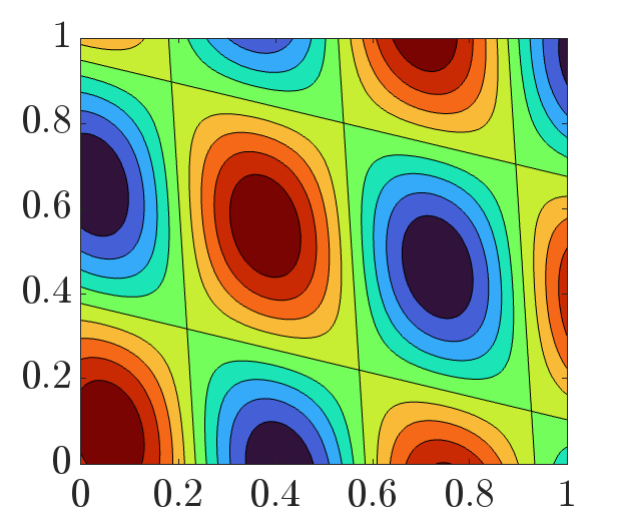}
        \caption[]%
        {{\footnotesize $u(x,y,t=1/2)$}}
        \label{2D u frame 2}
    \end{subfigure}
    \begin{subfigure}[b]{0.24\textwidth}
        \centering
        \includegraphics[width=\textwidth]{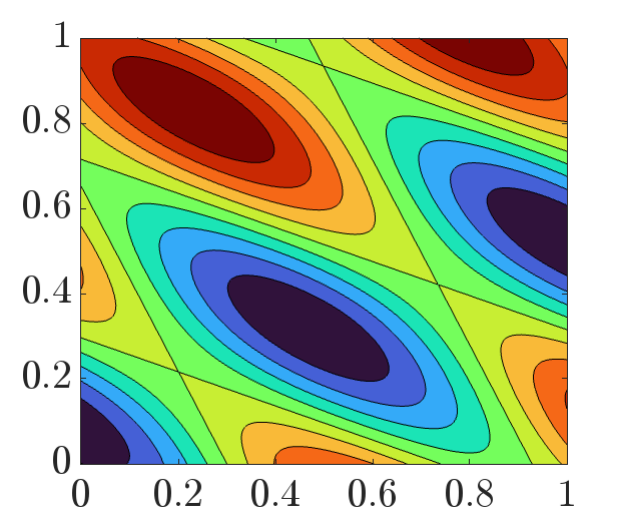}
        \caption[]%
        {{\footnotesize $v(x,y,t=1/2)$}}
        \label{2D v frame 2}
    \end{subfigure}
    \begin{subfigure}[b]{0.24\textwidth}
        \centering 
        \includegraphics[width=\textwidth]{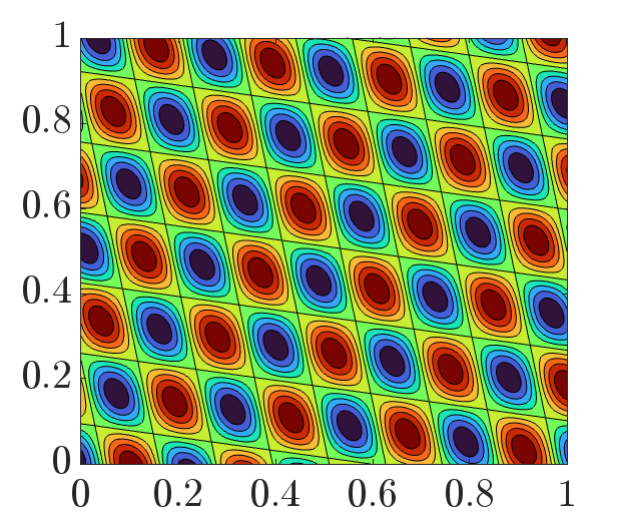}
        \caption[]%
        {{\footnotesize $\xi(x,y,t=1/2)$}}
        \label{2D ksi frame 2}
    \end{subfigure}
    \vskip\baselineskip
    \begin{subfigure}[b]{0.24\textwidth}
        \centering 
        \includegraphics[width=\textwidth]{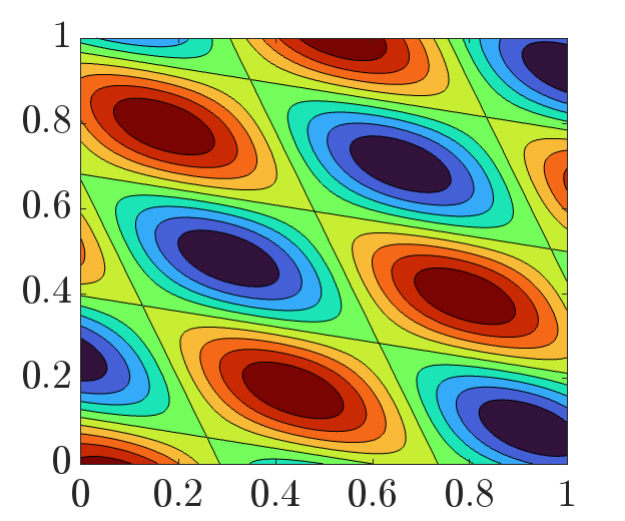}
        \caption[]%
        {{\footnotesize $\eta(x,y,t=1)$}}
        \label{2D eta frame 3}
    \end{subfigure}
    \begin{subfigure}[b]{0.24\textwidth}
        \centering
        \includegraphics[width=\textwidth]{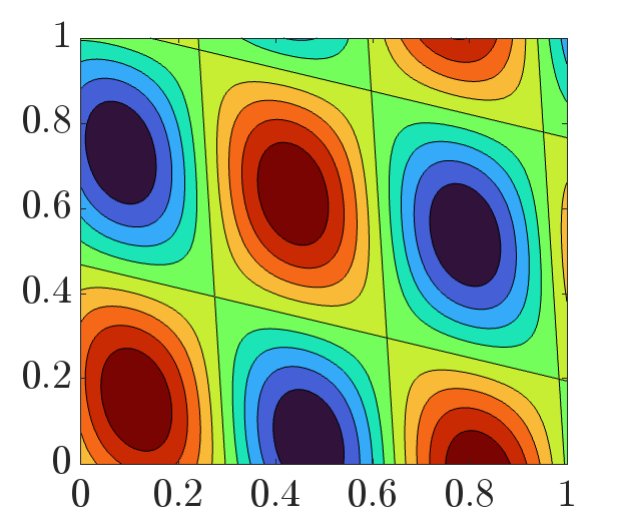}
        \caption[]%
        {{\footnotesize $u(x,y,t=1)$}}
        \label{2D u frame 3}
    \end{subfigure}
    \begin{subfigure}[b]{0.24\textwidth}
        \centering
        \includegraphics[width=\textwidth]{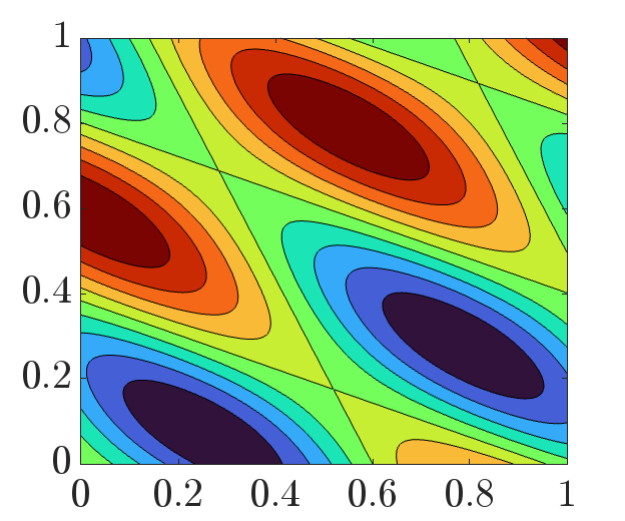}
        \caption[]%
        {{\footnotesize $v(x,y,t=1)$}}
        \label{2D v frame 3}
    \end{subfigure}
    \begin{subfigure}[b]{0.24\textwidth}
        \centering 
        \includegraphics[width=\textwidth]{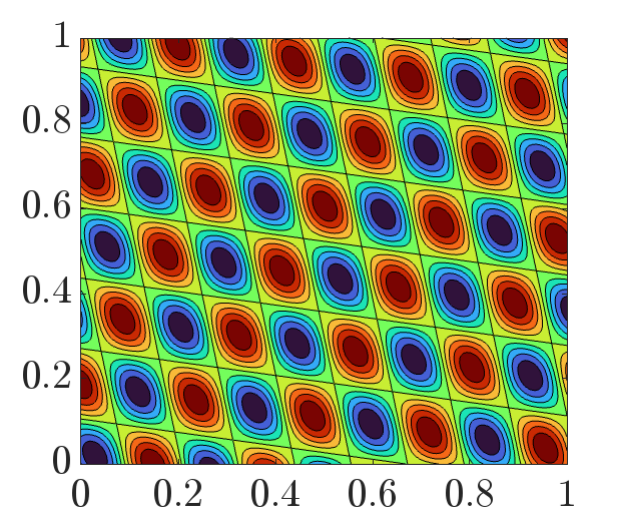}
        \caption[]%
        {{\footnotesize $\xi(x,y,t=1)$}}
        \label{2D ksi frame 3}
    \end{subfigure}
    \caption{\emph{Convergence studies.} Snapshots at various times of the 2D manufactured solutions corresponding to \autoref{MMS 2D} and~\autoref{MMS xi 2D}. The color scale ranges from -1 (dark blue) to 1 (dark red). }
    \label{fig:2D sinusFast}
\end{figure*}

\begin{figure}[!t]
	\centering
	\includegraphics[width=.9\textwidth]{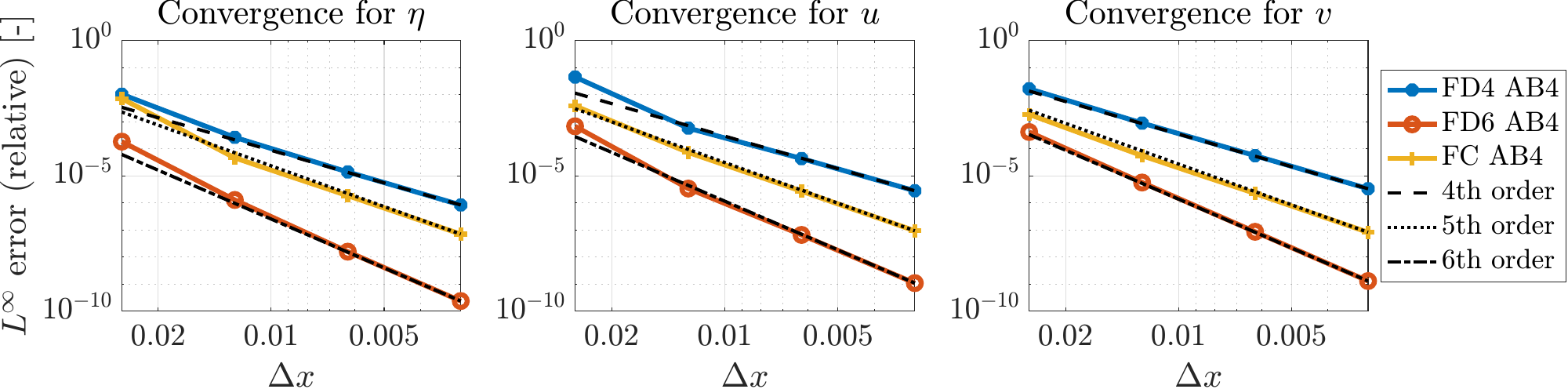}
	\caption{\emph{Convergence studies.} Maximum errors over all space and time in $\eta$ (left), $u$ (middle), and $v$ (right), as a function of increasing spatial resolution, for FD4 (blue), FD6 (orange), and FC (yellow) applied to the 2D problem whose exact solution is given by \autoref{MMS 2D} and \autoref{MMS xi 2D}. Here, $\Delta y = \Delta x$.}
	\label{fig:2D convergence}
\end{figure}

% \subsection{Convergence of the 2D solver}
% Convergence of the 2D solver can be similarly assessed as in \autoref{sec:1Dconvergence}. Consider a 2D manufactured solution given by
Convergence in 2D can be similarly assessed: consider a 2D manufactured solution given by
\begin{equation} \label{MMS 2D}
    \begin{cases}
        \eta(x,y,t) = \sin(7x+3y-2t)\sin(2x+11y-1.2t),
        \\ u(x,y,t) = \cos(1.5x+5.5y-t)\cos(9x+0.5y-1.1t),
        \\ v(x,y,t) = \sin(5x+2.3y-3t)\cos(3x+7.5y-1.3t).
    \end{cases}
\end{equation}
Similarly to the 1D case, a highly dynamic ground motion, given by
\begin{equation}\label{MMS xi 2D}
    \xi(x,y,t) = \sin(3x+19y-13t)\sin(27x+5y-15t),
\end{equation}
is imposed on the sea floor (chosen, again, to be faster and of higher frequency than the solutions of \autoref{MMS 2D} in order to mimic realistic seismic scenarios). Snapshots in time of the manufactured solutions $\eta, u$ and $\xi$ are presented in \autoref{fig:2D sinusFast}. Simulations are again performed employing $g=1$ for a sequence of finer discretizations in both $x$ and $y$, i.e., $\Delta x= \Delta y = 0.025, 0.0125, 0.00625, 0.003125$. Similarly to the 1D case, the final time is taken to be a complete cycle, i.e., $t_\text{max}=1$, employing a timestep that is one-fifth of the maximum allowable timestep from the CFL condition (\autoref{remark:FCparams}) so that errors are dominated by the spatial discretization. \autoref{fig:2D convergence} presents the $L^\infty$ errors (maximum errors over all space and time) for both the 2D FC solver as well as fourth- and sixth-order finite differences (all three employing the same time integration method). As expected by construction (using $R=L=5$ matching points for the Gram polynomial basis), the FC method converges at the expected fifth-order rate, verifying the 2D solver implementation.

% Again,  Although sixth-order FD converges at an (expectedly) higher rate, the following section demonstrates that it still incurs numerical pollution errors that compound and make it less performant for long-time or long-distance propagation. 

% The timestep corresponds to the finest spatial discretization to ensure that the convergence is only in space: $\Delta t = \text{CFL}\times (\Delta x)_\text{min}/(2c)$ with $\text{CFL}=0.05$. The convergence of the 2D solver is presented \autoref{fig:2D convergence}.

\subsection{Dispersion experiments}\label{sec:dispersion}
\begin{figure}[!t]
	\centering
	\includegraphics[width=.85\textwidth]{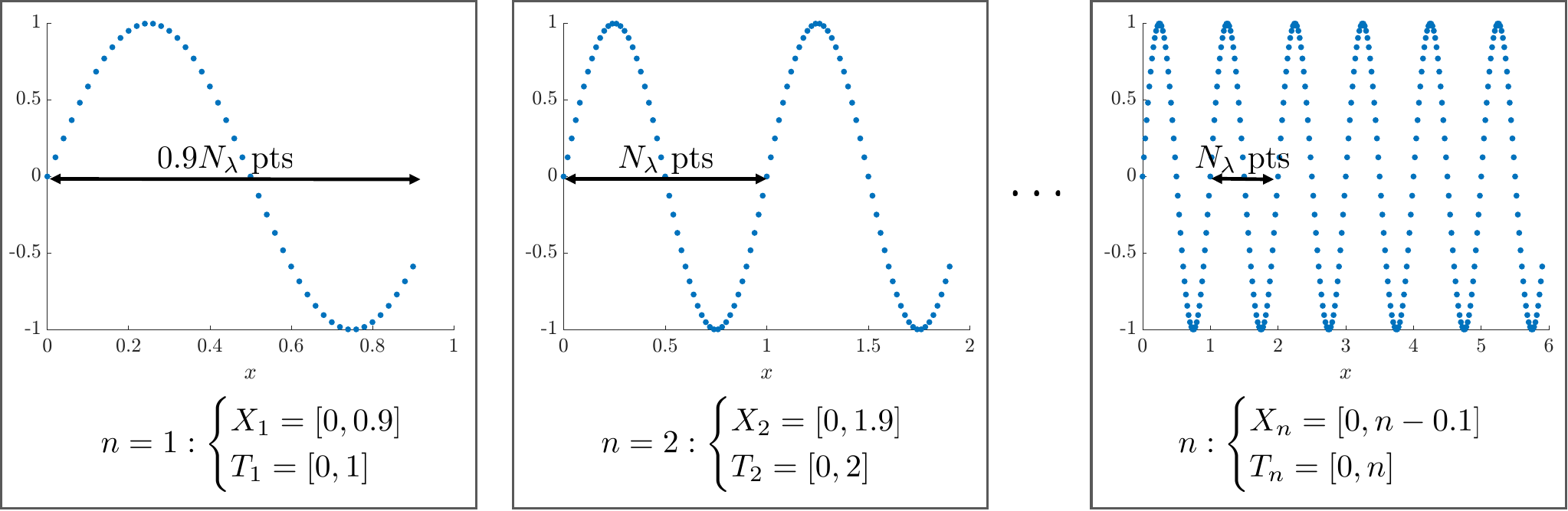}
	\caption{\emph{Dispersion studies.} Illustrations of the manufactured free surface displacement given by \autoref{eq:dispersion} employed for assessing numerical dispersion, as a function of increasing problem size $n$, for a fixed discretization of $N_\lambda$ points per wavelength.}
	\label{fig:dispersion mms}
\end{figure}

Numerical dispersion (or "pollution") errors of the FC-based solver can be studied and compared with other methods by considering a manufactured solution in the form of a propagating sine wave given by
\begin{equation}\label{eq:dispersion}
\begin{cases}
    \eta(x,t)= \sin(2\pi(x-t)),\\
    u(x,t)=\cos(2\pi(x-t)),
    \end{cases}
\end{equation}
on the space-time intervals corresponding to $x\in X_n = \left[0,n-0.1\right]$ and $t\in T_n = \left[0,n\right]$ for a number $n \in \mathbb{N}^*$ of wavelengths in the spatial domain, where each wavelength is discretized with a number $N_\lambda$ of points per wavelength. The temporal interval is chosen to correspond to a full temporal cycle of the solution given by \autoref{eq:dispersion}, and the spatial domain is defined on the slightly smaller interval ($[0,n-0.1]$) than the interval of periodicity ($[0,n]$) in order to ensure the relevance of employing the FC procedure. The still water depth is taken to be $h_0(x)=5$ throughout the domain, and no ground motion is considered (i.e., $\xi(x,t)=0$ at all times). Snapshots of the solutions corresponding to increasing values of $n$ at a fixed discretization of $N_\lambda$ points per wavelength are presented in \autoref{fig:dispersion mms}. 

\begin{figure}[!t]
	\centering
	\includegraphics[width=.65\textwidth]{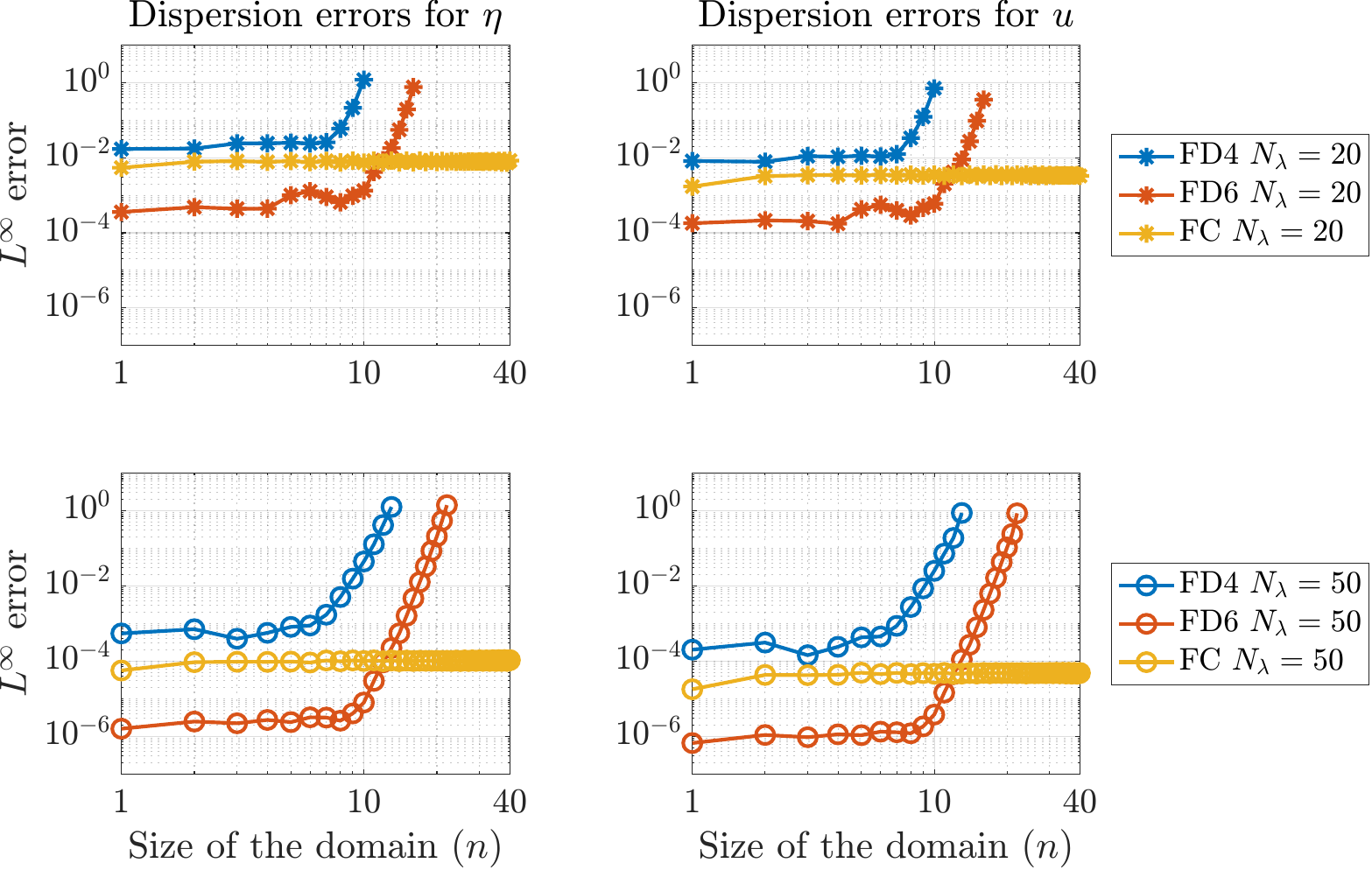}
	\caption{\emph{Dispersion studies.} Maximum errors over all space and time in $\eta$ (left) and $u$ (right), as a function of increasing problem size $n$, for FD4 (blue), FD6 (orange), and FC (yellow) applied to the problem whose exact solution is given by \autoref{eq:dispersion}. Two different numbers of fixed points per wavelength, $N_\lambda$, are considered: $N_\lambda = 20$ (top row) and $N_\lambda = 50$ (bottom row). }
	\label{fig:1D dispersion}
\end{figure}

As $n$ increases, the solution propagates over a longer distance (and a correspondingly longer time interval). For the following tests, two fixed numbers of points per wavelength (i.e., resolution) are considered: $N_\lambda \in \{20, 50 \}$ (corresponding to a spatial discretization of $\Delta x=1/N_\lambda$). The timestep is chosen accordingly for each problem size via the CFL condition (see \autoref{remark:FCparams}). \autoref{fig:1D dispersion} presents the maximum $L^\infty$ errors as a function of increasing problem size (i.e., $n$, for up to $n=40$ wavelengths), comparing both FC and FD. Such errors are calculated over all space and all time by the $L^{\infty}$ norm given by \autoref{eq:error_dispersion}.

 Clearly, the error of FD increases with increasing $n$, eventually leading to instability when errors reach levels where the corresponding solutions violate the constraint $h>0$ of the shallow water system. Such errors begin to grow after a certain size, even when increasing the FD stencil from fourth to sixth order (a higher order than FC and much higher than those commonly-used in FD-based tsunami solvers) or by increasing the number $N_\lambda$ of points per wavelength. Indeed, FD and FEM-based methods are well-known to incur numerical dispersion/pollution errors~\cite{Durran1999, Deraemaeker1999}, requiring a substantial refinement of the corresponding meshes in order to capture solutions within a given accuracy for long-time and long-distance propagation. On the other hand, the error from the proposed FC-based solver remains constant, implying that it successfully avoids these pollution errors. This results in a level of error that is solely determined by the number of points per wavelength, no matter the problem size, which is very appealing for treating large-scale wave propagation problems.
%  The conclusion is that Fourier Continuation allows to solve the PDE while the wave propagates without additional efforts to refine the mesh.

\autoref{table:dispersion} presents values of the errors and computational times for FC-based simulations corresponding to different combinations of problem sizes $n$ and characteristic discretizations $N_\lambda$ (i.e., three of the points found in the FC error curves of \autoref{fig:1D dispersion}). Additionally presented are the corresponding discretizations and computational times that are required for FD4-based and FD6-based simulations to achieve similar errors as FC. For $n=20$ wavelengths at a discretization of $N_\lambda=20$ points per wavelength ("Test 1"), FC achieves its errors with far fewer points (4.5 times fewer than FD6 and 31.5 times fewer than FD4) and significantly shorter runtimes (around $211$ times faster than FD4 and around $6.5$ times faster than FD6). Even when increasing the number of points per wavelength $N_\lambda$, a similar speedup is observed ("Test 2"). The qualities of FC with respect to numerical dispersion errors are further evident when increasing the problem size to a higher frequency, from $n=20$ wavelengths to $n=30$ wavelengths ("Test 3"): its speedup factor increases from $9$ to almost $165$ even when compared to an ostensibly higher-order method such as FD6 (which, for this test, still incurs five times larger errors than FC at the given discretization). 

% Such speedups persist even  

% \red{STOPPED HERE.   }

% Even when increasing $N_\lambda$ for a fixed problem size $n$ (``Test 3"), such a speedup per

% and $16$ times faster than FD6 (corresponding to $n=20$ wavelengths, each discretized by $N_\lambda = 20$) and 4In each case, the qualities of FC with respect to numerical dispersion errors are evident. The first part of the table (three first rows) represents the computation from \autoref{fig:1D dispersion} corresponding to $N_\lambda=20$ at the abscissa $n=20$: FC's error remains constant and the objective is to observe what refinement is necessary for FD4 and FD6 to reach similar levels of error. Then, this study was repeated for $n=30$ where the refinements needed are even more important. Augmenting the number of points per wavelength to $N_\lambda=50$ (on the last two rows) leads to error reduction with FC, which is already much lower than the FDs at $n=20$ as observed qualitatively on the curve \autoref{fig:1D dispersion} and quantitatively highlighted by the refinements factors used in \autoref{table:dispersion}.

% An additional comparison of the computations leading to similar levels of error was performed (see \autoref{table:dispersion}). Three points of the FC curves presented \autoref{fig:1D dispersion} were taken as baseline for different values of $n,N_\lambda$ while the high-order FD methods were refined in order to match FC errors.

\begin{table}[!t]
    \centering
    \small
    \begin{tabular}{|c |c |c |r |c |c |c |r |} 
         \hline
          $\text{Test}$ & $\text{Method}$ & $n$ & \multicolumn{1}{c|}{$N_\lambda$} & CFL & $\eta \ \text{error}\left[ \% \right]$ & $u \ \text{error}\left[ \% \right]$  & \multicolumn{1}{c|}{$T_{\text{solve}} \left[\text{s} \right]$} \\ [0.5ex]
          \hline
          \rowcolor{gray!20}
          $\text{1}$ & {FC} & ${20}$ & ${20}$ & 0.17 & ${7.74 \times 10^{-1}}$ & ${3.45 \times 10^{-1}}$ & ${2}$ \\
         $ $ & $ $FD4 & 20 & 650 & 0.21 & $7.76 \times 10^{-1}$ & $5.75 \times 10^{-1}$ & 422 \\
         $ $ & $ $FD6 & 20 & 89 & 0.18 & $8.71 \times 10^{-1}$ & $4.35 \times 10^{-1}$ & 13 \\ 
         \hline
              \rowcolor{gray!20}
         $\text{2}$ & {FC} & ${20}$ & ${50}$ & 0.17 & ${1.10 \times 10^{-2}}$ & ${4.70 \times 10^{-3}}$  & ${3}$ \\
         $ $ & $ $FD6 & 20 & 189 & 0.18 & $1.10 \times 10^{-2}$ & $5.00 \times 10^{-3}$ & 27 \\ 
         \hline
         \rowcolor{gray!20}
         $\text{3}$ & {FC} & ${30}$ & ${20}$ & 0.17 & ${8.03 \times 10^{-1}}$ & ${3.47 \times 10^{-1}}$ & ${3}$ \\
         $ $ & $ $FD6 & 30 & 520 & 0.18 & $7.15\times 10^{-1}$ & $3.56\times 10^{-1}$ & 494 \\ 
         \hline
    
    \end{tabular}
    \caption{\emph{Dispersion studies}. Numerical discretization sizes (in terms of the number of points per wavelengths, $N_\lambda$) and computational times (in seconds) required by FD4 or FD6 to achieve similar errors as the FC-based solver (relative to the exact solution in \autoref{eq:dispersion}) for problem sizes of various numbers of $n$ wavelengths.}
    \label{table:dispersion}
\end{table}

\section{Benchmarks} \label{Benchmarks}
In what follows, the numerical implementation of the FC-based methodology is further validated, and its physical accuracy established, through applications of the solver to a number of (classical and non-classical) benchmark problems for a wide range of realistic configurations previously proposed by others.

%  in order to demonstrate its efficacy 
% All of the results that follow are produced using the methodology and parameters outlined in \autoref{SWE solver}.

\subsection{Benchmark 1: Perturbation of a lake at rest}\label{sec:benchmark1}
\begin{figure}[!t]
	\centering
	\includegraphics[width=0.45\textwidth]{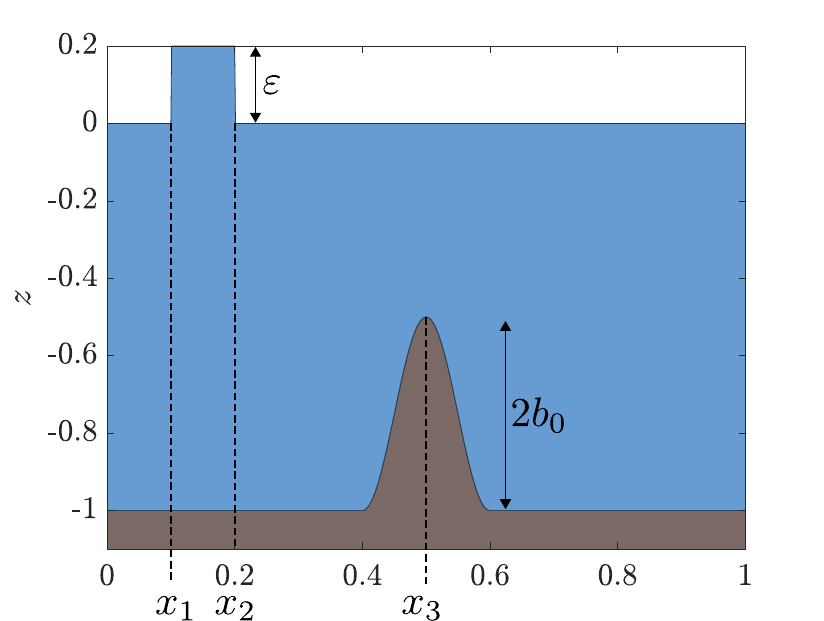}
	\caption{\emph{Benchmark 1}. The perturbed lake-at-rest configuration proposed by LeVeque \cite{leveque1998}. The total water depth is represented in blue, and the ground topography is represented in brown.}
	\label{fig:Leveque config}
\end{figure}

A classical validation of SWE solvers is to verify that the model satisfies a ``lake at rest" configuration \cite{lundgren2020,michel2016,ricchiuto2007}. Such a benchmark consists of defining a sea floor topography with a still water height and verifying that the free surface remains exactly zero. In the context of the present solver, the shallow water system does not generate perturbations coming from a topography associated with a flat (zero initial water height) free surface. Hence, in a similar spirit, the lake at rest benchmark can be modified and adapted for other methods by employing, instead, an initial perturbation of the free surface, as proposed by LeVeque \cite{leveque1998} and illustrated in \autoref{fig:Leveque config}. An initial water height perturbation $\epsilon>0$ is imposed on a portion of the free surface, i.e.,
\begin{equation}
    \eta(x)=
    \begin{cases}
        \varepsilon, & x_1 \leq x \leq x_2 \\
        0, & \text{otherwise},
    \end{cases}
\end{equation}
with a corresponding seafloor topography $s(x)$ given by
\begin{equation}
    s(x)=\begin{cases}
    s_0 \left(\cos\left( \pi (x-x_3) /\alpha \right)+1 \right), & |x-x_3|<0.1 \\
    0, & \text{otherwise}.
    \end{cases}
\end{equation}
The problem formulation is non-dimensionalized, where the still water height is given by $h_0(x)=D-s(x)$ for a constant $D\in \mathbb{R}$. 

\begin{figure*}[!t]
    \centering
    \captionsetup[subfigure]{oneside,margin={0cm,0cm}}
    \begin{subfigure}[b]{0.49\textwidth}
        \centering
        \includegraphics[width=.75\textwidth]{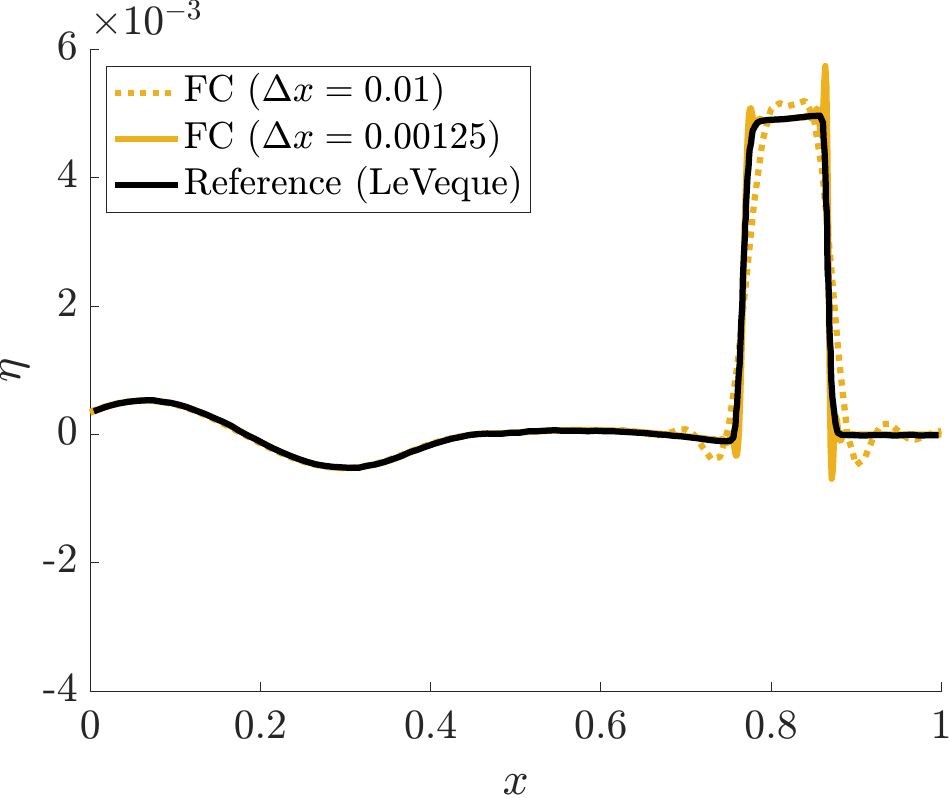}
        \caption[]%
        {{\footnotesize FC $\varepsilon=0.01$;}}
        \label{Leveque FC 0.001}
    \end{subfigure}
    \begin{subfigure}[b]{0.49\textwidth}
        \centering 
        \includegraphics[width=.8\textwidth]{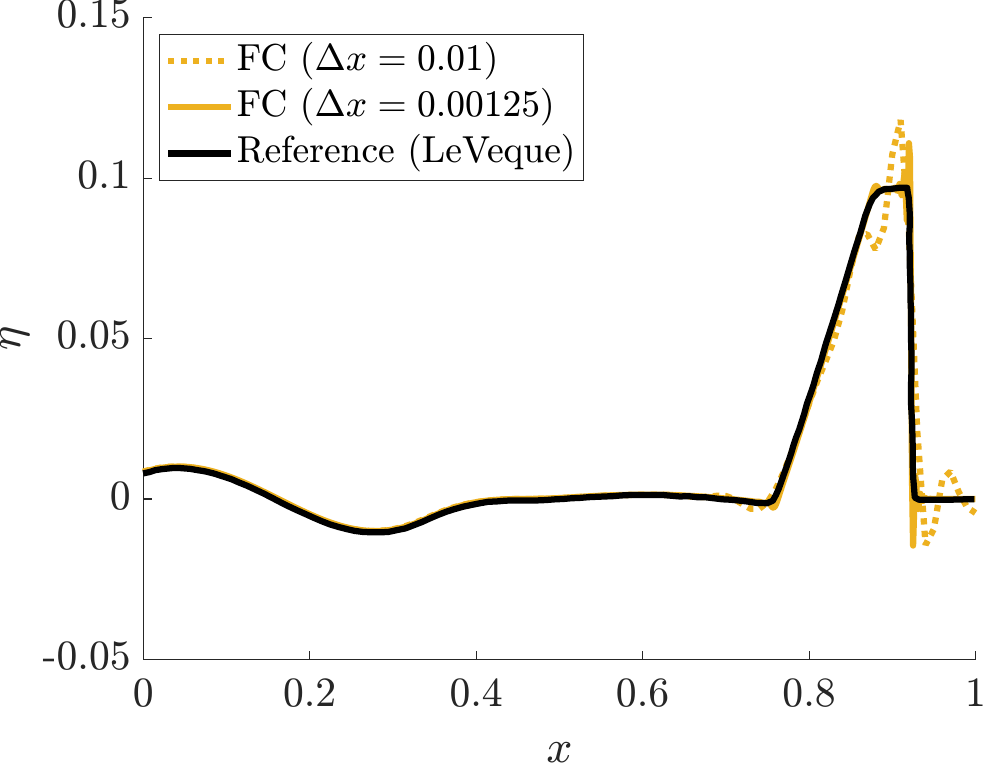}
        \caption[]%
        {{\footnotesize FC $\varepsilon=0.2$;}}
        \label{Leveque FC 0.2}
    \end{subfigure}
\medskip

    \begin{subfigure}[b]{0.49\textwidth}
        \centering
        \includegraphics[width=.75\textwidth]{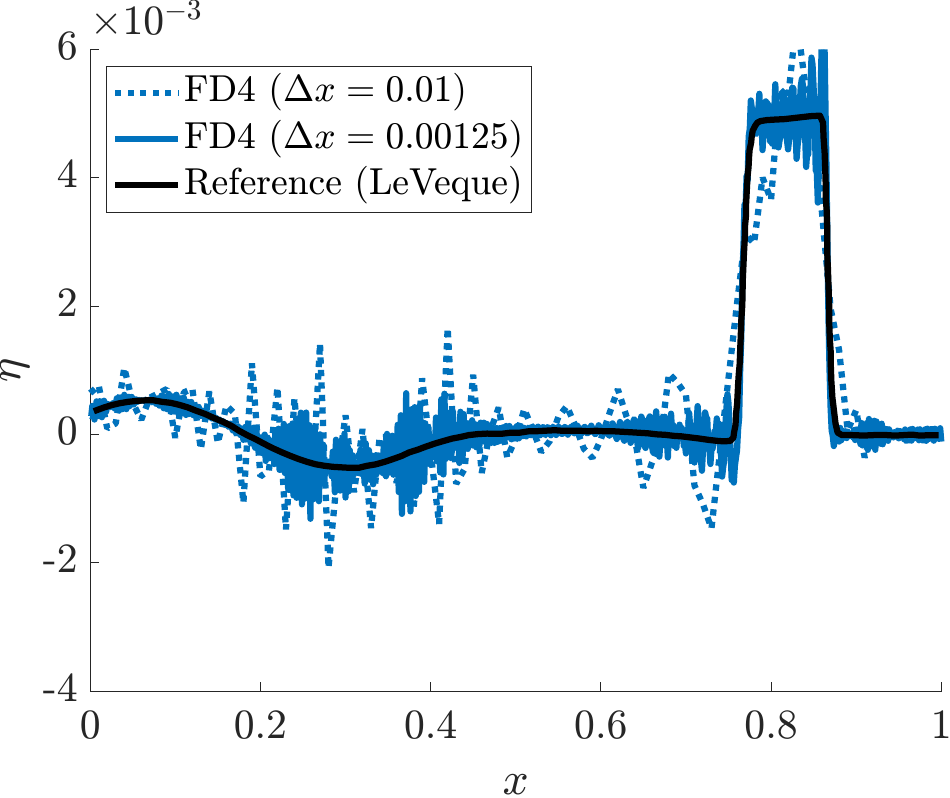}
        \caption[]%
        {{\footnotesize FD4 $\varepsilon=0.01$;}}
        \label{Leveque FD4 0.001}
    \end{subfigure}
    \begin{subfigure}[b]{0.49\textwidth}
        \centering 
        \includegraphics[width=.75\textwidth]{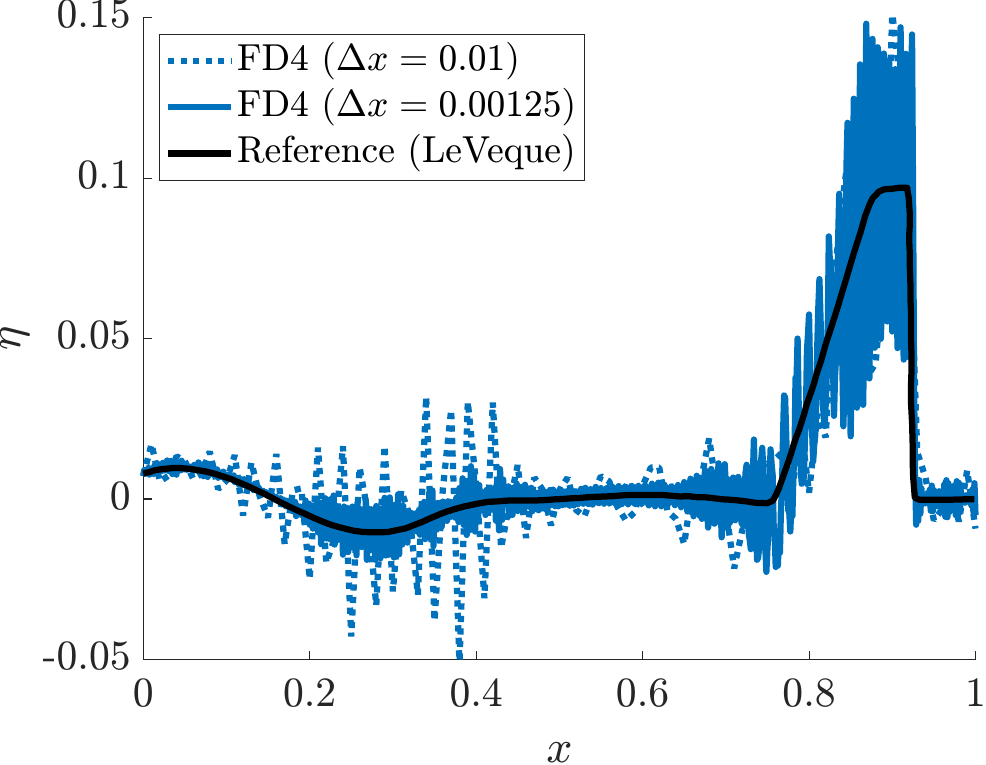}
        \caption[]%
        {{\footnotesize FD4 $\varepsilon=0.2$;}}
        \label{Leveque FD4 0.2}
    \end{subfigure}
\medskip
    \begin{subfigure}[b]{0.49\textwidth}
        \centering
        \includegraphics[width=.75\textwidth]{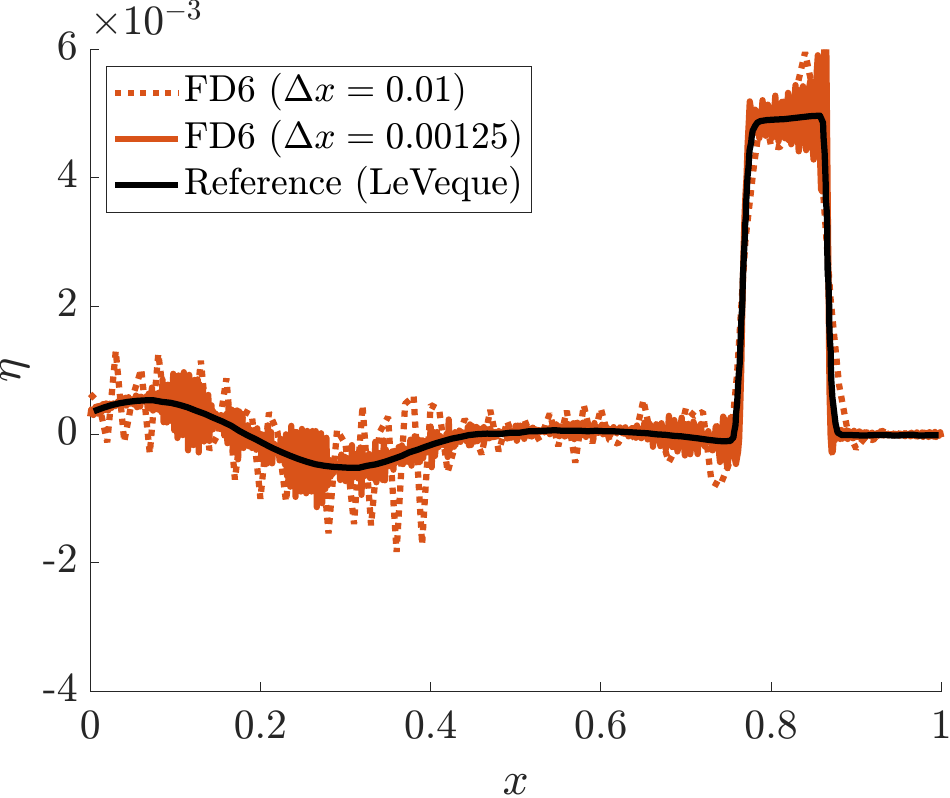}
        \caption[]%
        {{\footnotesize FD6 $\varepsilon=0.01$;}}
        \label{Leveque FD6 0.001}
    \end{subfigure}
    \begin{subfigure}[b]{0.49\textwidth}
        \centering 
        \includegraphics[width=.75\textwidth]{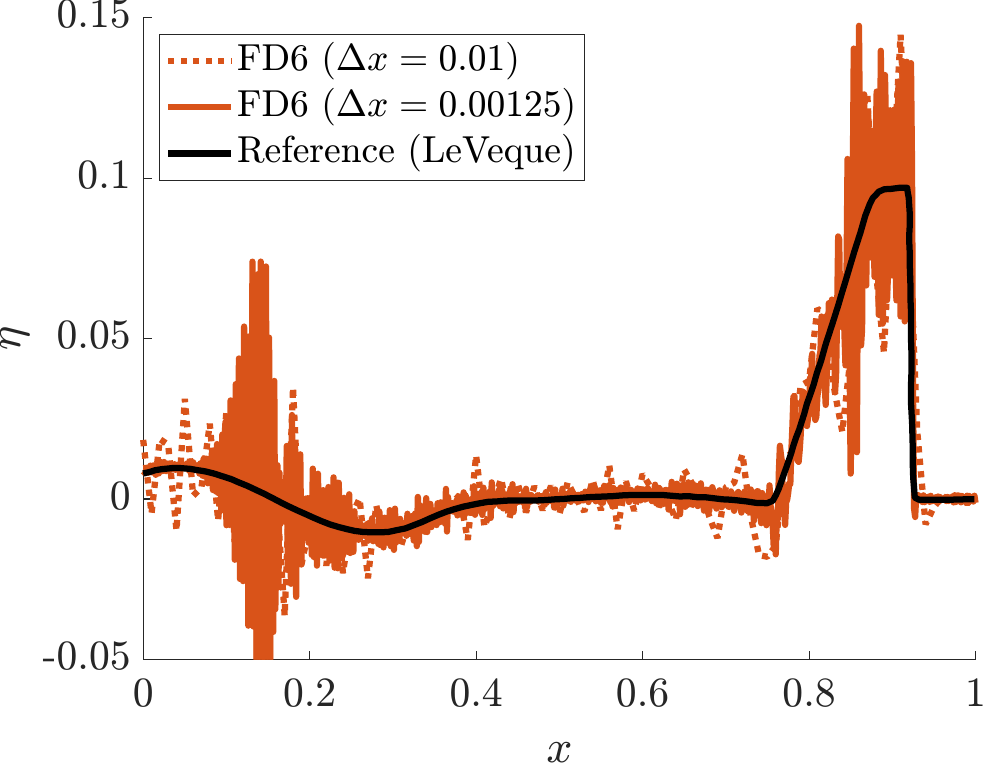}
        \caption[]%
        {{\footnotesize FD6 $\varepsilon=0.2$.}}
        \label{Leveque FD6 0.2}
    \end{subfigure}
    \caption{\emph{Benchmark 1}.  Free surface displacement snapshots at $t=0.7$ for FC (top), FD4 (middle), and FD6 (bottom) at various discretizations. Reference solutions are given by the second-order finite volume solver of \cite{leveque1998}.}
    \label{fig:LeVeque results}
\end{figure*}

\autoref{fig:LeVeque results} presents snapshots at various times of the free-surface displacement produced by the FC-based solver, an FD4 a FD4 solver, and an FD6 solver for perturbations of size $\epsilon = 0.01$ (left) and $\epsilon=0.2$ (right) and benchmark parameters corresponding to $D=1$, $g=1$, $s_0=0.25$, $x_1 = 0.1$, $x_2 = 0.2$, $x_3 = 0.5$, and $\alpha = 0.1$. All simulations are conducted on a spatial interval $x\in[-1,1]$ employing $\Delta x = 0.01$ and $\Delta x = 0.00125$. The boundary conditions are set as walls (see \autoref{eq:1Dwallbc}). The timesteps correspond to the maximum allowable by the respective CFL conditions (see \autoref{remark:FCparams}) of each method employed. Overlaid in black for each figure is the reference free surface obtained by the quasi-steady wave-propagation method extracted from LeVeque's results~\cite{leveque1998}. Here, both FD4 and FD6 suffer from spurious oscillations throughout the domain due to the initial condition and subsequently discontinuous solution. Even after refining the solution (eight times), such oscillations persist. However, even without any special treatment of such jumps, the FC-based solver, for the same discretizations as FD4 and FD6, still correctly captures the overall shape and form of the traveling perturbation after its interaction with the topographic bump (the non-flat region of the bathymetry), with no discernible oscillations away from it (treatment of discontinuities/shocks is a subject of future work and beyond the scope of the present solver). Indeed, overall results agree well with those presented in the original work~\cite{leveque1998}, and only minor oscillations, that are confined very locally near the regions of sharpest curvature and are due to the truncation of the discrete Fourier series, are observable.

% \begin{table}[!t]
%     \centering
%     \begin{tabular}{|c |c |c |c |c |c |c |c|} 
%          \hline
%           $D$ & $g$ & $b_0$ & $x_1$ & $x_2$ & $x_3$ & $\alpha$ & $\varepsilon$\\ 
%           $\left[ - \right]$ & $\left[- \right]$ & $\left[-\right]$ & $\left[-\right]$ & $\left[-\right]$ & $\left[ - \right]$ & $\left[ - \right]$ & $\left[ - \right]$\\ [0.5ex]
%          \hline
%          1 & 1 & 0.25 & 0.1 & 0.2 & 0.5 & 0.1 & $\{ 0.001, 0.2 \}$ \\ [1ex]
%          \hline
%     \end{tabular}
%     \caption{Parameters of the Benchmark 1.}
%     \label{table:LeVeque params}
% \end{table}

\subsection{Benchmark 2: 1D tsunami}\label{sec:benchmark2}
\begin{figure}[!t]
	\centering
	\includegraphics[width=0.375\textwidth]{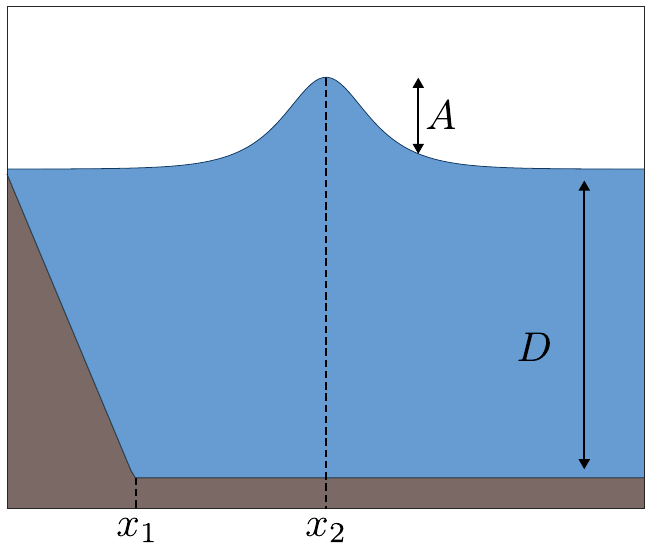}
	\caption{\emph{Benchmark 2}. The 1D tsunami configuration proposed by Setiyowati and Sumardi \cite{Setiyowati2019}. The total water depth is represented in blue, and the ground topography is represented in brown. The free surface displacement is amplified for visualization ($\times50$).}
	\label{fig:1D toy problem}
\end{figure}

A second classical benchmark problem~\cite{Synolakis2008, Setiyowati2019} concerns the shallow water propagation of a tsunami towards a sloping beach in order to mimic a more realistic tsunami configuration. An illustration of the problem setup is presented in \autoref{fig:1D toy problem}.
\ Over a domain $x\in [0, x_{max}]$, the initial condition is a solitary wave centered at $x = x_2$ and is given by
\begin{equation}
    \eta(x,0)=A \ \text{sech} \left(\sqrt{\frac{3}{4D^2}} \left(x-x_2 \right) \right)
\end{equation}
for a constant $D>0$ and wave amplitude (from the free surface) $A>0$. The fluid is considered to be initially at rest, i.e, $u(x,0)=0$, and the corresponding bathymetry is given by
\begin{equation}\label{slopingbeach}
    h_0(x)=\begin{cases}
    D+\frac{B}{x_1} \left(x-x_1 \right), & 0\leq x \leq x_1 \\ D, & x\geq x_1,
    \end{cases}
\end{equation}
where $x_1\in\mathbb{R}$ demarcates the beginning of the sloping beach. Since such a bathymetric formulation is only piecewise continuous, the results that follow employ a standard $C^\infty$ rounding of \autoref{slopingbeach} through the use of hyperbolic tangents. The corresponding boundary conditions are considered to be radiation boundary conditions (see \autoref{eq:1Dradiationbc} at both the left and right boundaries.

\begin{figure*}[!t]
    \centering
    \begin{subfigure}[b]{0.4\textwidth}
        \centering
        \includegraphics[width=\textwidth]{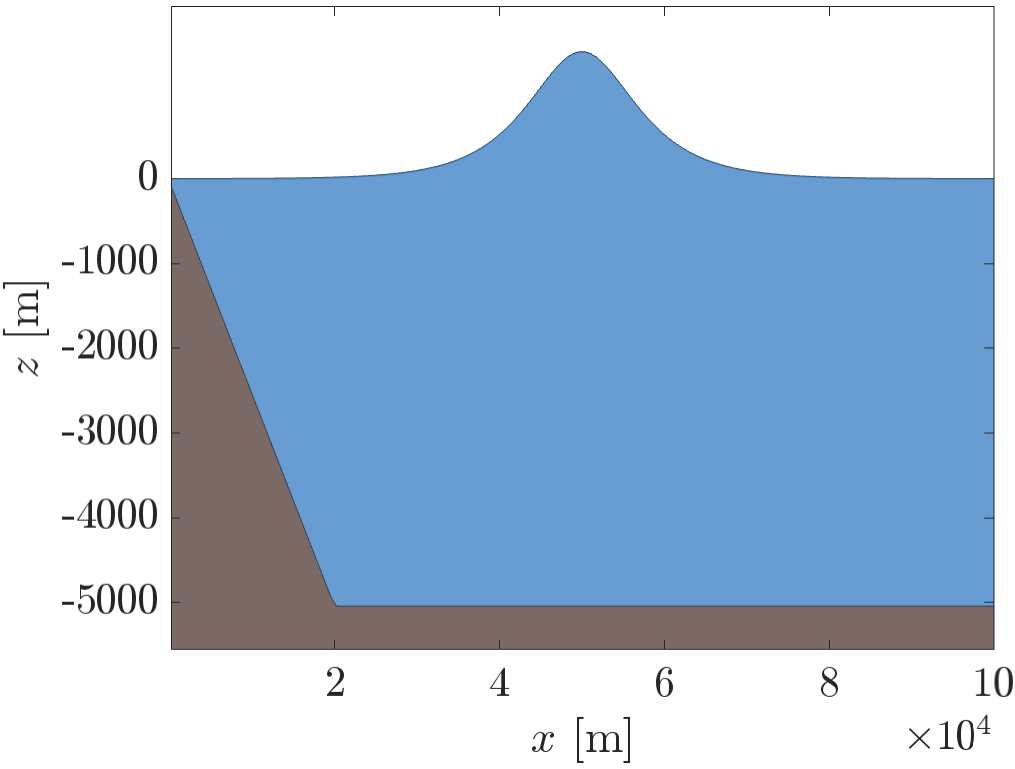}
        \caption[eta]%
        {{\footnotesize $t=0$ s;}}
        \label{init cond}
    \end{subfigure}
    \begin{subfigure}[b]{0.4\textwidth}
        \centering 
        \includegraphics[width=\textwidth]{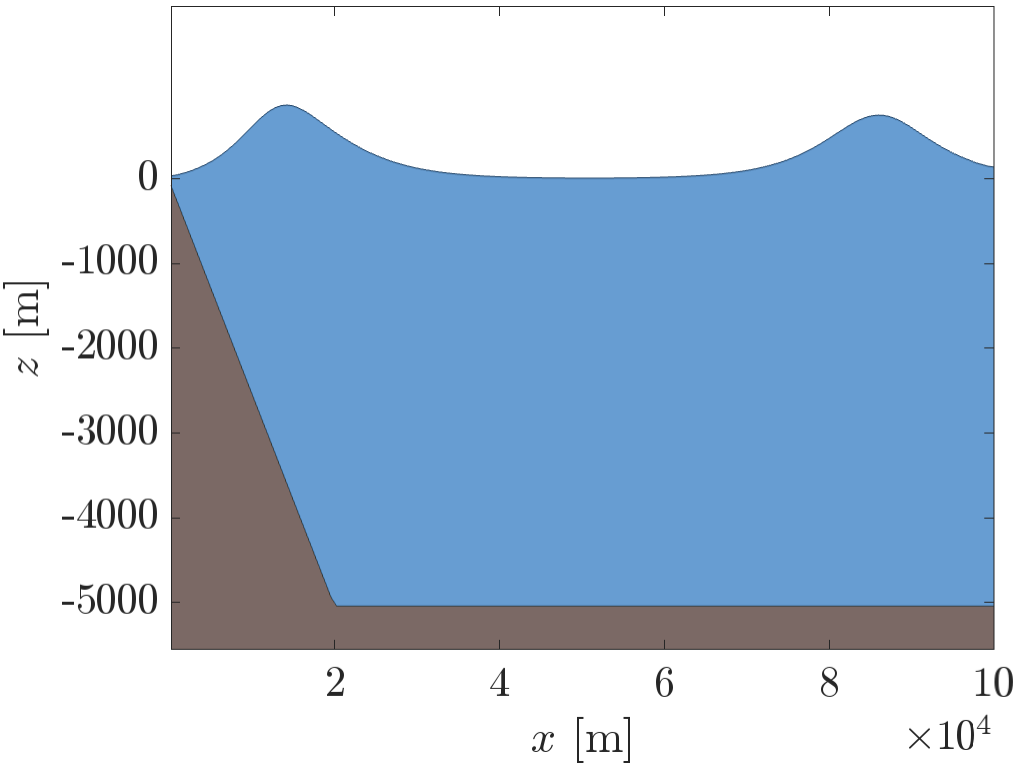}
        \caption[]%
        {{\footnotesize $t=161$ s;}}
        \label{split}
    \end{subfigure}
    
    \smallskip
    
    \begin{subfigure}[b]{0.4\textwidth}  
        \centering 
        \includegraphics[width=\textwidth]{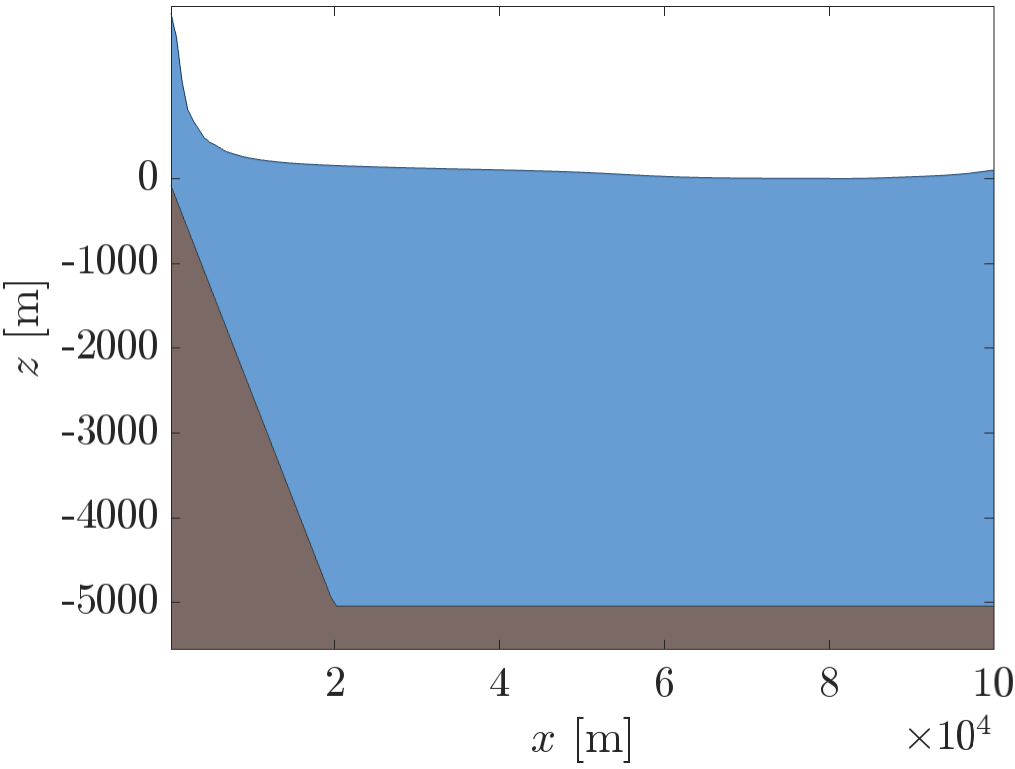}
        \caption[]%
        {{\footnotesize $t=295$ s;}}
        \label{max amplitude}
    \end{subfigure}
    \begin{subfigure}[b]{0.4\textwidth}   
        \centering 
        \includegraphics[width=\textwidth]{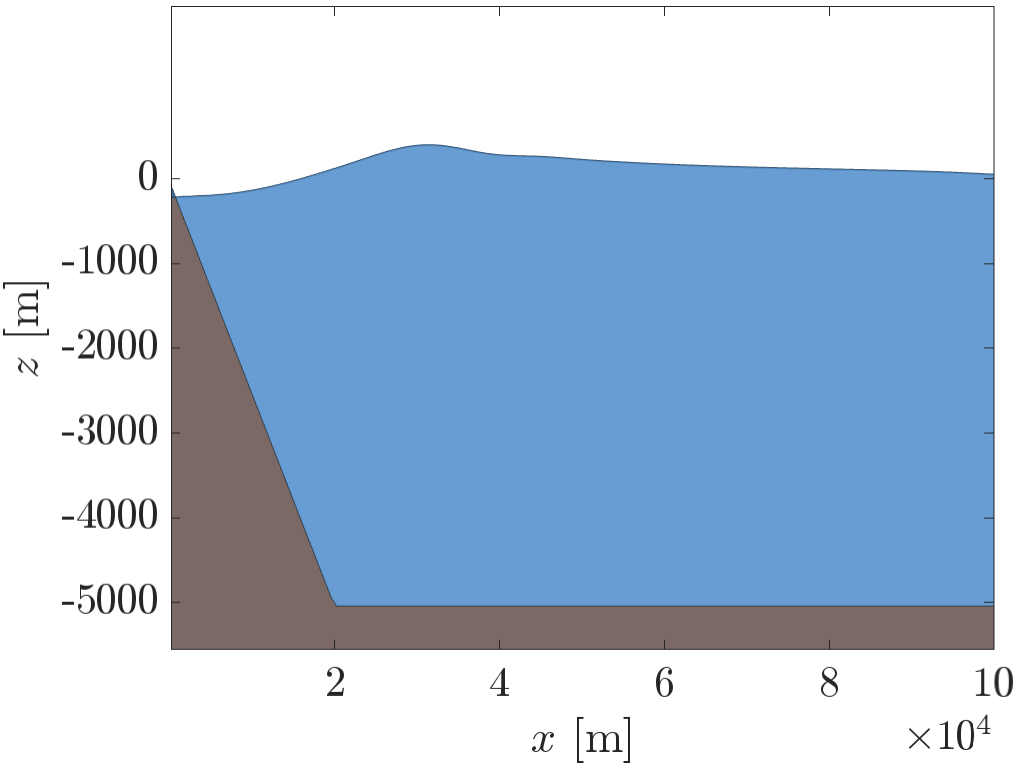}
        \caption[]%
        {{\footnotesize $t=498$ s.}}
        \label{wave escapes}
    \end{subfigure}
    \caption{\emph{Benchmark 2}. Snapshots at various times of the water perturbation $\eta$ around the mean free surface $z=0$ produced by the FC solver. The total water depth is represented in blue, and the ground topography is represented in brown. The free surface displacement is amplified for visualization ($\times50$).}
    \label{fig:Setiyowati}
\end{figure*}

\autoref{fig:Setiyowati} presents solution snapshots at various times produced by the FC solver for values of $A=30$ m, $D = 5.04$ km, $B = 5.00$ km, $x_1 = 20$ km , and $x_2 = 50$ km (all parameters adopted from \cite{Setiyowati2019}). The gravitational constant is taken to be $g=9.81$ m/s$^2$. The simulation is conducted on a spatial interval $x\in [0.2,100.0]$ km employing a discretization size of $\Delta x = 671$ m ($N_x=150$ points, similarly to \cite{Setiyowati2019}). The corresponding timestep is taken to be $\Delta t=0.511$ s, and the solution is evolved up to a final time of $t = 500$ s. The initial condition of \autoref{init cond} splits into two waves propagating towards each end of the domain (\autoref{split}). The left-going wave propagates until it encounters the sloping beach, reaching maximal amplitude at the left boundary as illustrated in \autoref{max amplitude}, while the right-going wave passes out of the region of interest (here, after 100 km \cite{Setiyowati2019}). \autoref{wave escapes} illustrates the reflected wave from the left boundary propagating back towards the right into the water. 
% \begin{table}[!t]
%     \centering
%     \begin{tabular}{|c |c |c |c |c |c|} 
%          \hline
%           $D$ & $g$ & $x_1$ & $x_2$ & $x_{\text{max}}$ & $A$\\ 
%           $\left[ \text{m} \right]$ & $\left[\text{m.s}^{-2} \right]$ & $\left[ \text{km} \right]$ & $\left[ \text{km} \right]$ & $\left[ \text{km} \right]$ & $\left[ \text{m} \right]$ \\ [0.5ex]
%          \hline
%          5040 & 9.81 & 20 & 50 & 100 & 30 \\ [1ex] 
%          \hline
%     \end{tabular}
%     \caption{Parameters of the Benchmark 2.}
%     \label{table:Setiyowati params}
% \end{table}

% One can observe that this sea floor function is not smooth, because its derivative is not continuous. In this type of situations, two options are possible: modeling the sea floor as it is and being aware of the consequent loss of convergence order of FC, or smoothing the function using hyperbolic tangents for example. For the present case the sea floor was modelled as it is, with its inner derivative discontinuities.
% \\Here the dry/wet interface has not been modelled, so the left boundary condition is computed as a wall interface. The right boundary is a wall boundary condition, however it is too far to modify the solution in the considered inner domain.
% \\The evolution of the water wave is shown \autoref{fig:Setiyowati}.

% \begin{figure}[!t]
% 	\centering
% 	\includegraphics[width=0.7\textwidth]{Setiyowati solution.pdf}
% 	\caption{Time-space diagram for the propagation of the tsunami wave.}
% 	\label{fig:1D spacetime}
% \end{figure}

\begin{figure}[!t]
	\centering
% 	\begin{subfigure}[b]{0.49\textwidth}
        % \centering
  	\includegraphics[width=.4\textwidth]{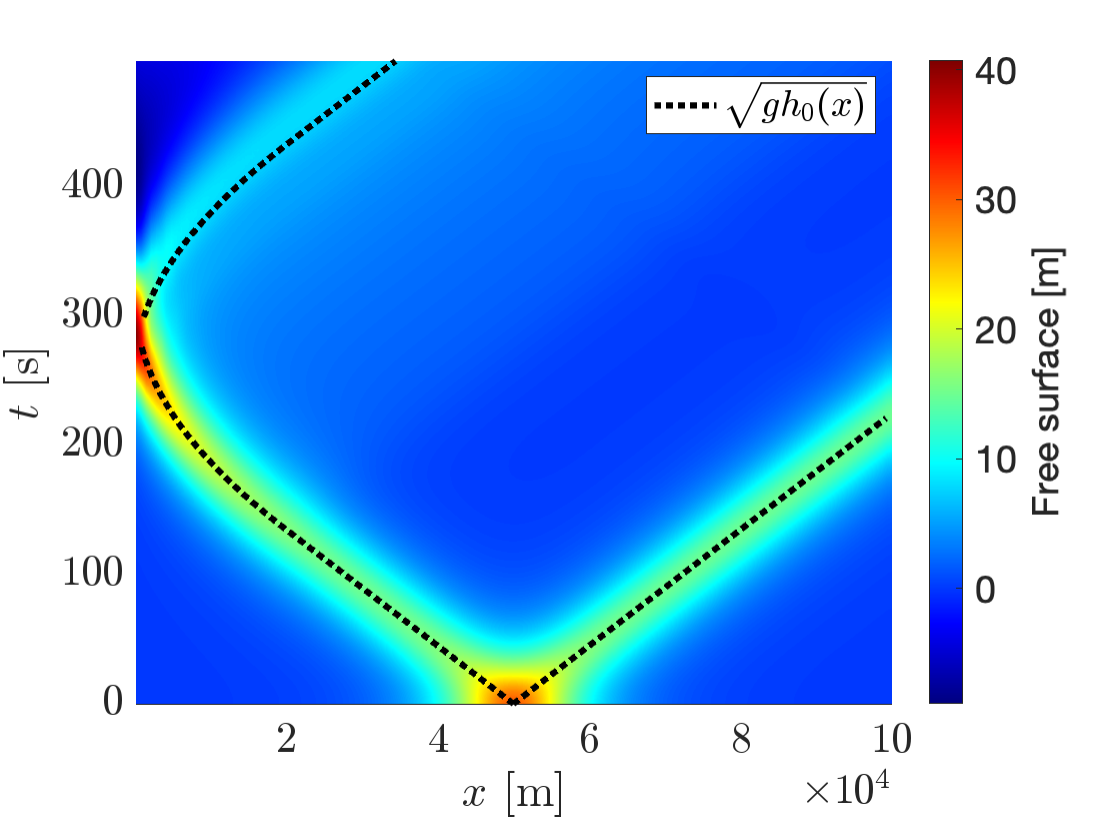}
        % \caption[eta]%
        % {{\footnotesize Space-time evolution.}}
        % \label{fig:1Dwavespeeda}
    % \end{subfigure}
%     \begin{subfigure}[b]{0.49\textwidth}
%         \centering 
%  	\includegraphics[width=\textwidth]{wavespeed.pdf}
%         \caption[]%
%         {{\footnotesize Wave speed as a function of distance from the left boundary.}}
%         \label{fig:1Dwavespeedb}
%     \end{subfigure}
	\caption{\emph{Benchmark 2}. The complete space-time evolution of the tsunami waves (left-going towards the sloping beach) produced by the FC-based solver. Overlaid in black is the expected approximate wave speed,  $c=\sqrt{gh_0(x)}$,  given by linear theory~\cite{Kundu2007}.}
	\label{fig:1Dwavespeeda}
\end{figure}

\autoref{fig:1Dwavespeeda} presents the complete space-time evolution of the left-going and right-going traveling waves in terms of a pseudocolor plot. The FC-based simulation agrees well with the reference study \cite{Setiyowati2019} in terms of propagation times and wave heights. Overlaid on the figure is the theoretical approximation of the speeds of the left-going and right-going waves (given by  linear wave theory~\cite{Kundu2007} as $c=\sqrt{gh_0}$ for still water depth $h_0$). Here, the FC-based solution tracks well with the expected linear theory, including with the expected slow down of the left-going wave as it approaches and climbs the sloping beach.

% The wave speed agrees well with the linear theory, near the shore, where nonlinear behavior are known to happen.

A more quantitative assessment relates to mass conservation. The area (or "volume") displaced at a time $t$ is given by
\begin{equation}
    V_\text{1D}(t)=\int_{x_{\text{min}}}^{x_{\text{max}}} \eta(x,t) dx.
\end{equation}
The corresponding percentage error from an initial time $t_\text{min} = 0 $ to a final time $t_\text{max}$ (when the left-going wave is entirely reflected) is hence given by 
\begin{equation}
    \varepsilon_{mc}=\left(\frac{V_\text{1D} \left(t_\text{max} \right)}{V_\text{1D} \left(t_\text{min} \right)}-1\right)\times 100.
    \label{eq:masserror}
\end{equation}
The FC-based solution using the same $\Delta x = 671$ m discretization (corresponding to $N_x = 150$ points per $100$ km) yields a mass conservation error $\varepsilon_{mc} = 0.058\%$. For reference, increasing the discretization to $N_x = 5\,000$ points per $100$ km reduces this error to a value of $\varepsilon_{mc} = 0.0000025 \%$, implying high-order convergence. In both cases (coarse and fine), mass conservation is effectively preserved numerically and converges. Indeed, such values fall significantly below the mass conservation tolerances as set forth in the standards of the National Oceanic and Atmospheric Administration (NOAA) \cite{synolakis2007standards}, which defines appropriate numerical conservation as when the final volume is within  5\% of the initial displaced volume. 

% total displaced volume at the end of the computation, i.e., when the initial wave is entirely reflected and offshore.'}

% m yields an error 
% The errors are summarized in \autoref{table:mass conservation}.
% \begin{table}[!t]
%     \centering
%     \begin{tabular}{|c |c|} 
%          \hline
%           $N_x \ \left[- \right]$ & $\varepsilon_{mc} \ \left[ \% \right]$ \\ [0.5ex]
%          \hline
%          $2\times 150$ & $5.8\times 10^{-2}$ \\ [1ex] 
%          \hline
%          $2\times 5000$ & $2.5\times 10^{-6}$ \\ [1ex] 
%          \hline
%     \end{tabular}
%     \caption{Mass conservation for the two discretizations.}
%     \label{table:mass conservation}
% \end{table}
% For each discretization, the standard criterion of mass conservation within 5\% error is respected, and refining leads to a reduction of the error which gives an indication of the convergence quality.
% \\One can observe that the sea floor function is not $\mathrm{C}^\infty$, which reduces the order of convergence of our solver in this case. However, if concerned about the high-order properties of the method, one can replace the slope function using hyperbolic functions to approximate the slope with $\mathrm{C}^\infty$ functions.

\subsection{Benchmark 3: 2D steady vortex}
\begin{figure*}[!t]
    \centering
     \begin{subfigure}[b]{0.45\textwidth}
        \centering 
        \includegraphics[width=.95\textwidth]{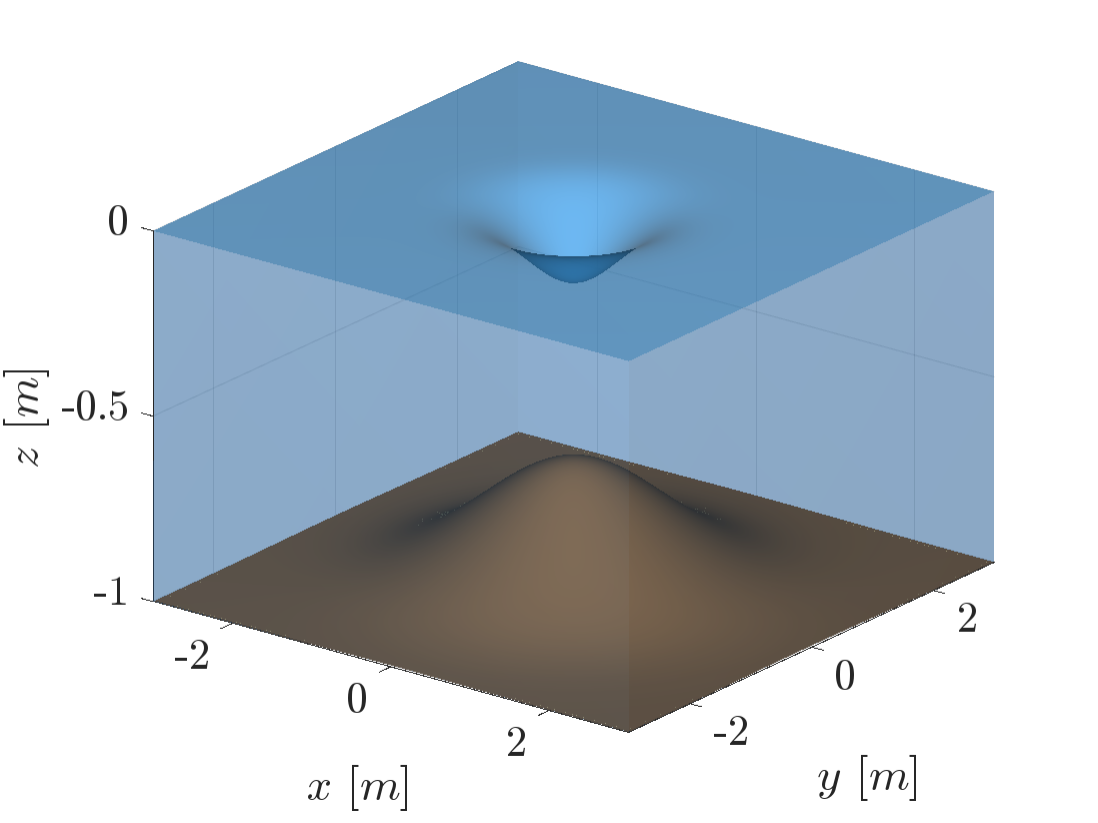}
        \caption[]%
        {{\footnotesize Water height and bathymetry;}}
        \label{first frame and bathy}
    \end{subfigure}
        \begin{subfigure}[b]{0.45\textwidth}
        \centering 
        \includegraphics[width=\textwidth]{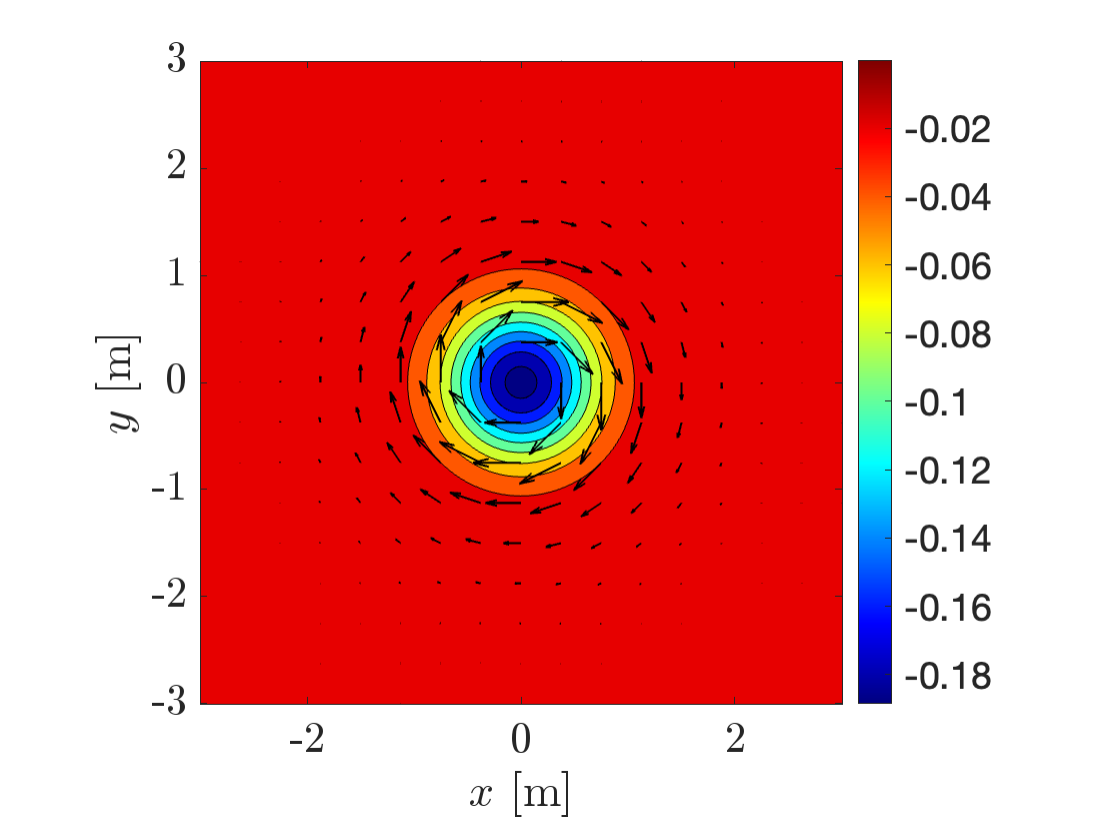}
        \caption[]%
        {{\footnotesize Steady free surface $\eta(x,y)$ and velocity field;}}
        \label{First eta frame}
    \end{subfigure}
    \begin{subfigure}[b]{0.45\textwidth}
        \centering 
        \includegraphics[width=\textwidth]{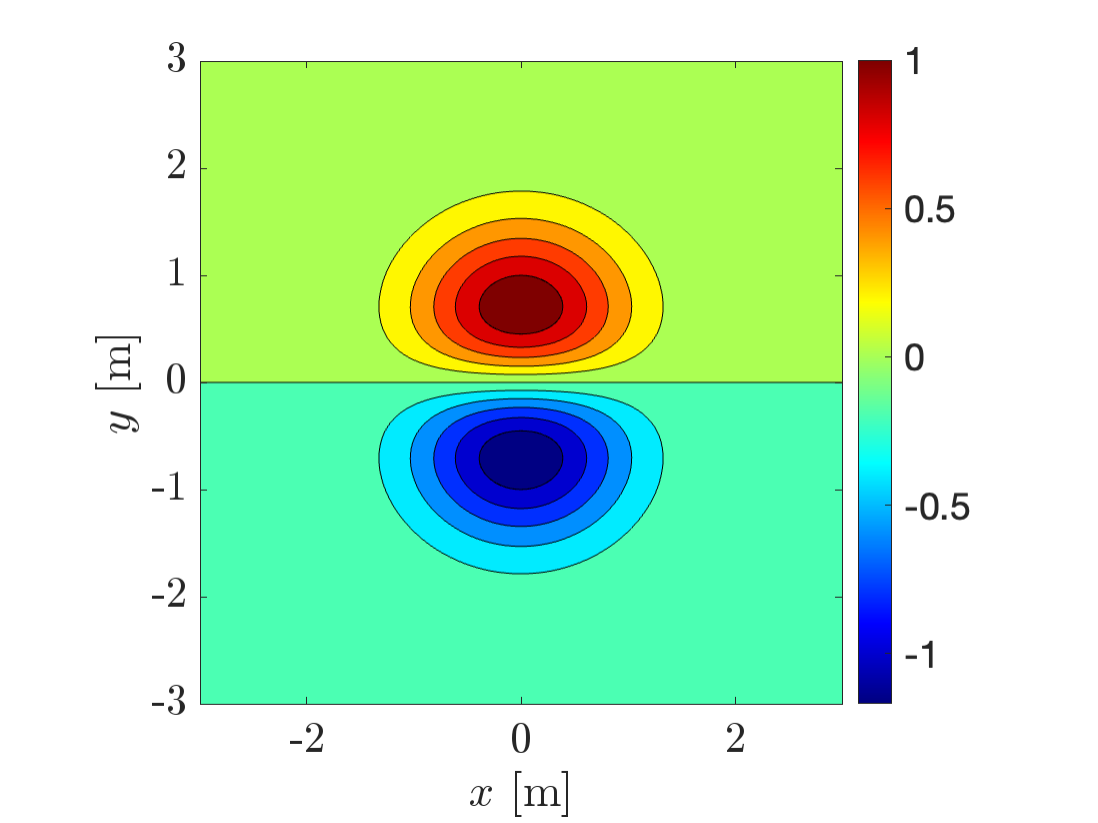}
        \caption[]%
        {{\footnotesize Steady velocity field $u(x,y)$};}
        \label{First u frame}
    \end{subfigure}
    \begin{subfigure}[b]{0.45\textwidth}
        \centering 
        \includegraphics[width=\textwidth]{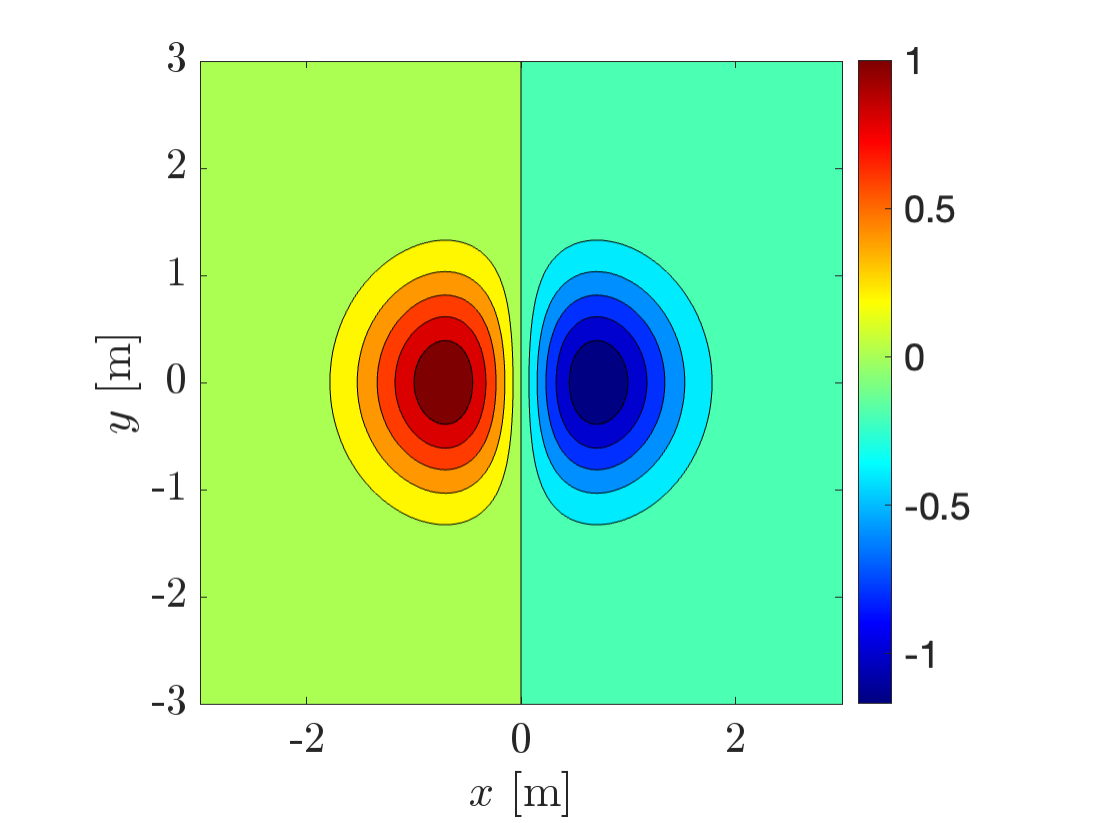}
        \caption[]%
        {{\footnotesize Steady velocity field $v(x,y)$.}}
        \label{First v frame}
    \end{subfigure}
    % \begin{subfigure}[b]{0.49\textwidth}
    %     \centering 
    %     \includegraphics[width=\textwidth]{vortex_velocity_frame.pdf}
    %     \caption[]%
    %     {{\footnotesize Steady velocity field}}
    %     \label{velocity}
    % \end{subfigure}
   
    \caption{\emph{Benchmark 3}. Illustration of the overall configuration from~\cite{michel2016} corresponding to the analytical solution given by~\autoref{eq:vortex}.}
    \label{fig:Vortex frames}
\end{figure*}

A third benchmark \cite{michel2016} concerns a 2D FC-based evaluation of the shallow water equations applied to a problem of a steady vortex and its comparison to a high-order finite volume (FV) treatment from \cite{michel2016}. Over a spatial domain $(x,y) \in \left[-3,3 \right]\times \left[-3,3 \right]$ and employing a gravitational constant $g=9.81$ m/s$^2$, the initial free surface and velocity field are given by
\begin{equation}
    \begin{cases}
        \eta_0(x,y)=-\frac{1}{4g}\exp(2(1-(x^2+y^2))),\\
        u_0(x,y)=y\exp(1-(x^2+y^2)),\\
        v_0(x,y)=-x\exp(1-(x^2+y^2)),
    \end{cases},
    \label{eq:vortex}
\end{equation}
where boundary conditions are chosen to be the same. The considered bathymetry is a Gaussian bump given by
\begin{equation}
    h_0(x,y)=-0.2\exp(0.5(1-(x^2+y^2))).
\end{equation}
The solution to such a problem is the (exact) steady solution given by
\begin{equation}
    \begin{cases}
        \eta(x,y,t)=\eta_0(x,y),\\
        u(x,y,t)=u_0(x,y),\\
        v(x,y,t)=v_0(x,y).
    \end{cases}
    \label{eq:vortexsolution}
    \end{equation}
\autoref{fig:Vortex frames} presents the overall problem setup alongside solution snapshots for the free-surface displacement $\eta(x,y)$ (\autoref{First eta frame}) and velocity field $(u(x,y)~v(x,y))^\text{T}$ (Figures~\ref{First u frame}-\ref{First v frame}) given by  \autoref{eq:vortexsolution}.

 \autoref{table:vortex cvg} presents the absolute $L^\infty$ errors (over all space and time up to $t_\text{max}=1$ s) for the free surface $\eta(x,y,t)$ of simulations produced by the FC-based solver applied to increasing spatial discretizations of size $N= N_x N_y$. For reference, the corresponding errors from~\cite{michel2016} employing FV are additionally reported using the very same discretizations. The high-order convergence of FC can be observed (up to seventh order), leading to errors that are orders-of-magnitude lower than the reference FV values from~\cite{michel2016}  (at the same fine discretization of $N=262\,144$, around 37 million times more accurate). Similarly low errors can be observed in \autoref{table:annex vortex cvg} for the velocity fields simulated by the FC solver (errors for such components are not provided in~\cite{michel2016}), which additionally presents the corresponding computation times in MATLAB (not provided for the reference FV solutions in~\cite{michel2016}).

% The errors have been calculated on the free surface at the same discretization points as the study by Michel-Donsac et al \cite{michel2016}, with $N=(N_x-1)(N_y-1)$ grid cells. \autoref{table:vortex cvg} gives the corresponding $L^\infty$ errors computed as $\max_{x,y,t} \left( \eta(x,y,t)-\eta_0(x,y) \right)$ and compare the solver errors to the highest accuracy version of the finite volume scheme given by Michel-Donsac et al \cite{michel2016}.

% The solver proves 7th order accuracy on this case, which is extremely efficient. The FC-solver errors are order of magnitude lower than the finite volume code. The solver's order of convergence decreases at the most converged cases, but remains more important than the finite volumes while proving an accuracy 30 million times finer on this last case.

\begin{table}[!t]
    \centering
    \begin{tabular}{|c |c |c | c  |c|} 
         \hline
          $N$ & FC $L^\infty$ error & FC order & FV $L^\infty$ error (from \cite{michel2016}) & FV order (from \cite{michel2016})  \\ [0.5ex]
         \hline
         1\,024 & $1.13\times 10^{-3}$ & - &  $5.52\times 10^{-2}$ & - \\ 
         \hline
         4\,096 & $9.87\times 10^{-6}$ & 6.99 & $1.68\times 10^{-2}$ & 1.72 \\
         \hline
         16\,384 & $7.75\times 10^{-8}$ & 7.07 &  $5.00\times 10^{-3}$ & 1.75 \\ 
         \hline
         65\,536 & $6.01\times 10^{-10}$ & 7.05 &  $1.40\times 10^{-3}$ & 1.84 \\ 
         \hline
         262\,144 & $1.02\times 10^{-11}$ & 5.90 &  $3.77\times 10^{-4}$ & 1.89 \\ 
         \hline
    \end{tabular}
    \caption{\emph{Benchmark 3.} Maximum errors and convergence rates of the free surface displacement $\eta(x,y,t)$ for the steady vortex problem treated by the FC-based solver and reference values corresponding to the highest order of finite volume (FV) employed in~\cite{michel2016}.}
    \label{table:vortex cvg}
\end{table}

\begin{table}[!t]
    \centering
    \small
    \begin{tabular}{|r  |c |c |r|} 
         \hline
          \multicolumn{1}{|c|}{$N$} &  FC $L^\infty$ $u$-error & FC $L^\infty$ $v$-error & \multicolumn{1}{c|}{Time [s]}  \\ 
         \hline
         1\,024 & $5.06\times 10^{-3}$ & $5.09\times 10^{-3}$ & 0.2 \\ 
         \hline
         4\,096 & $4.58\times 10^{-5}$ & $4.60\times 10^{-5}$ & 0.5 \\ 
         \hline
         16\,384 & $3.74\times 10^{-7}$ & $3.75\times 10^{-7}$ & 2.1 \\ 
         \hline
         65\,536 &  $1.26\times 10^{-8}$ & $1.26\times 10^{-8}$ & 14.0 \\ 
         \hline
         262\,144 & $2.07\times 10^{-9}$ & $2.07\times 10^{-9}$ & 102.1 \\ 
         \hline
    \end{tabular}
    \caption{\emph{Benchmark 3}. Maximum errors and computing times (MATLAB, see \autoref{remark:FCparams}), as a function of increasing spatial discretization, for the two velocity components $u(x,y,t)$ and $v(x,y,t)$ produced by the FC-based solver.}
    \label{table:annex vortex cvg}
\end{table}

\subsection{Benchmark 4: 2D tsunami propagation}
A fourth benchmark problem~\cite{Kanoglu2013} concerns the 2D propagation of an N-wave, which is a typical (realistic) shape of observed tsunami waves. Such a benchmark is used as a reference at the International Atomic Energy Agency~\cite{iaea2022} (IAEA) for evaluating the performance of academic tsunami solvers \cite{NAMIDANCE, Jagurs}. Over a domain $(x,y) \in [-100,100] \times [-100,100]$ (meters), the initial condition has the form of an N-wave given by
\begin{equation}
    \eta_0(x,y) = \frac{\varepsilon H}{2} \left[ \tanh \left(\gamma \left(x-x_0 \right) \right) - \tanh \left(\gamma \left(x-x_0-L \right) \right) \right] \left(y-y_2 \right) \text{sech}^2 \left(\gamma \left(y-y_1 \right) \right)
\end{equation}
for positions $x_0, y_1, y_2 \in \mathbb{R}$ and coefficients $\varepsilon, \gamma, H, L \in \mathbb{R}^+$. The problem admits an exact analytical solution to the \emph{linearized} shallow water equations~\cite{Kanoglu2013} that is given by
\begin{equation}
    \eta(x,y,t) = \frac{1}{(2\pi)^2} \int_{-\infty}^{\infty} \int_{-\infty}^{\infty} \hat{\eta}_0(k,l)e^{i(kx+ly)}\cos(\omega t) \text{ d}k \text{d}l
    \label{eq:linearsolution}
\end{equation}
for $\omega=\sqrt{k^2+l^2}$. Here, $\hat{\eta}_0$ is the Fourier transform of the initial condition and is given by~\cite{Kanoglu2013}
\begin{equation}
    \hat{\eta}_0(k,l) = i\frac{4\varepsilon H}{\pi} \alpha^3 \left(e^{-ikL}-1 \right)e^{-i(kx_0+ly_1)} \left[ \left(y_1-y_2 \right)l+i(1-\alpha l \text{coth}(\alpha l)) \right] \text{cosech}(\alpha k) \text{cosech}(\alpha l)
\end{equation}
for $\alpha = \pi/(2\gamma)$.

\begin{figure*}[!t]
    \centering
    \begin{subfigure}[b]{0.45\textwidth}
        \centering 
        \includegraphics[width=\textwidth]{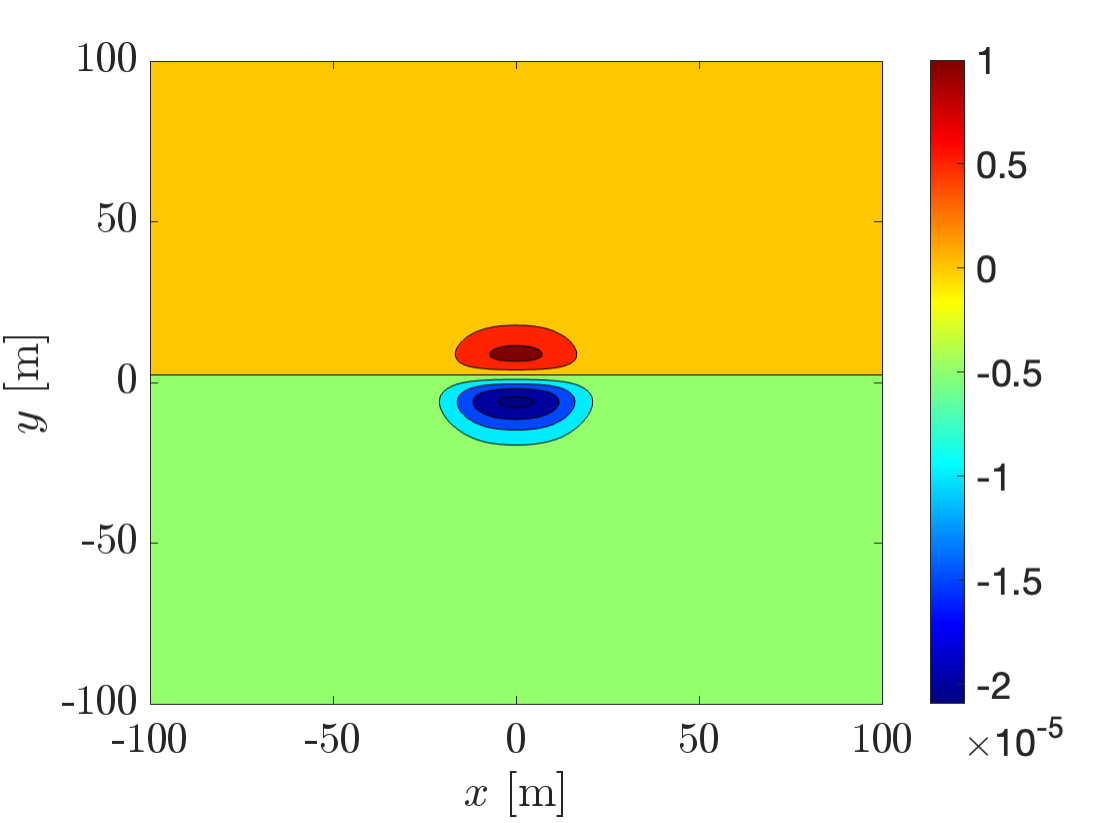}
        \caption[]%
        {{\footnotesize $t=0$ s;}}
        \label{First frame}
    \end{subfigure}\quad
    \begin{subfigure}[b]{0.45\textwidth}
        \centering 
        \includegraphics[width=\textwidth]{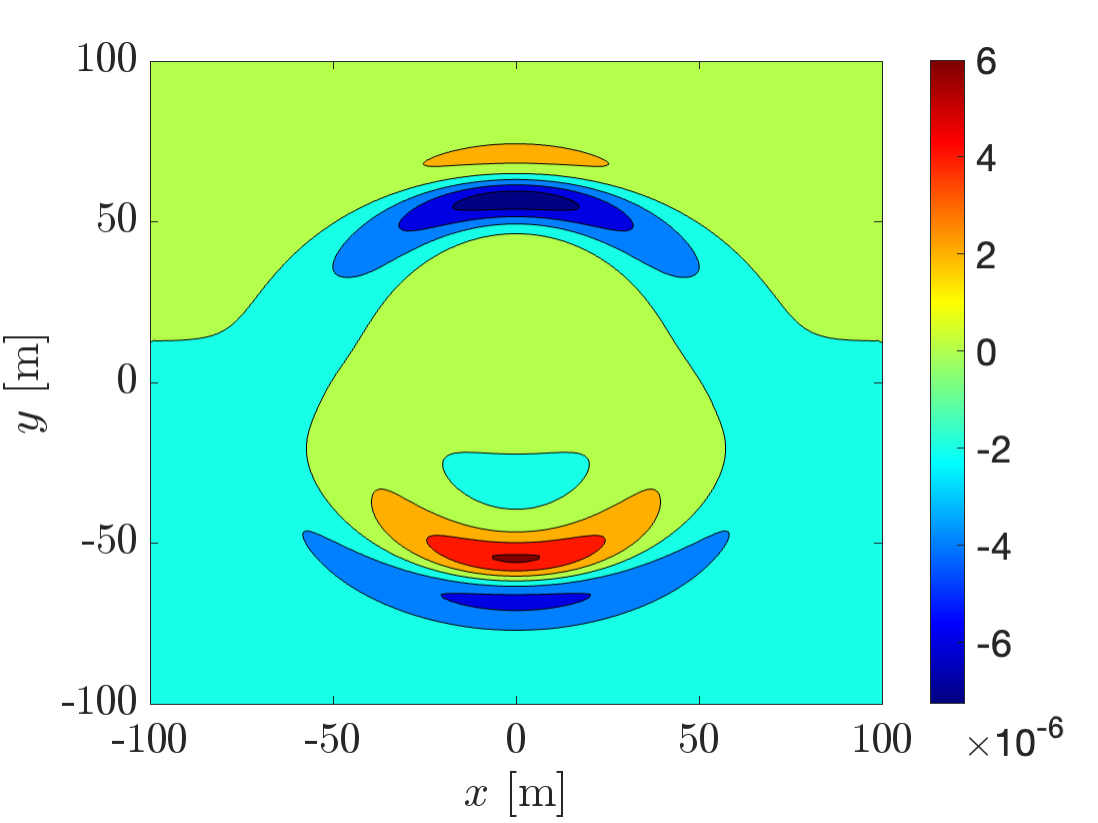}
        \caption[]%
        {{\footnotesize $t=60$ s.}}
        \label{Last frame}
    \end{subfigure}
    \caption{\emph{Benchmark 4.} Snapshots of initial condition ($t=0$ s) and last timestep ($t=60$ s) produced by the FC-based solver for the N-wave configuration.}
    \label{fig:Kanoglu contours}
\end{figure*}

\autoref{fig:Kanoglu contours} presents the initial and final temporal snapshots in 2D of the solutions obtained by the FC-solver for a simulation conducted on a discretization of size $N_x=N_y=201$ and advanced up to a time $t=60$ s. The boundary condition are set as radiation (see \autoref{eq:2Dradiationbc}).  All parameters of the initial condition are adopted from \cite{iaea2022, NAMIDANCE, Jagurs} and correspond to $h_0 =1$ m, $g = 1$ m/s$^2$, $H=0.001$ m, $L=30$ m, $y_1=0$ m, $y_2 = 2.3$ m, $x_0 = -15$ m, $\epsilon = 0.004$, and $\gamma=0.1061$ (unitless).

\begin{figure}[!t]
	\centering
	\includegraphics[width=0.5\textwidth]{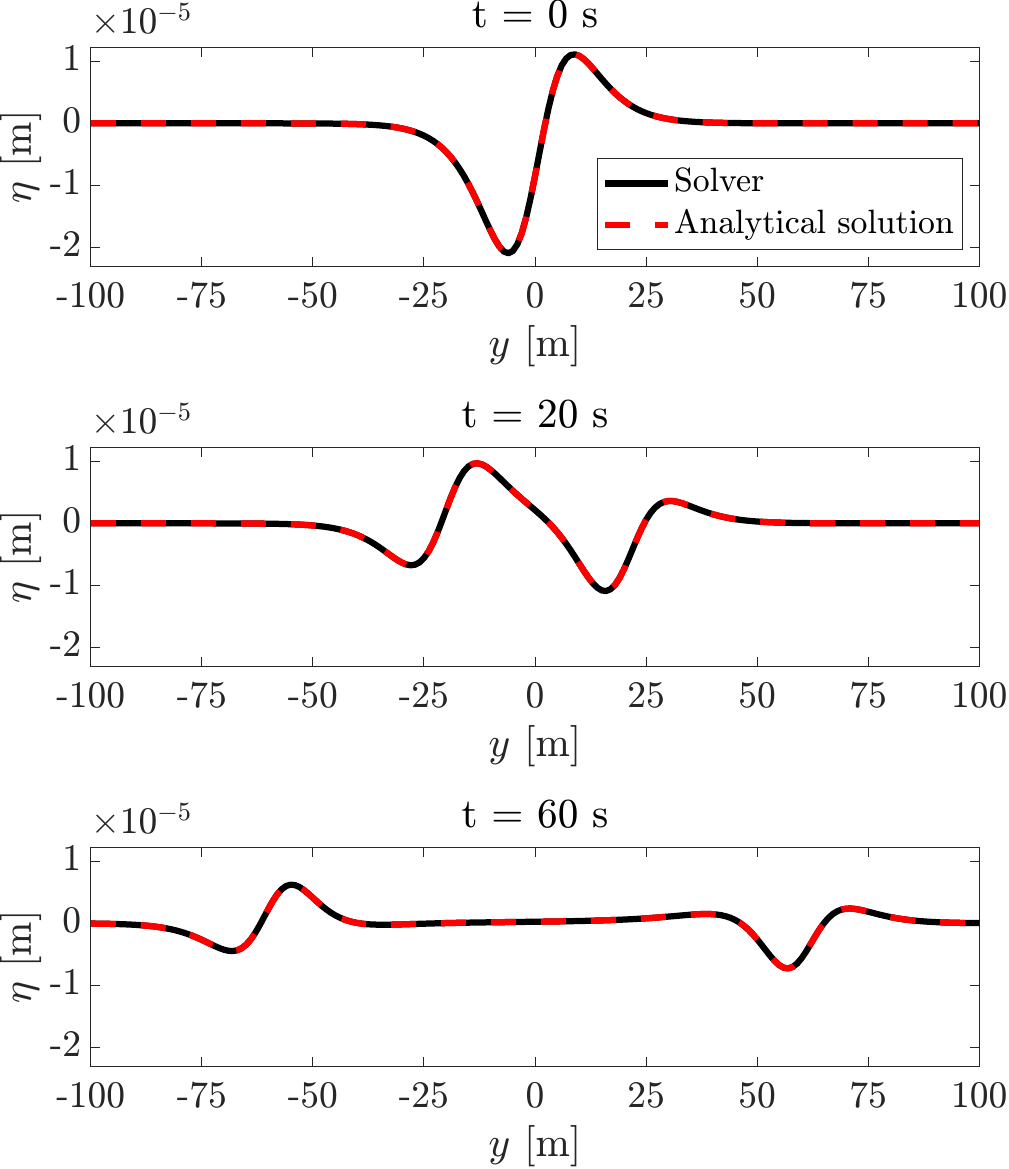}
	\caption{\emph{Benchmark 4.} Snapshots at various times of the free surface displacement $\eta$ produced by the FC-based solver (black) on the line corresponding to $x=0$. Overlaid in red is the analytical solution from \cite{Kanoglu2013}.}
	\label{fig:Kanoglu curves}
\end{figure}

\autoref{fig:Kanoglu curves} additionally presents temporal snapshots slice-wise at $x=0$ m for times $t=0, 20$ and $60$ s (corresponding to the same as those considered in \cite{NAMIDANCE, Jagurs}).  The analytical (linearized) solution given by \autoref{eq:linearsolution} is overlaid in red, demonstrating excellent qualitative agreement with the solutions produced by the FC solver in black. The corresponding maximum $L^\infty$ errors over all time and space, given by (in percentage)
\begin{equation}\label{eq:error}
    \varepsilon_\%(t)=\frac{\max_{x,y} \left(|\eta(x,y,t)-\eta_{\text{ref}}(x,y,t)| \right)}{\max_{x,y,t'} \left(|\eta_{\text{ref}}(x,y,t')| \right)}\times 100,
\end{equation}
are found to be $\varepsilon_\% = 5.62 \times 10^{-3}~\%$ and $\varepsilon_\% = 2.82 \times 10^{-2}~\%$ up to time  $t=20$ s and $t= 60$ s, respectively. It should be noted that these errors are computed with respect to an exact solution of the \emph{linearized} shallow water equations~\cite{Kanoglu2013}, and hence cannot be reasonably expected to converge using the  nonlinear shallow water equations that are treated in this work. In addition, the simplified non-reflecting boundary condition limits the order of convergence given by FC at the boundary. Nevertheless, the errors are below $0.03\%$ over all space and all time.

% The location of these snapshots 

% and \autoref{fig:Kanoglu curves}

% Assessing the performance of the FC solver more quantitatively, 
%  Therefore, we calculated the free surface at these three times. The parameters used for this benchmark are summarized in \autoref{table:kanoglu params}. The other constants $\gamma$ and $\alpha$ are deduced from the parameters as $\gamma=\sqrt{3HP_0/4}$ and $\alpha=\pi/(2\gamma)$. The number of points in $x$ and in $y$ is $N_x=N_y=201$. $\gamma=.1061$
% \begin{table}[!t]
%     \centering
%     \begin{tabular}{|c |c |c |c |c |c |c |c |c |c |c |c |c|} 
%          \hline
%           $h_0$ & $g$ & $H$ & $L$ & $y_1$ & $y_2$ & $x_0$ & $\varepsilon$ & $P_0$ & $x_{\text{min}}$ & $x_{\text{max}}$ & $y_{\text{min}}$ & $y_{\text{max}}$\\ 
%           $[\text{m}]$ & $[\text{m.s}^{-2}]$ & $[\text{m}]$ & $[\text{m}]$ & $[\text{m}]$ & $[\text{m}]$ & $[\text{m}]$ & $[-]$ & $[\text{m}^{-1}]$ & $[\text{m}]$ & $[\text{m}]$ & $[\text{m}]$ & $[\text{m}]$ \\ [0.5ex]
%          \hline
%          1 & 1 & 0.001 & 30 & 0 & 2.3 & -15 & 0.004 & 15 & -100 & 100 & -100 & 100 \\ [1ex] 
%          \hline
%     \end{tabular}
%     \caption{Parameters of Benchmark 4.}
%     \label{table:kanoglu params}
% \end{table}

% These results are calculated at the same timesteps as in the IAEA study, on the slice $y=0$. The corresponding curves are shown \autoref{fig:Kanoglu curves}. The mass conservation is respected within an error $\varepsilon_{mc}=4.8\times 10^{-11}$ \%. The solver results capture very well the analytical solution of this problem. To give a quantitative insight of the accuracy, we summarized the errors in the \autoref{Errors for Benchmark 4}. The errors are calculated as:
% \begin{equation}\label{eq:error}
%     \varepsilon_\%(t)=\frac{\max_{x,y} \left(|\eta(x,y,t)-\eta_{\text{ref}}(x,y,t)| \right)}{\max_{x,y,t'} \left(|\eta_{\text{ref}}(x,y,t')| \right)}\times 100
% \end{equation}
% The maximal error is of the order of 0.01\% which is a satisfactory accuracy.

% The errors of this computation also include the fact that the theoretical model was derived form the linear shallow water system, whereas the solver is solving the nonlinear version of these equations.
% \begin{table}[!t]
% \centering
% \begin{tabular}{|c |c |c|} 
%  \hline
%   & $\varepsilon_\% \ (t=20 \ \text{s})$ & $\varepsilon_\% \ (t=60 \ \text{s})$ \\ [0.5ex] 
%  \hline
%  Solver & $5.56\times 10^{-3}$ & $1.35\times 10^{-2}$ \\ [1ex] 
%  \hline
% \end{tabular}
% \caption{Errors of the solver on the 2D analytical benchmark. The percentage errors are computed as in \autoref{eq:error}.}
% \label{Errors for Benchmark 4}
% \end{table}

\subsection{Benchmark 5: waves generated by dynamic sea floor displacement} \label{piston benchmark}
\begin{figure}[!t]
    \centering
    % \begin{subfigure}[b]{0.49\textwidth}
    %     \centering 
    %     \includegraphics[width=\textwidth]{hammack_original.pdf}
    %     \caption[]%
    %     {{\footnotesize Experimental setup from Hammack (1973) \cite{Hammack1973}.}}
    %     \label{fig:Hammack original}
    % \end{subfigure}
    % \hfill
    % \begin{subfigure}[b]{0.49\textwidth}
        % \centering 
        % \vspace{15mm} % Adjust this value to vertically center the image
        \includegraphics[width=.425\textwidth]{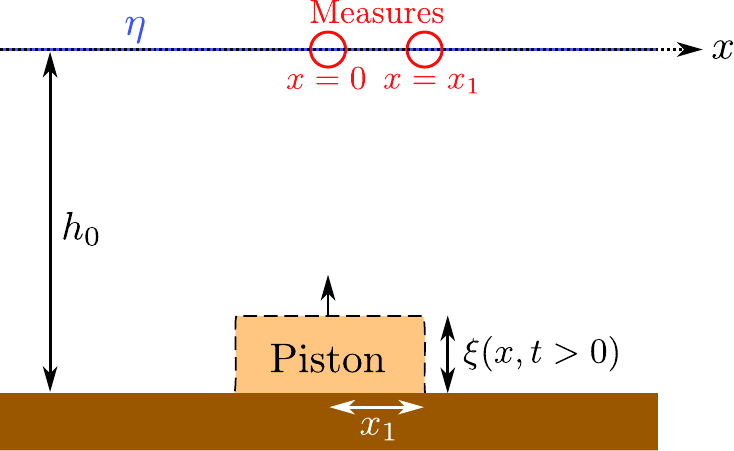}
        % \caption[]%
        % {%
        %     \centering
        %     {\small Sketch of the benchmark problem inspired by Hammack \cite{Hammack1973}.}
        % }
        % \label{fig:Hammack}
    % \end{subfigure}
    \caption{\emph{Benchmark 5}. Illustration of the moving piston configuration inspired by \cite{Hammack1973}.}
    \label{fig:Hammack}
\end{figure}

A final benchmark concerns the simulation of waves generated by fully dynamic motions of the seafloor (bathymetry movement) and their assessment versus experimental data. Such a configuration has been treated by other solvers \cite{Elbanna2021, Derakhti2019}, and is inspired by the pioneering work of Hammack, et al. \cite{Hammack1973}, which presented laboratory measurements of the water elevation caused by a moving source (a piston vertically displacing an interval of the seafloor). An illustration of the problem formulation is presented \autoref{fig:Hammack} (the reader is referred to \cite{Hammack1973} for the corresponding experimental setup). 

The unperturbed bathymetry $h(x) = h_0$ is considered initially flat and is then subjected to a dynamic perturbation that is given by
\begin{equation}
    \xi(x,t)=f(t)\mathcal{H} \left( x^2-x_1^2 \right),
\end{equation}
where $x_1 \in \mathbb{R}$ is the half-length of the piston and $\mathcal{H}$ is the Heaviside function given by
\begin{equation}
    \mathcal{H}(x)=
    \begin{cases}
    1, & x\geq0 \\ 0, & x<0.
    \end{cases}
\end{equation} 
The temporal function $f(t)$ is given by
\begin{equation}
    f(t)=
    \begin{cases}
    0, & t<0\\ \xi_m \left(1-\cos \left(\pi t/t_r\right)\right)/2, & 0\leq t \leq t_r
    \\ \xi_m, & t>t_r,
    \end{cases}
\end{equation}
for perturbation amplitude $\xi_m \in \mathbb{R}^+$ (i.e., the maximum piston height) and rising time $t_r \in \mathbb{R}^+$. The latter is given by a non-dimensionalized characteristic ratio $\tau = t_r/t_c \ll 1$ between the rise time $t_r$ and the shallow water wave propagation time $t_c = x_1/\sqrt{gh_0}$ for characteristic length $x_1$ (the half-length of the piston). The initial conditions correspond to water at rest, i.e. $\eta(x,0)=0$, $u(x,0)=0$, and the boundary conditions are considered to be walls, i.e., $\eta_x=0$ and $u=0$ at the boundary (see \autoref{eq:1Dwallbc}).

\begin{remark}
In order to retain the high-order nature of the FC-based solver introduced in this work, a smoothing of the rectangular perturbation (piston) is employed by approximating the Heaviside function $\mathcal{H}$ by means of hyperbolic tangents, i.e.,
\begin{equation*}
    \mathcal{H} \left( x^2 - x_1^2 \right) \approx \frac{1}{2} \left(\text{tanh} \left(a \left(x+x_1 \right) \right)-\text{tanh} \left(a \left(x-x_1 \right) \right) \right),
\end{equation*}
with $a \ll 1$ (a value of $a=0.01$ m$^{-1}$ is employed in this section). 
\end{remark}

\begin{figure*}[!t]
    \centering
    \captionsetup[subfigure]{oneside,margin={0.6cm,0cm}}
    \begin{subfigure}[b]{0.325\textwidth}
        \centering 
        \includegraphics[width=\textwidth]{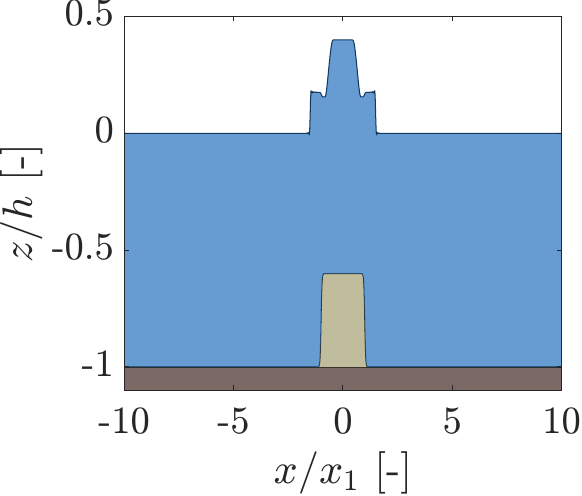}
        \caption[]%
        {{\footnotesize FC solver, $t \sqrt{g/h}=6.0$;}}
        \label{FC frame 1}
    \end{subfigure}
    \begin{subfigure}[b]{0.325\textwidth}
        \centering
        \includegraphics[width=\textwidth]{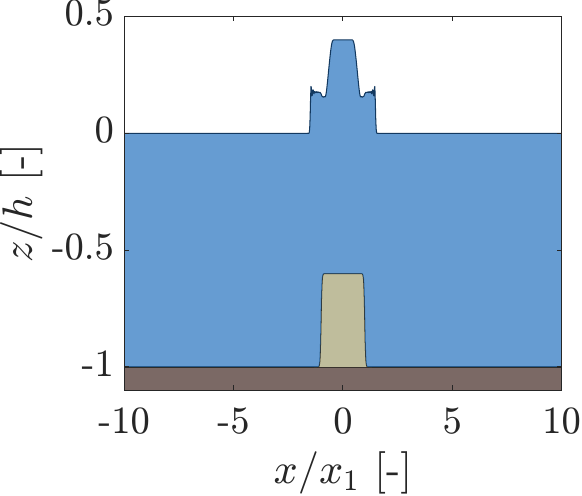}
        \caption[]%
        {{\footnotesize FD4, $t \sqrt{g/h}=6.0$;}}
        \label{FD4 frame 1}
    \end{subfigure}
    \begin{subfigure}[b]{0.325\textwidth}
        \centering 
        \includegraphics[width=\textwidth]{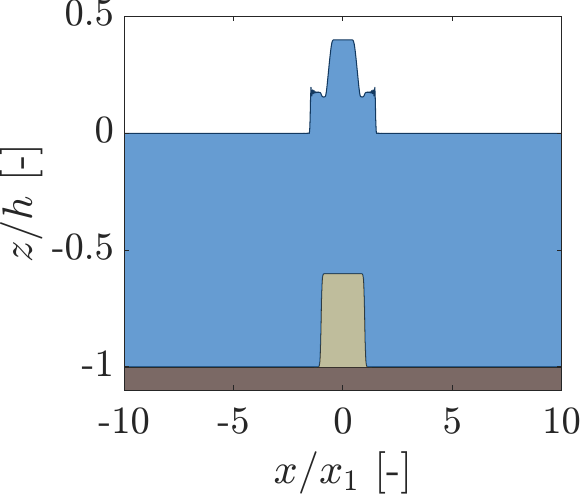}
        \caption[]%
        {{\footnotesize FD6, $t \sqrt{g/h}=6.0$;}}
        \label{FD6 frame 1}
    \end{subfigure}
    \vskip\baselineskip
    \begin{subfigure}[b]{0.325\textwidth}
        \centering 
        \includegraphics[width=\textwidth]{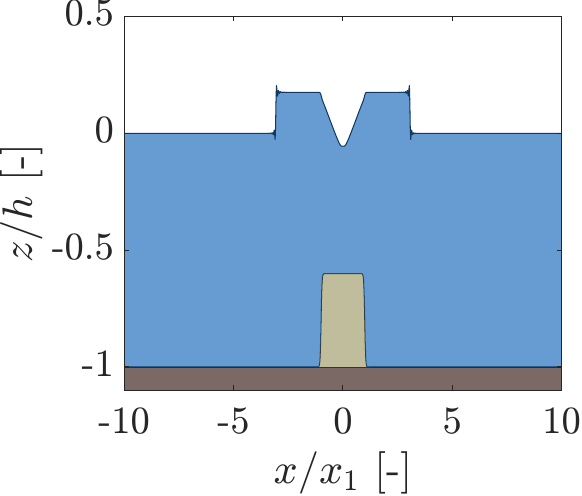}
        \caption[]%
        {{\footnotesize FC solver, $t \sqrt{g/h}=23.0$;}}
        \label{FC frame 2}
    \end{subfigure}
    \begin{subfigure}[b]{0.325\textwidth}
        \centering
        \includegraphics[width=\textwidth]{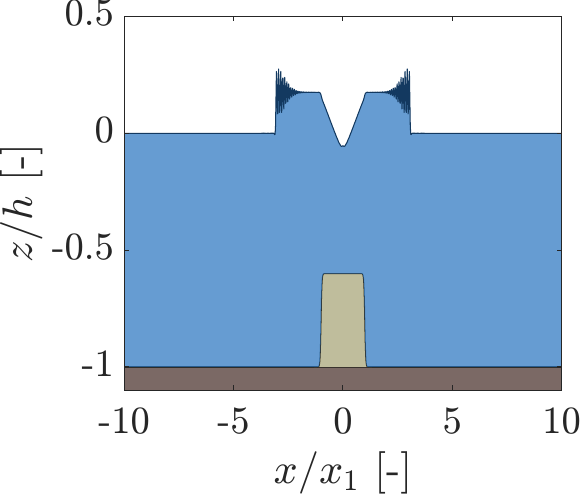}
        \caption[]%
        {{\footnotesize FD4, $t \sqrt{g/h}=23.0$;}}
        \label{FD4 frame 2}
    \end{subfigure}
    \begin{subfigure}[b]{0.325\textwidth}
        \centering 
        \includegraphics[width=\textwidth]{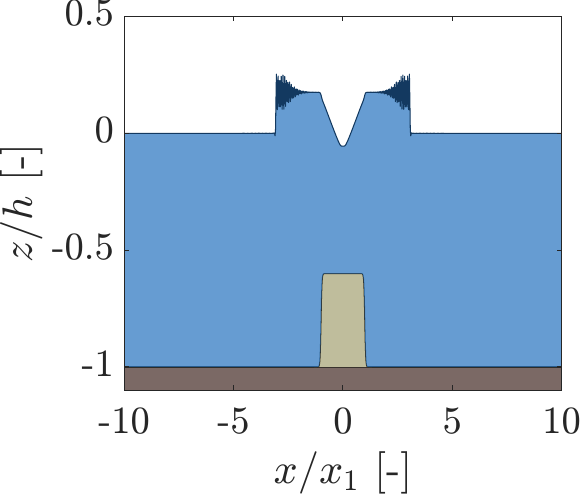}
        \caption[]%
        {{\footnotesize FD6, $t \sqrt{g/h}=23.0$;}}
        \label{FD6 frame 2}
    \end{subfigure}
    \vskip\baselineskip
    \begin{subfigure}[b]{0.325\textwidth}
        \centering 
        \includegraphics[width=\textwidth]{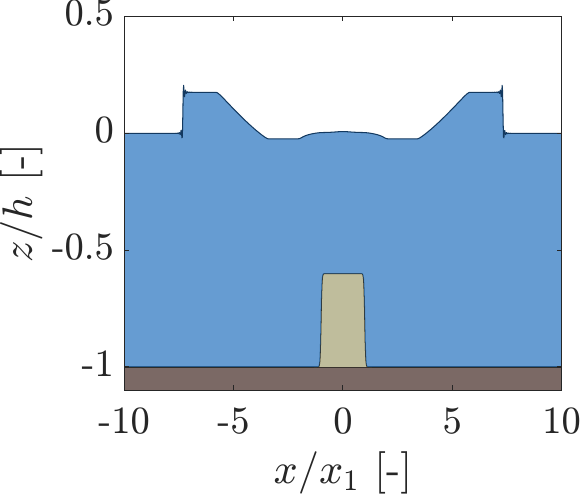}
        \caption[]%
        {{\footnotesize FC solver, $t \sqrt{g/h}=69.0$;}}
        \label{FC frame 3}
    \end{subfigure}
    \begin{subfigure}[b]{0.325\textwidth}
        \centering
        \includegraphics[width=\textwidth]{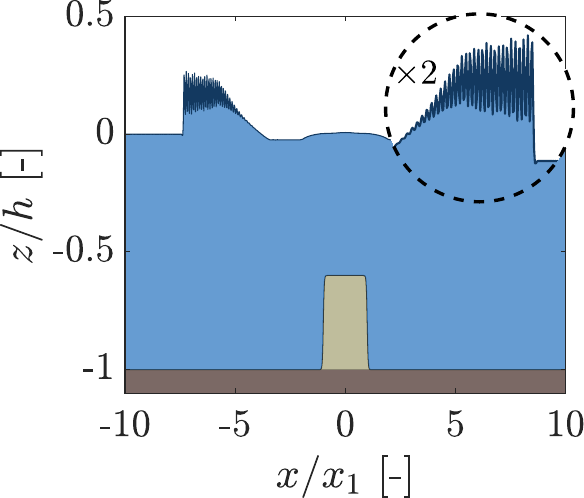}
        \caption[]%
        {{\footnotesize FD4, $t \sqrt{g/h}=69.0$;}}
        \label{FD4 frame 3}
    \end{subfigure}
    \begin{subfigure}[b]{0.325\textwidth}
        \centering 
        \includegraphics[width=\textwidth]{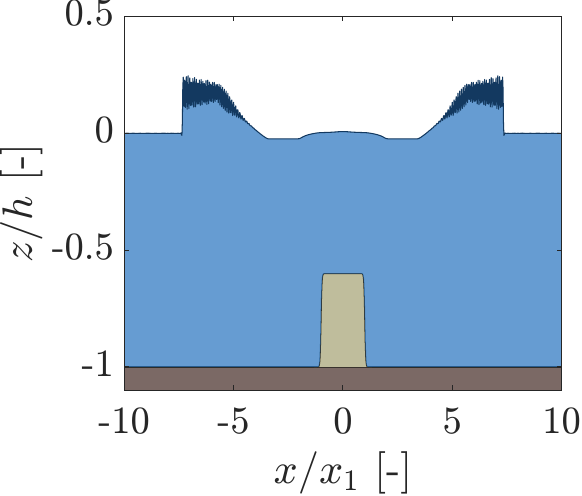}
        \caption[]%
        {{\footnotesize FD6, $t \sqrt{g/h}=69.0$.}}
        \label{FD6 frame 3}
    \end{subfigure}
    \caption{\emph{Benchmark 5}. Snapshots at various times (each row) of the free surface displacement $\eta$ produced by FC (left column), FD4 (middle column), and FD6 (right column). All simulations employ the same discretization of size $N_x=1001$. The dynamic piston movement is already is already at maximum displacement at these times (yellow). No amplification scaling of the solution is presented here for visualization.}
    \label{fig:Derakhti curves}
\end{figure*}

Employing $h_0 = 200$ m, $g = 9.81$ m/s$^2$, and non-dimensionalized values $x_1/h_0 = 12.2$ and $\xi_m/h = 0.4$, \autoref{fig:Derakhti curves} presents temporal snapshots of the solutions produced by FC (left column), FD4 (middle column), and FD6 (right column). All simulations are conducted over a (non-dimensionalized) domain $x/x_1 \in [-10,10]$ up to a final (non-dimensionalized) time $t = t_\text{max}\sqrt{g/h} = 69$, using the same discretization $N_x = 1001$ points and the maximum allowable  timestep by each of the method's respective CFL conditions (\autoref{remark:FCparams}). After being generated by the seafloor displacement caused by the moving piston, two waves propagate in each direction.  Without any special treatment of the discontinuity or shock, the FC-based solution incurs very reasonable errors despite the sharp wave front. A slight error (that does not grow in time)  can be observed at the water jump due to the lack of the highest-frequency Fourier modes in the discrete Fourier series. For finite differences, on the other hand, one can observe that the waves experience very large numerical oscillations at the sharp wave front (even though the stencils are of high order). Such high-frequency errors for FD4 and FD6 continue to grow in time as the waves propagate (\autoref{FD4 frame 3} and \autoref{FD6 frame 3}), illustrating that, in the case of such rapid sea floor motions, robust methods like FC may be better suited. Even with substantial refinements, such oscillations are not assured to be mitigated (see also the benchmark in \autoref{sec:benchmark1}). 

\begin{figure*}[!t]
    \centering
    \captionsetup[subfigure]{oneside,margin={0.6cm,0cm}}
    \begin{subfigure}[b]{0.49\textwidth}
        \centering 
        \includegraphics[width=\textwidth]{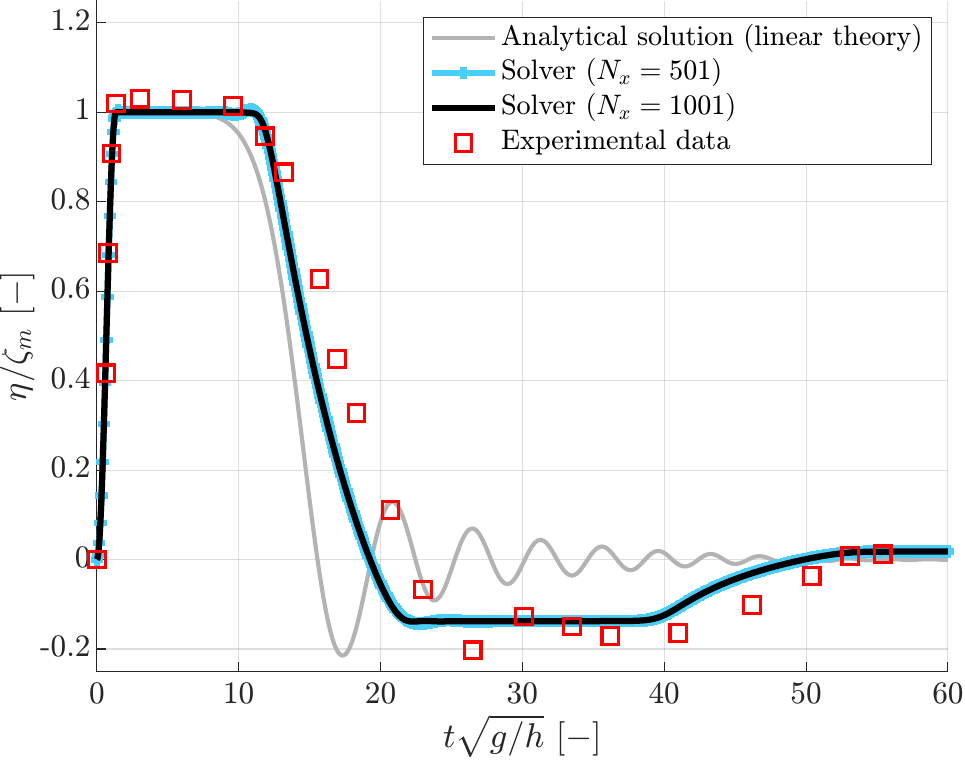}
        \caption[]%
        {{\footnotesize Piston center ($x=0$);}}
        \label{Piston center}
    \end{subfigure}\hfill
    % \vskip\baselineskip
    \begin{subfigure}[b]{0.49\textwidth}
        \centering 
        \includegraphics[width=\textwidth]{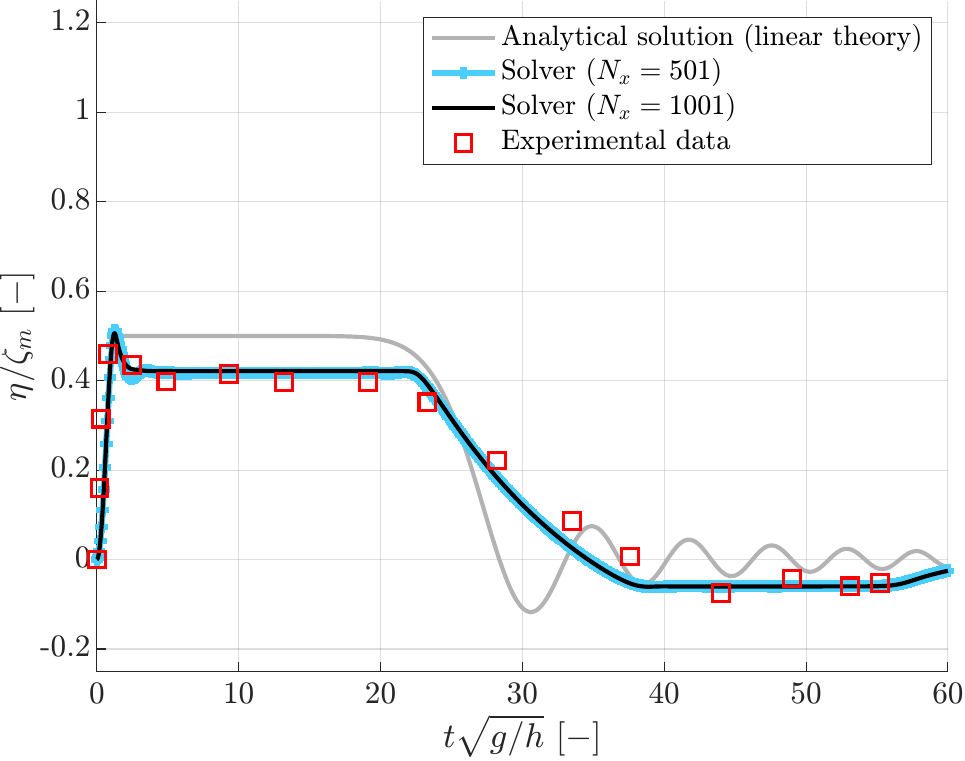}
        \caption[]%
        {{\footnotesize Piston edge ($x=x_1$).}}
        \label{Piston edge}
    \end{subfigure}
    \caption{\emph{Benchmark 5}. Time evolution of the normalized free surface displacement $\eta$ at the center (left) and edge (right) of the piston, produced by the FC-based solver. Two discretizations are presented (in blue and black), demonstrating convergence of the solver and good correspondence to the experimental results (red squares) from \cite{Hammack1973}. The linearized analytical solution from \cite{Hammack1973} is additionally presented and demonstrates its limitations in comparison with the data.}
    \label{fig:Piston}
\end{figure*}

\autoref{fig:Piston} presents the evolution of the free surface at the middle (left figure) and edge (right figure) of the piston, where available laboratory data are measured~\cite{Hammack1973} and are presented alongside in red. Since dimensions of the physical experiments differ from the proposed simulated configurations~\cite{Derakhti2019} by orders of magnitudes (initial bathymetry depth here is $h_0 = 200$ m, whereas in the laboratory depths ranged between 5 cm to 50 cm), all distances and times are presented in non-dimensionalized form. The results of the solver agree well globally with all observed measurements at both locations, with a maximum local error ($L^\infty$ over time) of 32\% and 20\%, respectively. For reference, an approximate analytical solution is additionally presented in grey, derived from linear wave theory in the original work~\cite{Hammack1973} and given by
\begin{equation}
    \eta(x,t)=\frac{\xi_m}{\pi}\int_0^{\infty} \frac{\cos(kx)\sin(kb)}{k\cosh(kh)} \frac{\kappa^2}{\kappa^2-\omega^2} \left[\cos(\omega t)-\cos(\kappa t)+ \mathcal{H}\left(t-t_r \right) \left(\cos \left(\omega \left(t-t_r \right)\right))+\cos(\kappa t)\right)\right]dk
\end{equation}
for $\kappa=\pi/t_r$ and $\omega=\sqrt{gk\tanh(kh)}$. Such a solution evolves like a damped oscillatory system and is unable to adequately capture the (nonlinear) dynamic response. Indeed, this linear theory matches experimental results only for $\xi_m<<h$~\cite{Derakhti2019}, which is not the case here given that $\xi_m/h=0.4$.

% Even without such smoothing, use of the original sharp Heaviside function has been observed to still generate relatively accurate solutions with the FC solver (although with very small spurious errors).

% \section{A case study of tsunamigenesis induced by dynamic earthquake ground motion}\label{sec:Application}
\section{A parametric case study on the effects of earthquake speeds on tsunamigenesis}\label{sec:Application}
The performance and physical accuracy of the proposed FC-based solver has been demonstrated through the rigorous numerical experiments of \autoref{Numerical study} and the benchmarks of \autoref{Benchmarks}. This section applies the solver on a realistic and understudied geophysical scenario in an effort towards further study of the dynamic effects of earthquakes (and their subsequent time-dependent velocities and displacements) on tsunamigenesis. Such a problem, introduced herein, is inspired by recent results that have suggested that some induced tsunami dynamics are linked to the corresponding time-dependencies of the earthquake sources~\cite{Elbanna2021,amlanibhat2022, dias2007dynamics}. The objective is to propose a novel problem configuration in order to investigate the influence of \emph{earthquake speed} on tsunami generation, as well as to compare it to classical tsunami approaches that neglect time-dependent seafloor movements. To this end, two variations are presented as follows: a 1D formulation and a 2D formulation, both utilizing physically-realistic parameters inspired by an idealized subduction fault (perpendicular to the sea surface).

\subsection{1D formulation}\label{sec:1D application}
\begin{figure*}[!t]
    \centering
    \captionsetup[subfigure]{oneside,margin={0.6cm,0cm}}
    \begin{subfigure}[b]{0.495\textwidth}
        \centering 
        \includegraphics[width=.8\textwidth]{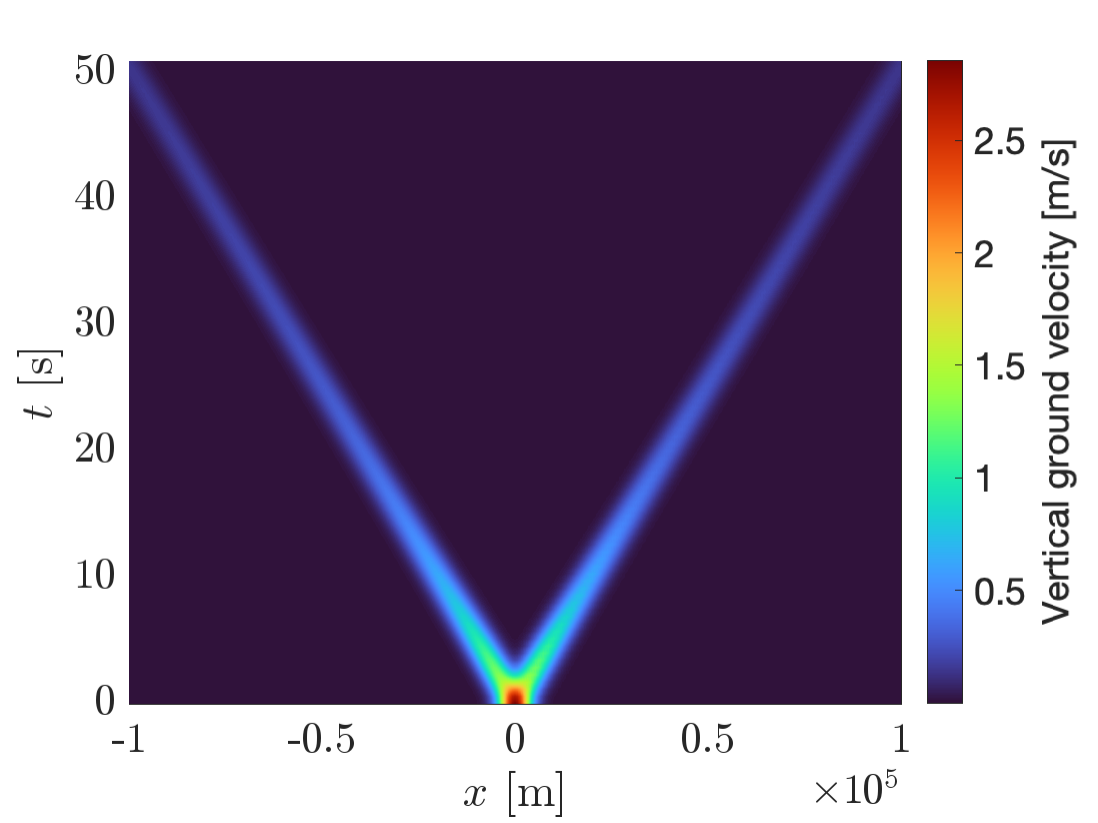}
        % \caption[]%
        % {{\footnotesize \red{Dynamic earthquake source velocity $\xi_t(x,t)$;}}}
        \label{fig:ksit surface waves}
    \end{subfigure}
    \begin{subfigure}[b]{0.495\textwidth}
        \centering 
        \includegraphics[width=.8\textwidth]{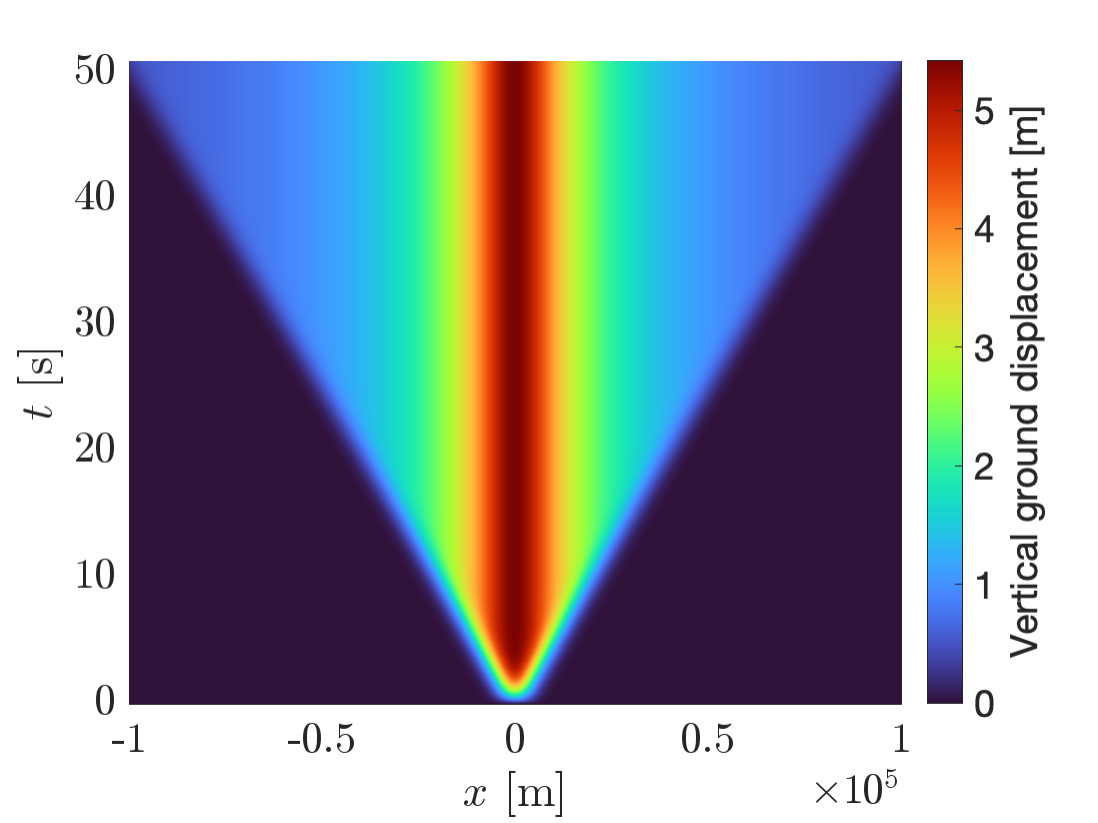}
        % \caption[]%
        % {{\footnotesize Corresponding total dynamic displacement of the sea floor $\xi(x,t)$.}}
        \label{fig:ksi surface waves}
    \end{subfigure}
    \caption{\emph{1D earthquake-induced tsunamigenesis}. The dynamic source for an earthquake speed of $c_\text{eqk}=10 \sqrt{gH_0}$, in terms of velocity $\xi_t(x,t)$ (left) and displacement $\xi(x,t)$ (right), corresponding to the proposed problem configuration defined by \autoref{eq:ksi surface waves} and \autoref{eq:dksidt surface waves}.}
    \label{fig:surface wave source}
\end{figure*}

The fundamental novelty of the overall problem lies in the analysis of the effects of earthquake speeds on tsunamigenesis by considering sources mimicking ruptures, perpendicular to the sea floor, that propagate seismic waves and ultimately generate surface (Rayleigh) waves of parameterized speeds given by
\begin{equation}\label{eq:ceqk}
    c_\text{eqk} = n \sqrt{g H_0}
\end{equation} 
for integers $n \in \mathbb{N}^*$ (i.e., integer multiples of the characteristic tsunami propagation speed $\sqrt{gH_0}$). Such sources can be modeled as time-dependent sea floor displacements $\xi(x,t)$ in the shallow water equations that are given by 
\begin{equation}\label{eq:ksi surface waves}
    \xi(x,t)=\int_{0}^{t} \xi_t(x,\tau) \text{d}\tau,
\end{equation}
where $\xi_t$ is the time-derivative of the vertical displacement (i.e., vertical seafloor velocity) given by
\begin{equation}\label{eq:dksidt surface waves}
    \xi_t(x,t)=\frac{nA}{\sqrt{(x/L)^2 + \varepsilon^2}} \left[ \exp \left(-\frac{(x - c_\text{eqk}t)^2} {2\sigma^2} \right)  + \exp \left(- \frac{(x + c_\text{eqk}t)^2} {2\sigma^2} \right) \right]
\end{equation}
for an amplitude scaling $A=0.2$ m/s, a unitless regularization parameter $\varepsilon = 0.5$, pulse width $\sigma = 3$ km, a decay parameter $L = 20$ km to account for geometric attenuation of the surface waves. The general formulation described by \autoref{eq:ksi surface waves} and \autoref{eq:dksidt surface waves} approximates these surface waves as the sum of two Gaussians propagating in opposite directions from a source centered at $x=0$ m (the considered rupture location) and at speeds that are a given $n$ integer multiples of the corresponding approximate tsunami speed ($\sqrt{g H_0}$). Such a problem configuration illustrates the two scales of tsunamigenesis considered in what follows: an earthquake surface wave propagating quickly with respect to the tsunami waves, but not quickly enough to neglect the earthquake dynamics.  \autoref{fig:surface wave source} presents an example of the time evolution of this dynamic source for $c_\text{eqk}=10 \sqrt{gH_0}$ (i.e., $n=10$ in \autoref{eq:ceqk}) given by \autoref{eq:ksi surface waves} and \autoref{eq:dksidt surface waves}, where numerical integration of the former is facilitated by a straightforward trapezoidal rule. For a depth of $H_0 = 4$ km (considered below), such an example yields $c_\text{eqk} \approx 2$ km/s, approaching the Rayleigh wave speed.

\begin{remark}\label{rem:source}Seismogenic tsunamis are classically modeled by the final (static) vertical displacement (i.e., neglecting its time-dependence as well as its velocity), applied as an initial condition at $t=0$~\cite{dynamicsource1}. Although this approximation may be valid in the far-field (i.e., for long distance propagation), such sourcing may not be well-suited in the near-field (the primary motivation of this work and of the problem introduced here), possibly leading to considerable errors with implications on hazard assessment~\cite{dynamicsource1, amlanibhat2022}. For comparison, in all the results that follow for 1D and 2D, statically-sourced solutions are presented alongside the dynamically-sourced results at different earthquake speeds (both produced by the FC-based solver). This includes static displacements applied conventionally at $t=0$ as well as instantaneous displacements at times $t=t_\text{eqk}$ (the latter following~\cite{amlanibhat2022}), where $t_\text{eqk}$ corresponding to the time at which the surface waves have achieved their final displacement in the domain for each earthquake speed $c_\text{eqk}$.\end{remark}

\begin{figure*}[!t]
    \centering
    \begin{subfigure}[b]{0.4\textwidth}
        \centering
        \includegraphics[width=.95\textwidth]{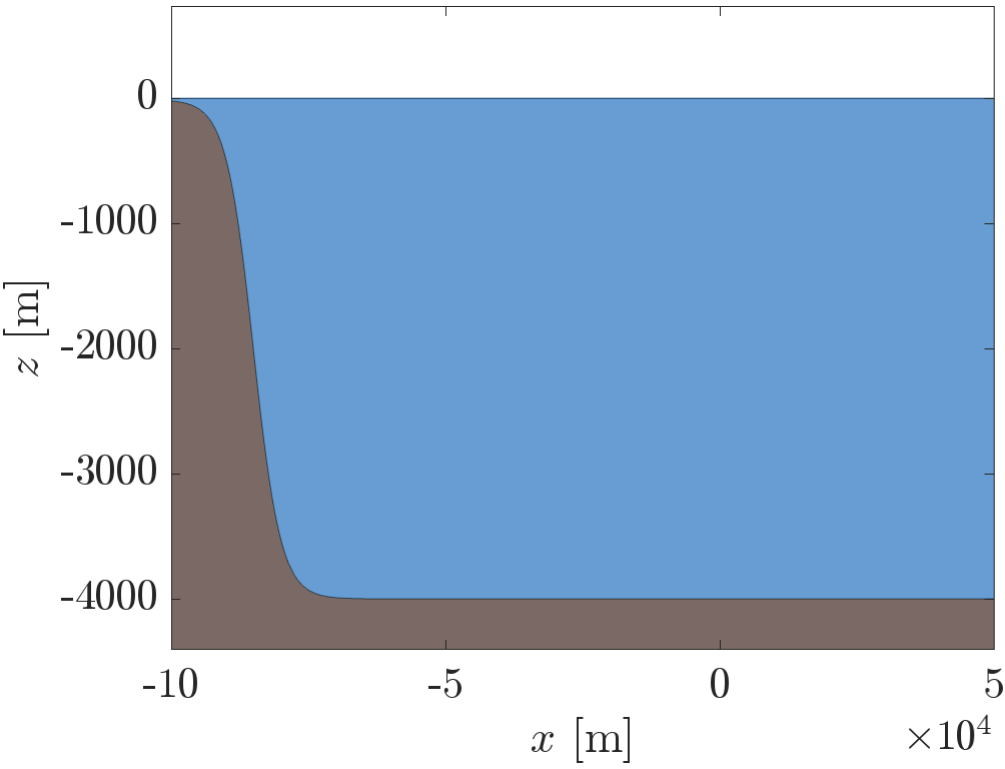}
        \caption[eta]%
        {{\footnotesize $t=0.0$ s;}}
        \label{subfig:1D surface wave init cond}
    \end{subfigure}
    \begin{subfigure}[b]{0.4\textwidth}
        \centering 
        \includegraphics[width=.95\textwidth]{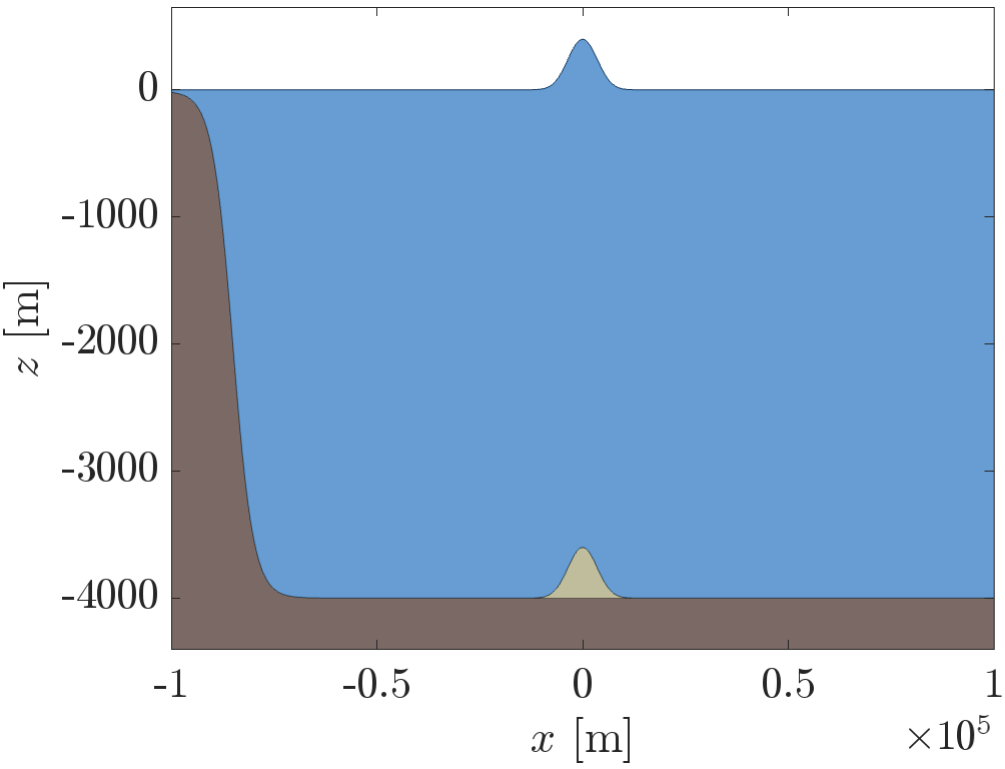}
        \caption[]%
        {{\footnotesize $t=1.7$ s;}}
        \label{subfig:1D surface wave propa1}
    \end{subfigure}
    
    \smallskip
    \begin{subfigure}[b]{0.4\textwidth}  
        \centering 
        \includegraphics[width=.95\textwidth]{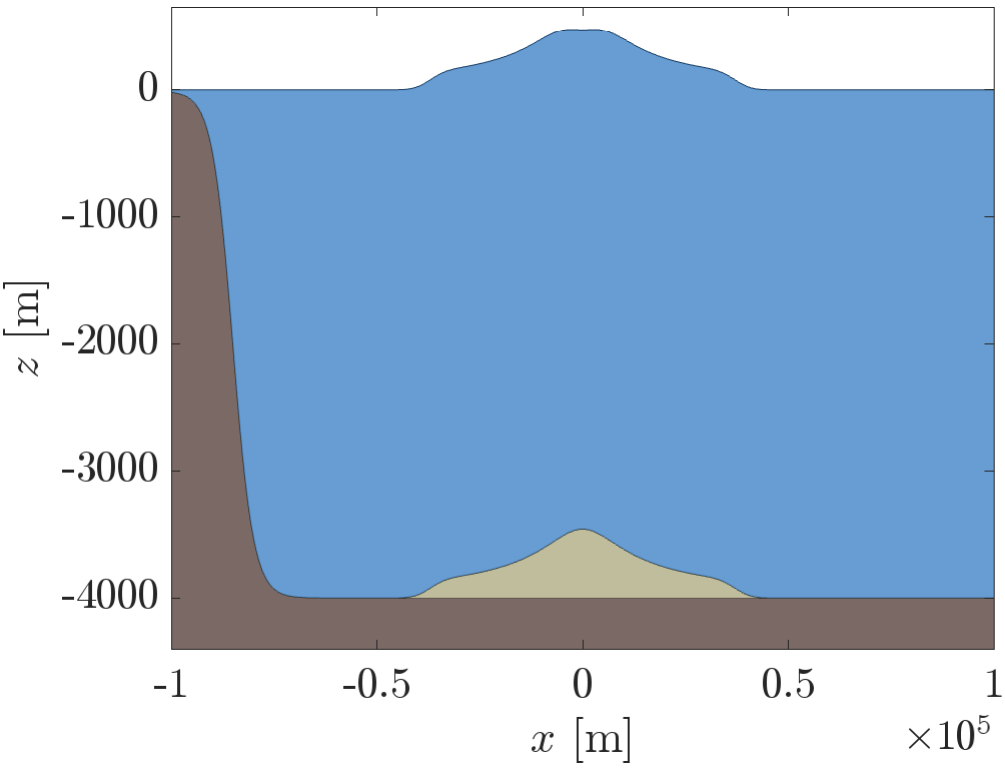}
        \caption[]%
        {{\footnotesize $t=18.9$ s;}}
        \label{subfig:1D surface wave propa2}
    \end{subfigure}
    \begin{subfigure}[b]{0.4\textwidth}   
        \centering 
        \includegraphics[width=.95\textwidth]{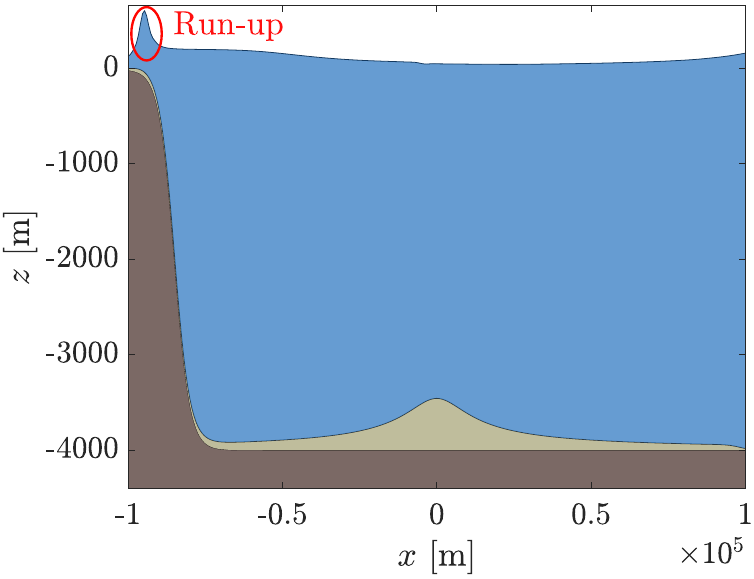}
        \caption[]%
        {{\footnotesize $t=578.4$ s.}}
        \label{subfig:1D surface wave runup}
    \end{subfigure}
    \caption{\emph{1D earthquake-induced tsunamigenesis}. Snapshots at various times of tsunami generation and propagation (blue), including the dynamic earthquake source of speed $c_\text{eqk}=10 \sqrt{gH_0}$ (yellow). The last snapshot indicates the time/location at which run-up is evaluated. The free surface displacement and the ground perturbation are amplified for visualization ($\times100$).}
    \label{fig:surface wave frames}
\end{figure*}

\begin{figure}[!t]
	\centering
	\includegraphics[width=0.4\textwidth]{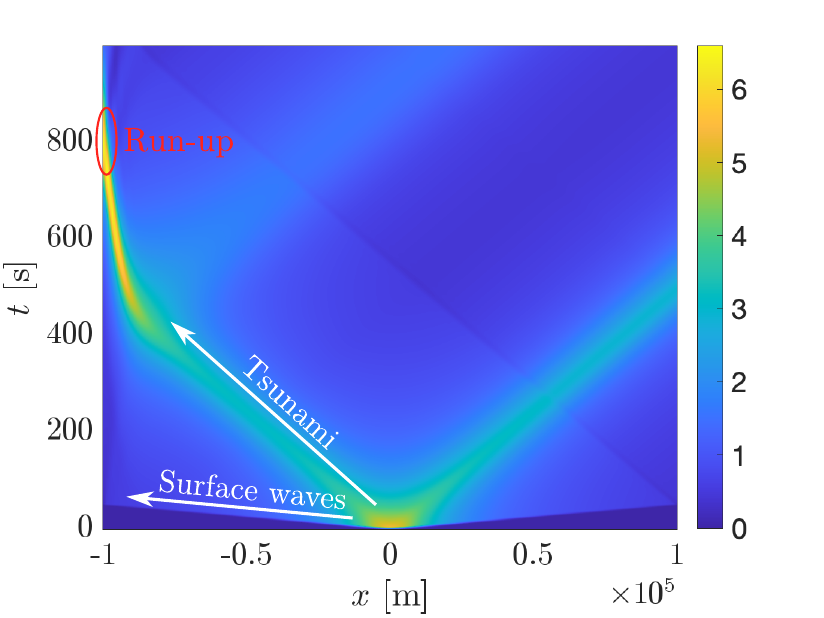}
	\caption{\emph{1D earthquake-induced tsunamigenesis}. The complete space-time evolution of the absolute value of the free surface displacement $\eta$ corresponding to the dynamic earthquake source  of speed $c_\text{eqk}=10 \sqrt{gH_0}$.}
	\label{fig:surface wave propa}
\end{figure}

Over a domain $x\in[-100,100]$ km, the initial (unperturbed) bathymetry is that of a hyperbolic tangent beach, given by $h_0(x)=H_0-s(x)$ for maximum depth $H_0=4$ km (close to typical mean ocean depth), where the initial seafloor topography $s(x)$ is given by
\begin{equation}\label{eq:bathymetry surface waves}
        \displaystyle s(x) = \frac{S}{2}\left(1-\text{tanh}\left(a\left(x-x_1\right)\right)\right)\left(x_1-x\right)
\end{equation}
for $x_1 = -85$ km, $a= 2\times 10^{-4}$ m$^{-1}$ and $S = 3.99$ km. The initial condition corresponds to that of still water (i.e., $\eta(x,t=0)=0$, $u(x,t=0)=0$), and the boundary conditions are prescribed as non-reflecting (see \autoref{eq:1Dradiationbc}).
\autoref{fig:surface wave frames} presents snapshots at various times. produced by the FC-based solver, of the tsunami propagation and its ground motion source corresponding to an example earthquake speed of $c_\text{eqk}=10 \sqrt{gH_0}$, along with a visualization of the unperturbed bathymetry in \autoref{subfig:1D surface wave init cond}. Such results are produced by a simulation that is advanced to a final time $t_\text{max}=1000$ s at a spatial discretization corresponding to $N_x = 300$ points.
\autoref{fig:surface wave propa} presents the complete space-time evolution of the corresponding tsunami displacement $\eta(x,t)$ produced by the FC solver for such a source.  One can observe two types of waves: the earthquake surface waves that propagate at $c_\text{eqk}=10 \sqrt{gH_0}$ and that reach the domain edge at around a time $t=50$ s, and the tsunami waves with their longer wavelengths and slower speeds ($\sqrt{gH_0}$). The left-going tsunami wave slows down as $h_0(x) = H_0-s(x)$ decreases with the approaching beach (similarly to the benchmark of \autoref{sec:benchmark2}).

\begin{figure*}[!t]
    \centering
    \captionsetup[subfigure]{oneside,margin={0.6cm,0cm}}
    \begin{subfigure}[t]{0.8\textwidth}
        \centering 
        	\includegraphics[width=\textwidth]{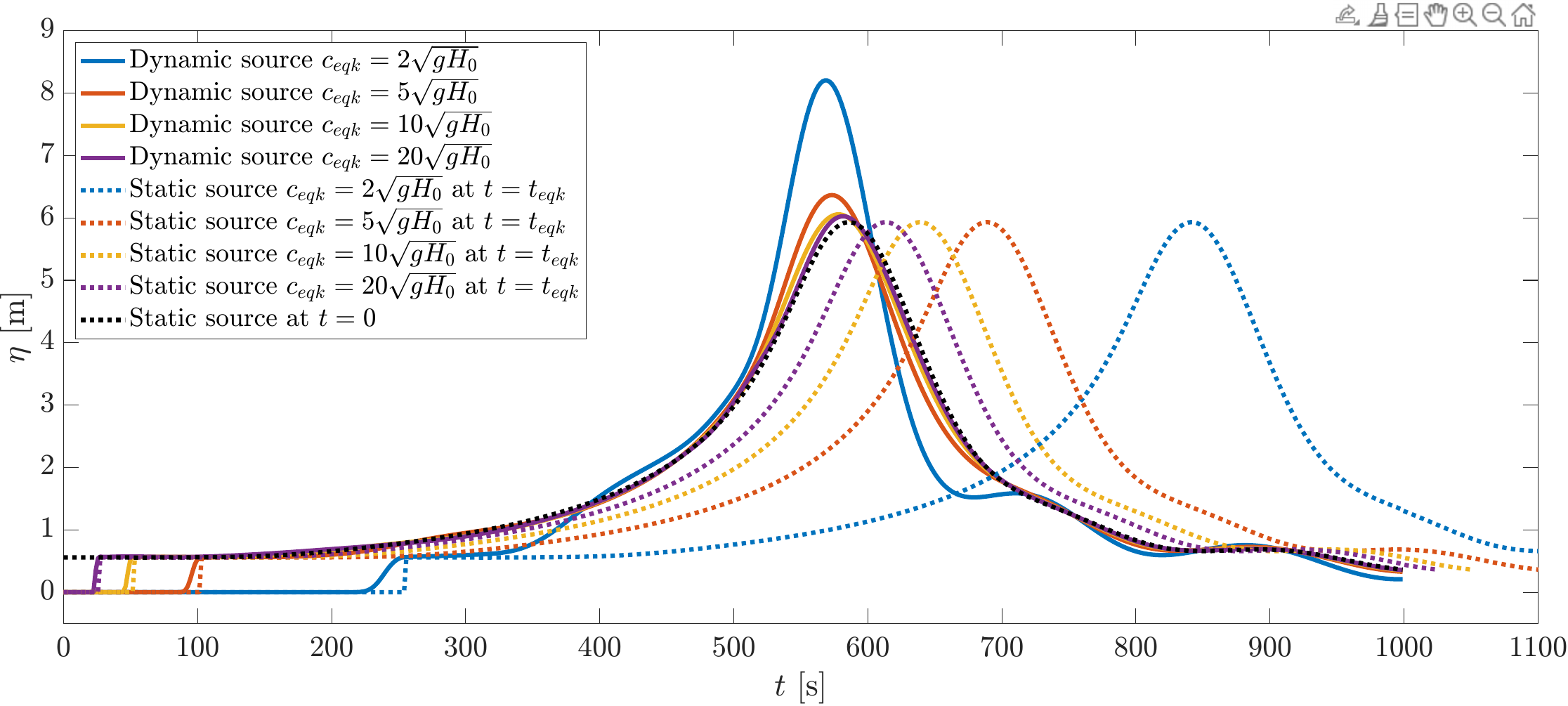}
        \caption[]%
        {{\footnotesize Water height at $x=-95$ km ;}}
        \label{1D Linear plot}
    \end{subfigure}
    % \bigskip
    \vskip\baselineskip
    \begin{subfigure}[t]{0.8\textwidth}
        \centering 
        \includegraphics[width=\textwidth]{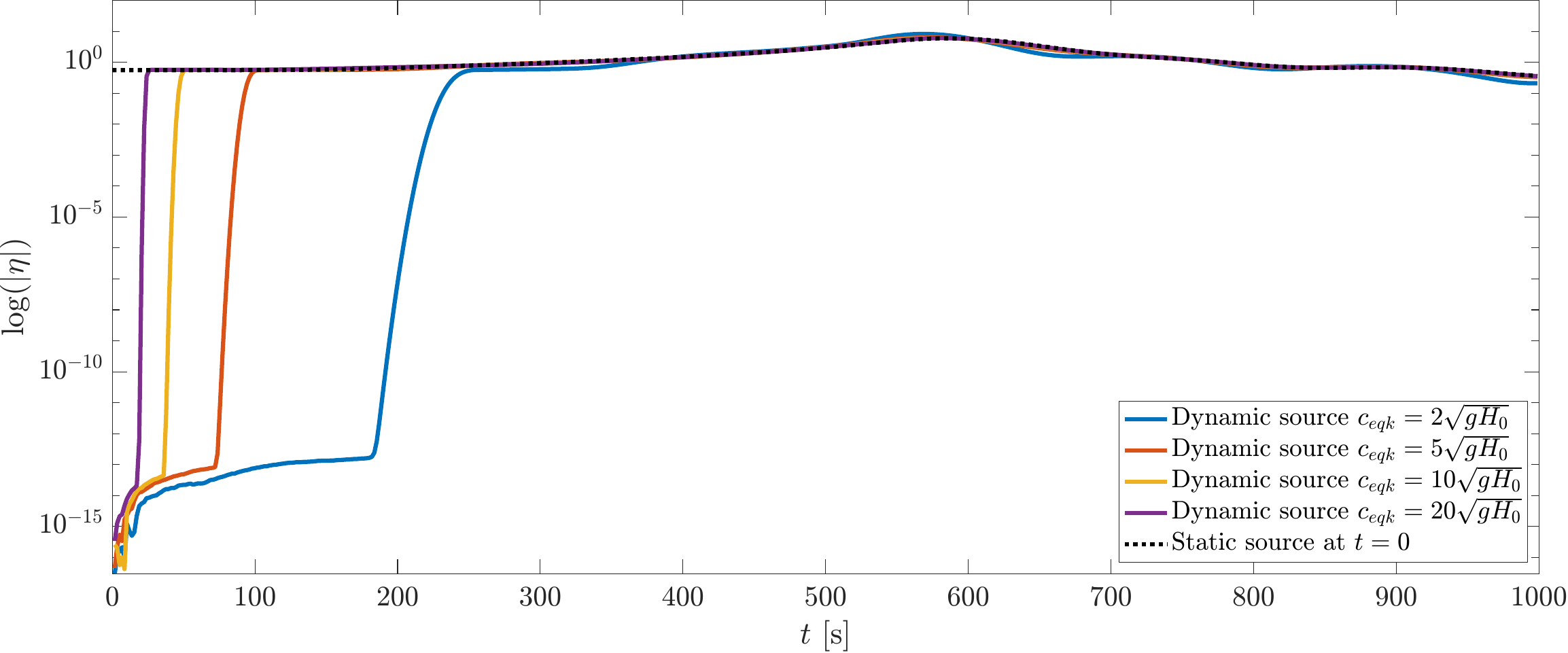}
        \caption[]%
        {{\footnotesize Log-scale water height at $x=-95$ km.}}
        \label{1D Log plot}
    \end{subfigure}
    \caption{\emph{1D earthquake-induced tsunamigenesis}. Time evolution of the free surface displacement $\eta(x,t)$ at the coastline for both time-dependent and classical static (instantaneous) seafloor perturbations at different sourcing earthquake propagation speeds. The final vertical displacement of the static sources are applied at $t=0$ or at $t=t_\text{eqk}$ (the final time of the corresponding earthquake motions, see~Remark~\ref{rem:source}).}
    \label{fig:1D surface wave curves}
\end{figure*}

 In order to emphasize the importance of considering earthquake dynamics in tsunami modeling (the main motivation of the solver introduced in this work), \autoref{1D Linear plot} presents the time evolution of water elevation $\eta(x_\text{min},t)$ for varying earthquake speeds $n = c_\text{eqk}/\sqrt{gH_0} \in \{2, 5, 10, 20 \}$ near the left boundary at a location corresponding to $x=-95$ km (indicated as "run-up" in the bottom right panel of \autoref{fig:surface wave frames}). It is notable here that no numerical tuning is necessary to achieve such results (all employ the same spatial discretization and CFL constant) and that the longest computational time (corresponding to $n=20$) is around two seconds. The curves illustrate the influence of the earthquake source speed $c_\text{eqk}$ on the displacement found at the highest point of the beach (left boundary), demonstrating that, when considering the full dynamics of the ground motion ("dynamic source"), larger water heights are observed for slower earthquakes. Additionally presented is the water elevation corresponding to the classical "static source" model associated with each earthquake speed. The static source is conventionally invoked as an instantaneous displacement representing the final or total displacement of the earthquake.
% imposed at the end of the earthquake duration The time 
 One can observe that the higher the earthquake propagation speed, the better the approximation of the static source (i.e., instantaneous displacement). However, as the earthquake source speed becomes more comparable to the characteristic propagation speed of tsunami waves ($\sqrt{gH_0}$, where $H_0$ is the maximum depth), the static source response leads to considerable underestimation of the water elevation at the coastline. This is quantified for each case in \autoref{table:1D tsunamigenesis amplitudes}, where one can observe that amplitude differences can be up to around $30\%$ for the slowest earthquake speed.

 \begin{table}[!t]
    \centering
    \small
    \begin{tabular}{|c |c |c|}
    \hline
        \multicolumn{3}{|c|}{\bfseries Tsunami height}\\
         \hline
          $c_\text{eqk}$ & Dynamic source& Static source \\ [0.5ex]
         \hline
         $2 \sqrt{gH_0}$ & 8.20 m & 5.93 m ($-27.7\%$)  \\ 
         \hline
         $5 \sqrt{gH_0}$ & 6.37 m & 5.93 m ($-06.8\%$) \\ 
         \hline
         $10 \sqrt{gH_0}$ & 6.06 m & 5.93 m ($-02.1\%$) \\ 
         \hline
         $20 \sqrt{gH_0}$ & 6.03 m & 5.93 m ($-01.6\%$) \\ 
         \hline
    \end{tabular}
    \caption{\emph{1D earthquake-induced tsunamigenesis}. Simulated tsunami amplitudes at coastline $x=-95$ km produced by the FC-based solver employing both the time-dependent behavior of the ground source as well as the classical static (instantaneous) source. The tsunami height underestimation due to the latter is indicated in parenthesis (\%).}
    \label{table:1D tsunamigenesis amplitudes}
\end{table}

 The static sources in \autoref{1D Linear plot} additionally overestimate the arrival times of their respective peaks in comparison to their corresponding dynamic source models (such differences become significantly more pronounced for slower earthquakes). This is in agreement with what has been similarly observed elsewhere~\cite{amlanibhat2022} and is expected since such static models apply their instantaneous displacements at the end-time of the earthquake. Such a time must be given as an input, which may be difficult to determine; no such input is required for a dynamic source. It is notable here (and also presented in \autoref{1D Linear plot}) that even in an overly-conservative scenario, where the static source is applied as an initial condition of the water elevation at $t=0$ (at the beginning of the earthquake instead of at the end), such a delay (purely due to neglecting the dynamics) is still observable when compared to the dynamic source results for slower earthquakes. To the best of the author's knowledge, the delay and its implications as a function of earthquake speed have not been explicitly demonstrated before in a parameter study for such configurations  (strike-slip ruptures have been studied only recently~\cite{Elbanna2021,amlanibhat2022}). Indeed, a quantification of arrival times in \autoref{table:1D tsunamigenesis arrival times} shows that the discrepancy in arrival times is at least 25 seconds for the fastest earthquake speed considered and can be more than four minutes for the slowest when the static displacement is applied at $t=t_\text{eqk}$. Even for a conservative $t=0$, the overestimation can still be observed.
 
%  these results demonstrate, for the first time, the possible importance of dynamics for subduction-like configurations (as opposed to the in-plane strike-slip ruptures considered in \cite{Elbanna2021,amlanibhat2022}).
%  here can be a delay, of purely dynamic nature, that would overestimate the time arrivals of tsunami waves.
 
% Here, a very novel and interesting remark is that the tsunami arrival is expected with a delay for the static source, even if the static source was set at the beginning of earthquake duration. This observation means that even if the static models were set at the worst settings (water elevation at the beginning of earthquake signal), there can be a delay, of purely dynamic nature, that would overestimate the time arrivals of tsunami waves.

% Classically, in tsunami modeling the source is a water elevation that corresponds to the final seafloor displacement, imposed at the end of earthquake duration so there is an expected delay with respect to a dynamic model, which is very inconvenient for hazard assessment because it overestimates the time available for the population in danger.

% The amplitude improvements corresponding to each surface wave speed from using the proposed dynamic source model instead of the classical static source is summarized \autoref{table:1D tsunamigenesis amplitudes}.

% The wave arrival improvements corresponding to each surface wave speed from using the proposed dynamic source model instead of the classical static source is summarized \autoref{table:1D tsunamigenesis arrival times}.

\begin{table}[!t]
    \centering
    \small
    \begin{tabular}{|c |c |c| c |} 
      \hline
    \multicolumn{4}{|c|}{\bfseries Tsunami arrival time}\\
         \hline
       
          $c_\text{eqk}$ & Dynamic source & Static source at $t=0$ & Static source at $t=t_\text{eqk}$   \\ [0.5ex]
         \hline
         $2 \sqrt{gH_0}$ & 568.1 s & 585.3 s ($+17.2$ s) & 837.7 s ($+269.6$ s) \\ 
         \hline
         $5 \sqrt{gH_0}$ & 573.3 s & 585.3 s ($+12.0$ s) & 686.3 s ($+113.0$ s)  \\ 
         \hline
         $10 \sqrt{gH_0}$ & 578.4 s & 585.3 s ($+06.9$ s) & 635.8 s ($+057.4$ s) \\ 
         \hline
         $20 \sqrt{gH_0}$ & 581.9 s & 585.3 s ($+03.4$ s) & 610.5 s ($+028.6$ s)  \\
         \hline
    \end{tabular}
    \caption{\emph{1D earthquake-induced tsunamigenesis}. Simulated arrival times at coastline $x=-95$ km produced by the FC-based solver employing both the time-dependent behavior of the ground source as well as the classical static (instantaneous) source applied at $t=0$ and at the end of the ground motion at $t=t_\text{eqk}$ (see Remark~\ref{rem:source}). The tsunami arrival time differences for static sourcing with respect to dynamic sourcing are indicated in parenthesis.}
    \label{table:1D tsunamigenesis arrival times}
\end{table}

% \autoref{Log plot} additionally presents the corresponding responses ina semilog plot 

% One can observe in \autoref{fig:surface wave curves}

\autoref{1D Log plot} presents the same water surface heights in a log-axis at the same left boundary, revealing the arrivals of the rupture waves and the initial uplift of the water. These possibly-measurable low amplitude signals in dynamic water elevations may be valuable for tsunami forecasting.
% since they are within measurement capabilities of underwater pressure sensors~\cite{RED}
Indeed, ongoing work entails full parametric studies in order to establish correlations between such signatures and the ultimate timing/amplitude of the resulting primary tsunami waves.  

% One can observe that they scale proportionally to the tsunami wave heights, so capturing these signals (for example, using the sea floor pressure sensors or buoys near the coasts) may allow to infer a quantification of the tsunami risk.

% Additionally, the static source overestimates the arrival time of the peak in comparison to the dynamic source models of the slower earthquakes. It is notable here that this static source was applied as an initial condition of the water elevation at $t=0$ instead of at the end of the earthquake duration as is classically taken; such a delay in wave arrival is even more pronounced 

% the methodology taking into account the sismogram of the signal has been compared with a \textit{static} approach where the initial condition is represented by the final (total) displacement caused by the earthquake. This comparison has been done for multiple earthquake speeds : $c_\text{eqk}/\sqrt{gh} \in \{2, 5, 10, 20 \}$ and is represented \autoref{fig:surface wave curves}.

% \begin{figure}[!t]
% 	\centering
% 	\includegraphics[width=0.7\textwidth]{Surface waves delay.pdf}
% 	\caption{Evolution of the free surface $\eta$ at the left domain edge (run-up) for dynamic and static seafloor movements defined \autoref{eq:dksidt surface waves} and \autoref{eq:ksi surface waves} considering different earthquake propagation speeds. The uplifts of the static seafloors are set at the end of earthquake duration (delay that is classically considered).}
% 	\label{fig:surface wave curves delay}
% \end{figure}

\subsection{2D formulation}
\begin{figure}[!t]
	\centering
	\includegraphics[width=0.5\textwidth]{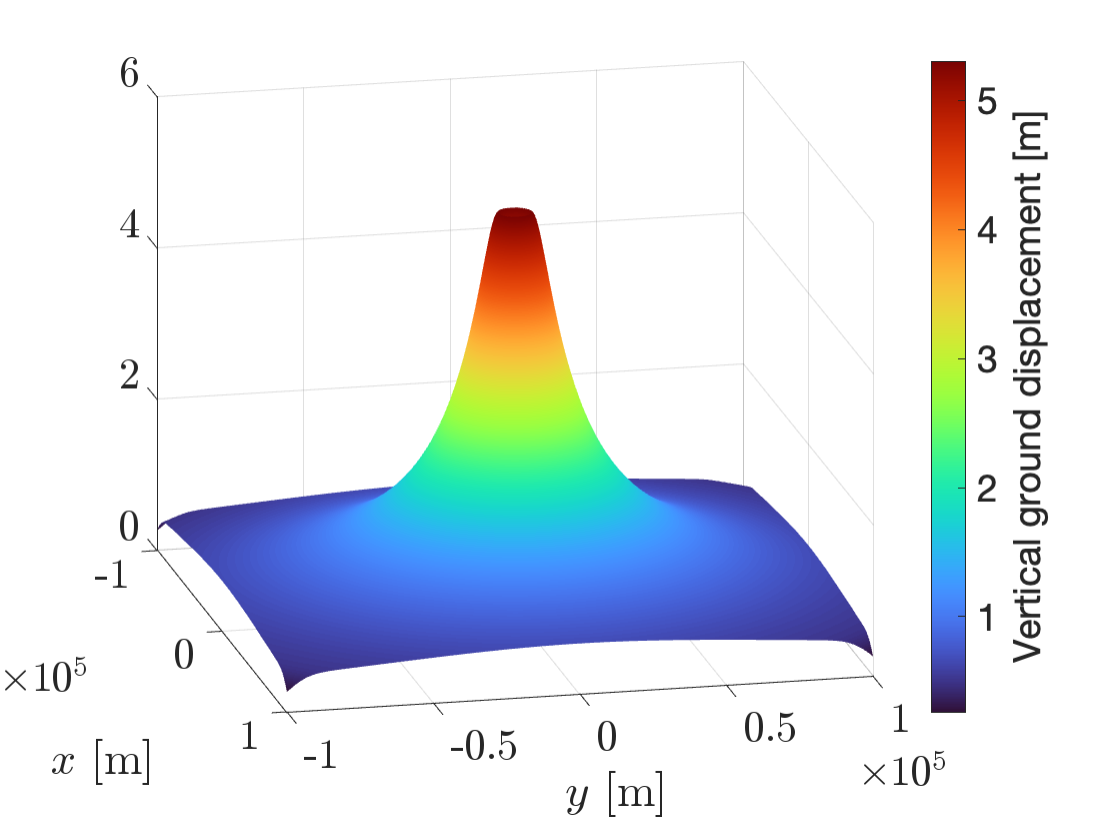}
	\caption{\emph{2D earthquake-induced tsunamigenesis}. The final 2D ground displacement for the source $\xi(x,y,t)$ corresponding to \autoref{eq:ksi 2D surface waves}.}       
	\label{fig:2D surface wave ksi final}
\end{figure}

The problem configuration in 1D can be extended to 2D with a similar time-dependent seafloor displacement, inspired by realistic scenarios, given by
\begin{equation}\label{eq:ksi 2D surface waves}
    \xi(x,y,t)=\int_{0}^{t} \xi_t(x,y,\tau) \text{d}\tau,
\end{equation}
where $\xi_t$ is the time-derivative of the vertical displacement (i.e., vertical seafloor velocity) given by
\begin{equation}\label{eq:dksidt 2D surface waves}
    \xi_t(x,y,t)=\frac{nA}{\sqrt{(r/L)^2 + \varepsilon^2}} \exp \left(-\frac{(r - c_\text{eqk}t)^2} {2\sigma^2} \right)
\end{equation}
for radial distance $r=\sqrt{x^2+y^2}$, an amplitude scaling $A=1.67 \times 10^{-1}$ m/s, a unitless regularization parameter $\varepsilon = 0.5$, pulse width $\sigma = 4.24$ km, and a decay parameter $L = 20$ km.
The final 2D displacement is illustrated in \autoref{fig:2D surface wave ksi final}.

\begin{figure*}[!t]
    \centering
    \begin{subfigure}[b]{0.4\textwidth}
        \centering
        \includegraphics[width=\textwidth]{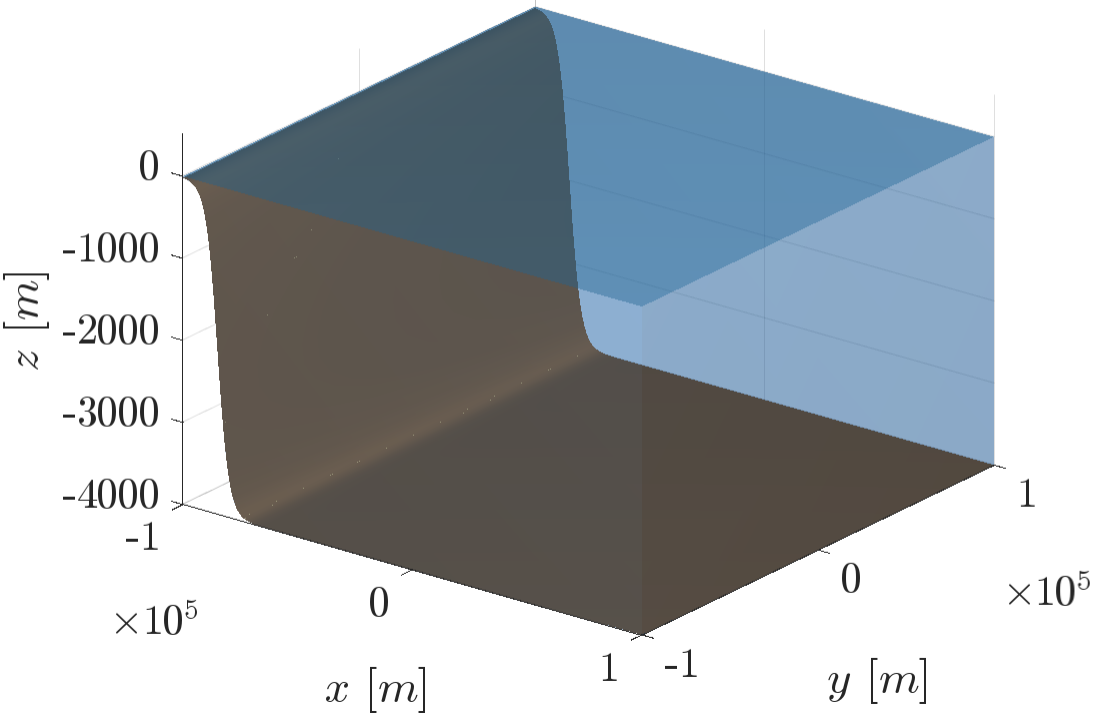}
        \caption[eta]%
        {{\footnotesize $t =0$ s;}}
        \label{fig:Surface waves 0}
    \end{subfigure}\quad
    \begin{subfigure}[b]{0.4\textwidth}
        \centering
        \includegraphics[width=\textwidth]{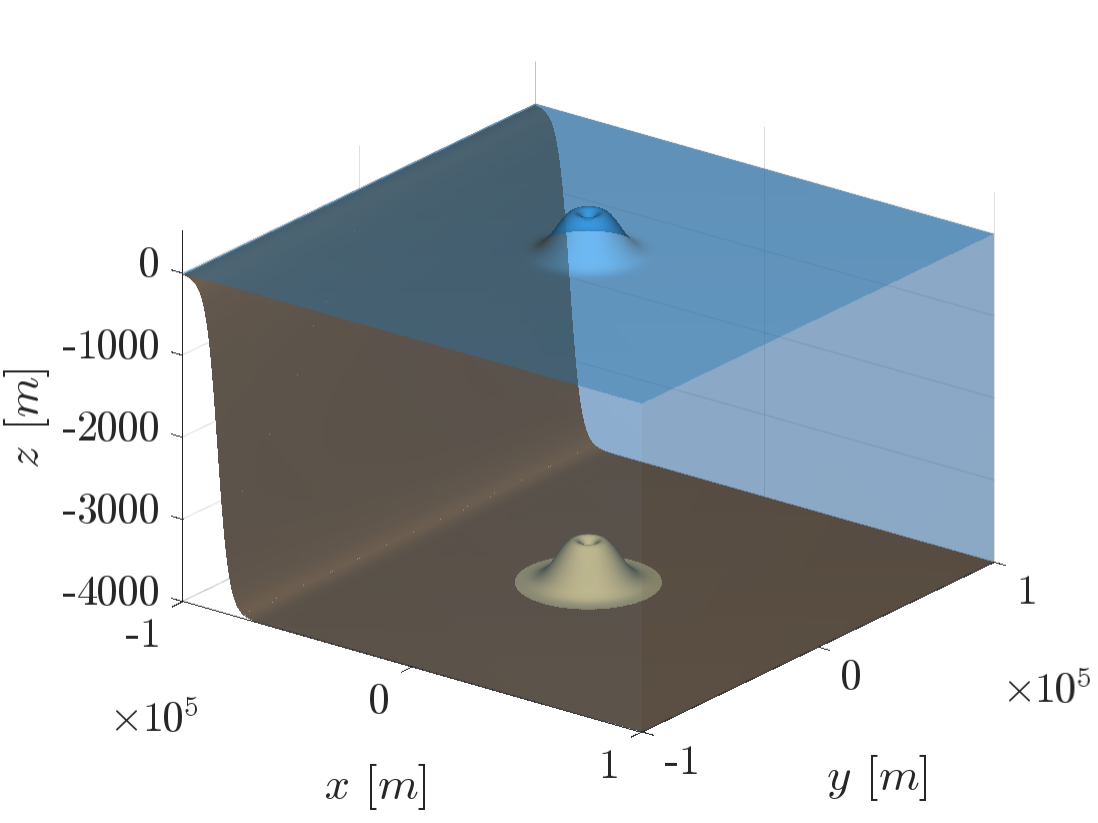}
        \caption[eta]%
        {{\footnotesize $t =7$ s;}}
        \label{fig:Surface waves 1}
    \end{subfigure}
 
 \medskip
 
    \begin{subfigure}[b]{0.4\textwidth}
        \centering 
        \includegraphics[width=\textwidth]{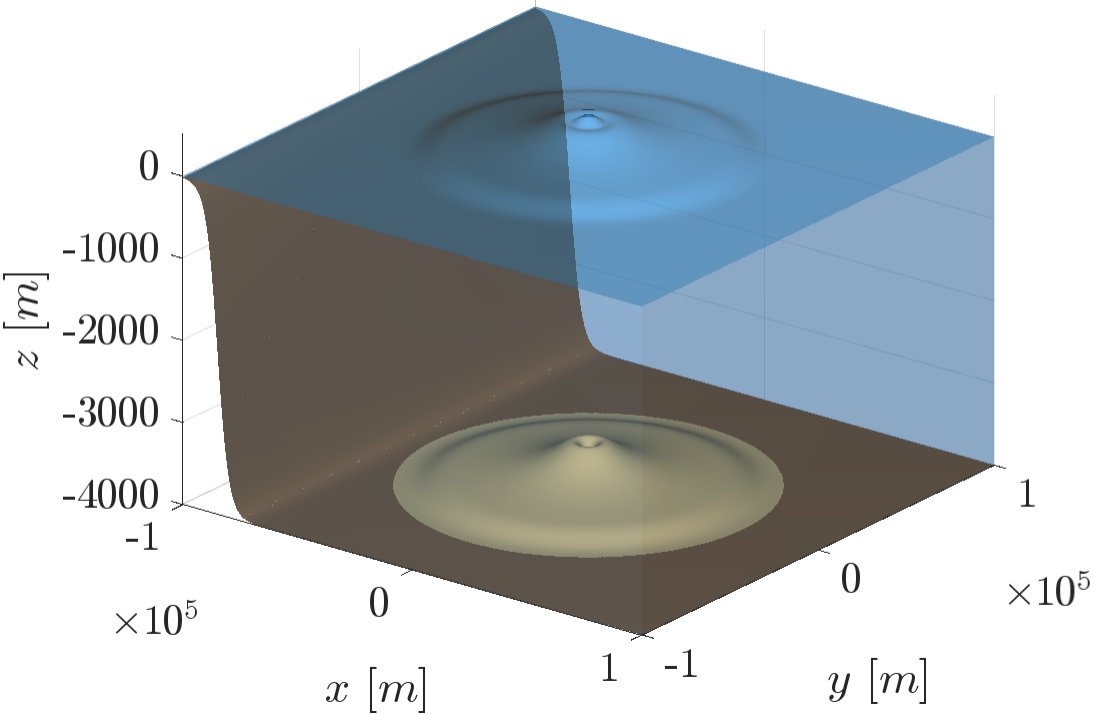}
        \caption[]%
        {{\footnotesize $t =29$ s;}}
        \label{fig:Surface waves 2}
    \end{subfigure}\quad
    \begin{subfigure}[b]{0.4\textwidth}
        \centering
        \includegraphics[width=\textwidth]{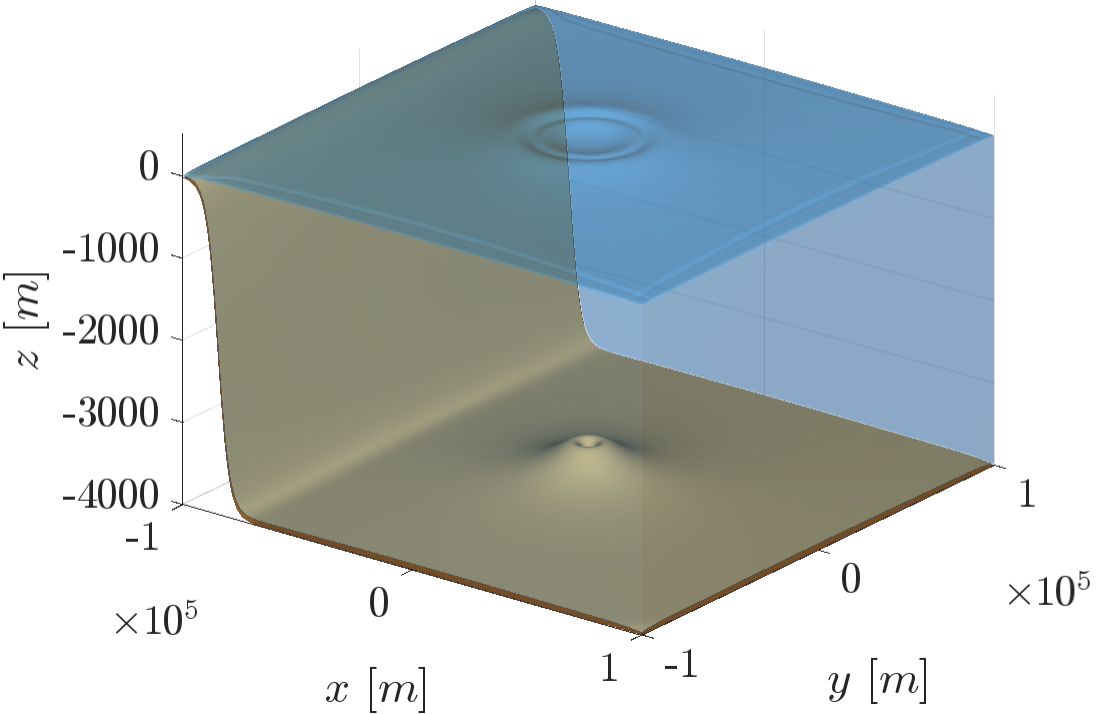}
        \caption[eta]%
        {{\footnotesize $t =99$ s;}}
        \label{fig:Surface waves 3}
    \end{subfigure}
    
    \medskip
    
    \begin{subfigure}[b]{0.4\textwidth}  
        \centering 
        \includegraphics[width=\textwidth]{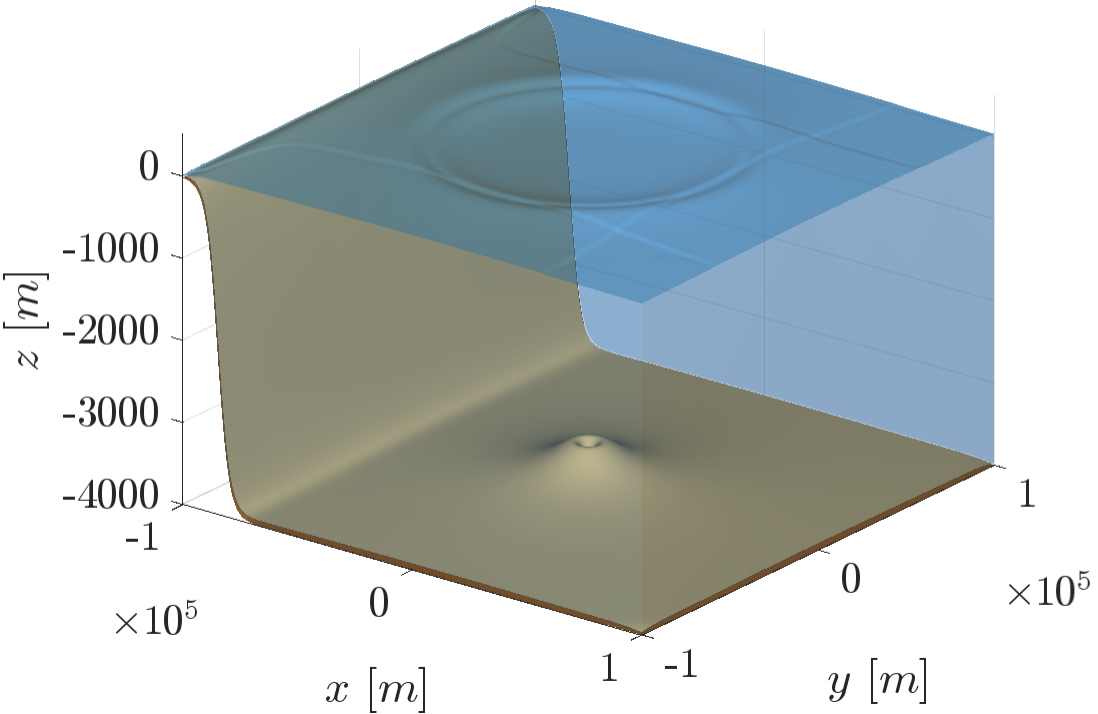}
        \caption[]%
        {{\footnotesize $t =284$ s;}}
        \label{fig:Surface waves 4}
    \end{subfigure}\quad
    \begin{subfigure}[b]{0.4\textwidth}   
        \centering 
        \includegraphics[width=\textwidth]{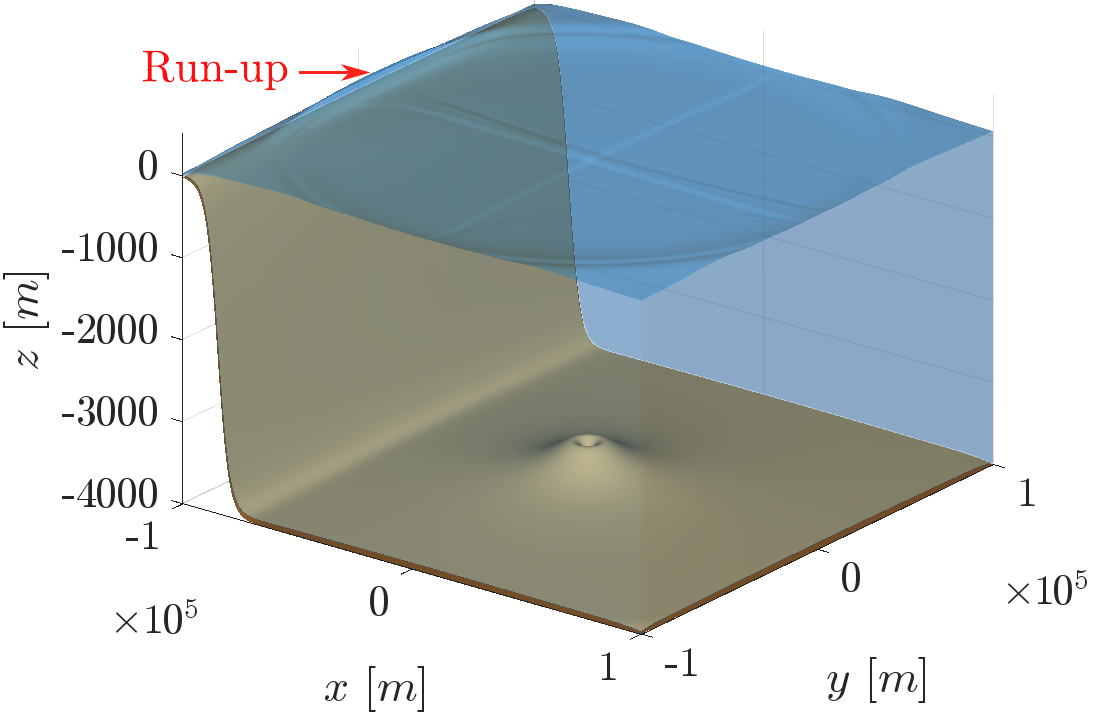}
        \caption[]%
        {{\footnotesize $t =538$ s (maximum height at $x=-95$ km).}}
        \label{fig:Surface waves 5}
    \end{subfigure}
    \caption{\emph{2D earthquake-induced tsunamigenesis}. Snapshots at various times of tsunami generation and propagation, including the dynamic earthquake source of speed $c_\text{eqk}=10 \sqrt{gH_0}$. The free surface displacement (blue) and the ground perturbation (yellow) are amplified for visualization ($\times100$).}
    \label{fig:Surface waves 2D}
\end{figure*}

Over a similarly-sized domain $x,y\in[-100,100]$ km, the initial (unperturbed) bathymetry is that of a hyperbolic tangent beach similar to~\autoref{sec:1D application}, similarly given by $h_0(x,y)=H_0-s(x,y)$ for maximum depth $H_0 = 4$ km, where the initial seafloor topography $s(x,y)$ given by 
\begin{equation}\label{eq:2D bathymetry surface waves}
        \displaystyle s(x,y) = \frac{S}{2}\left(1-\text{tanh}\left(a\left(x-x_1\right)\right)\right)\left(x_1-x\right)
\end{equation}
for $x_1 = -85$ km, $a= 2\times 10^{-4}$ m$^{-1}$ and $S = 3.99$ km. The initial condition corresponds to that of still water (i.e. $\eta(x,y,0)=0$, $u(x,y,0)=0$ and $v(x,y,0)=0$), and the boundary conditions are prescribed as non-reflecting (see \autoref{eq:2Dradiationbc}). All simulations that follow employ $N_x = N_y = 300$ discretization points, up to a final time $t_\text{max}=1\,000$ s. \autoref{fig:Surface waves 2D} presents snapshots at various times, produced by the FC-based solver, of the tsunami propagation and its ground motion source corresponding to an example earthquake speed of $c_\text{eqk}=10 \sqrt{gH_0}$. As illustrated in the figures, earthquake surface waves propagate on the sea floor, leading to a corresponding uplift of water height of the same value due to the non-compressibility hypothesis (\autoref{fig:Surface waves 1}). The ground motion source propagates at $c_\text{eqk}=10 \sqrt{gH_0}$ (\autoref{fig:Surface waves 2}) until its waves reach the domain boundaries (\autoref{fig:Surface waves 3}). Meanwhile the corresponding tsunami wave propagates at a longer time-scale (\autoref{fig:Surface waves 4}) until the maximal wave height is achieved near the coastline (\autoref{fig:Surface waves 5}).

\begin{figure*}[!t]
    \centering
    \captionsetup[subfigure]{oneside,margin={0.6cm,0cm}}
    \begin{subfigure}[t]{0.8\textwidth}
        \centering 
        	\includegraphics[width=\textwidth]{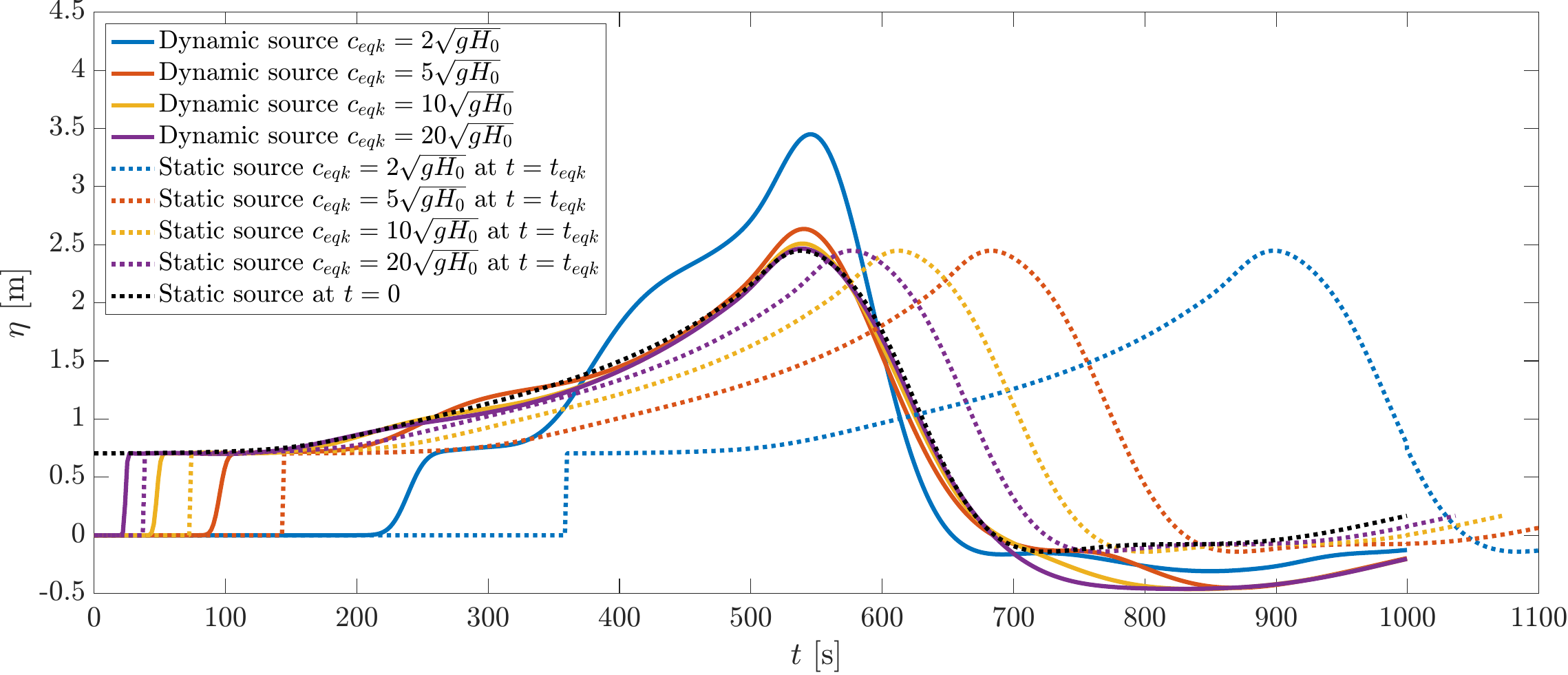}
        \caption[]%
        {{\footnotesize Water height at $(x=-95, y=0)$ km;}}
        \label{2D Linear plot}
    \end{subfigure}
    % \bigskip
    \vskip\baselineskip
    \begin{subfigure}[t]{0.8\textwidth}
        \centering 
        \includegraphics[width=\textwidth]{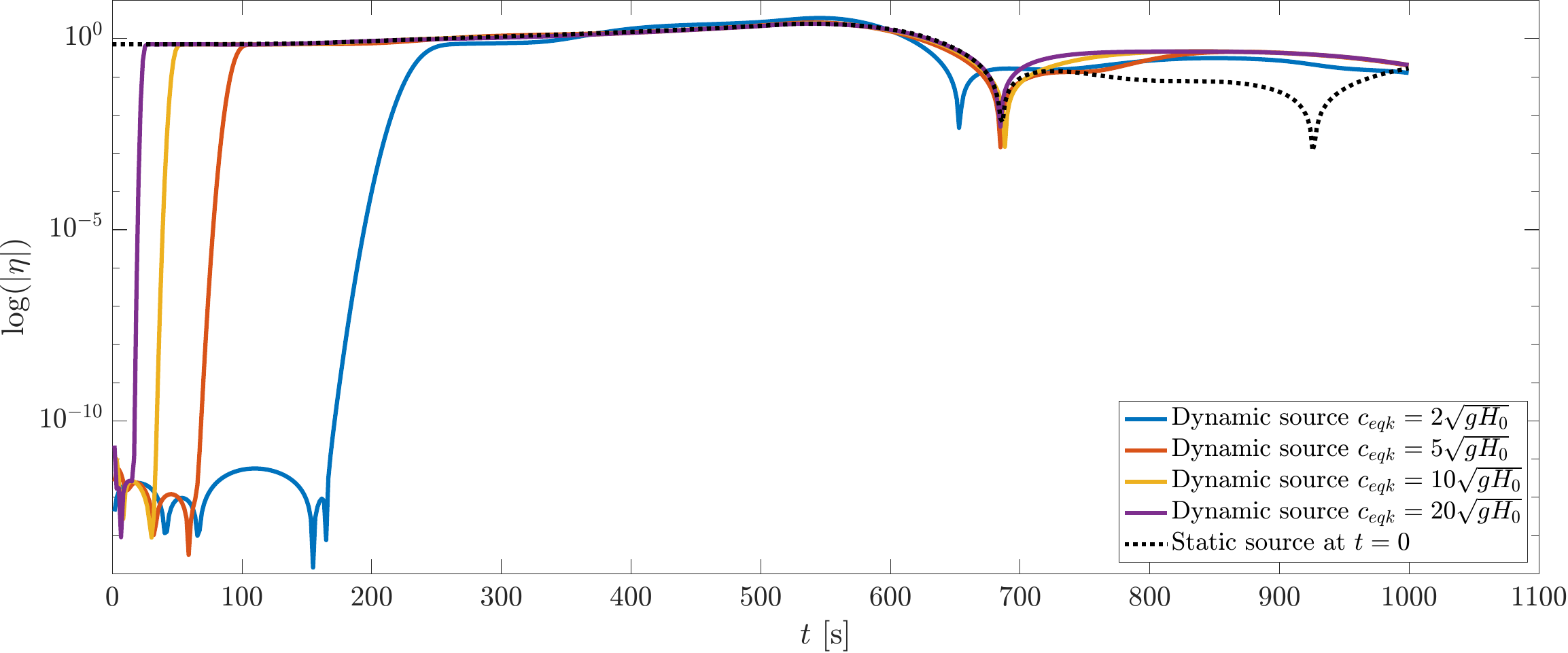}
        \caption[]%
        {{\footnotesize Log-scale water height at $(x=-95, y=0)$ km.}}
        \label{2D Log plot}
    \end{subfigure}
    \caption{\emph{2D earthquake-induced tsunamigenesis}. Time evolution of the free surface $\eta(x,y,t)$ at the left domain edge (run-up) $(x,y) = (-95,0)$ km for both time-dependent and classical static (instantaneous) seafloor perturbations at different sourcing earthquake propagation speeds.  The final vertical displacement of the static sources are applied at $t=0$ or at $t=t_\text{eqk}$ (the final time of the corresponding earthquake motions, see~Remark~\ref{rem:source}).}
    \label{fig:surface wave curves}
\end{figure*}

Similarly to the previous 1D formulation, the importance of considering earthquake dynamics in 2D is demonstrated in \autoref{2D Linear plot}, where the time evolution of water elevation $\eta(x_\text{min},t)$ is measured for varying earthquake speeds $n = c_\text{eqk}/\sqrt{gH_0} \in \{2, 5, 10, 20 \}$ near the left boundary at a location corresponding to $x=-95$ km and $y=0$ km (indicated as "run-up" in \autoref{fig:Surface waves 5}). The conclusions drawn for the 1D case remain valid for the 2D model: employing time-dependent seafloor motion, larger tsunami heights correspond to slower earthquakes. The "static source" models presented for each earthquake speed highlight the underestimation of the tsunami height and arrival time when employing such classical instantaneous displacements. These observations are quantified for each case in \autoref{table:2D tsunamigenesis amplitudes} and \autoref{table:2D tsunamigenesis arrival times}, which additionally present the differences in water heights and the differences in arrival times, respectively. Additionally, including in the 1D case, the FC-based solver is able to produce all solutions with very low computational cost: all 2D simulations are produced in under 44 seconds in MATLAB (without MEX files). Indeed, the computational times for all the 1D and 2D simulations of this section, presented in \autoref{table: tsunamigenesis computation times},  imply strong potential for real-time or early warning risk assessments (with no added cost for considering time-dependent sources).

\begin{table}[!t]
    \centering
    \small
    \begin{tabular}{|c |c |c|} 
    \hline
    \multicolumn{3}{|c|}{\bfseries Tsunami height}\\
         \hline
          $c_\text{eqk}$ & Dynamic source& Static source \\ [0.5ex]
         \hline
         $2 \sqrt{gH_0}$ & 3.45 m & 2.45 m ($-29.0\%$)  \\ 
         \hline
         $5 \sqrt{gH_0}$ & 2.64 m & 2.45 m ($-07.1\%$) \\ 
         \hline
         $10 \sqrt{gH_0}$ & 2.51 m & 2.45 m ($-02.4\%$) \\ 
         \hline
         $20 \sqrt{gH_0}$ & 2.47 m & 2.45 m ($-00.7\%$) \\ 
         \hline
    \end{tabular}
    \caption{\emph{2D earthquake-induced tsunamigenesis}. Simulated tsunami amplitudes at coastline $(x,y) = (-95, 0)$ km produced by the FC-based solver employing both the time-dependent behavior of the ground source as well as the classical static (instantaneous) source. The tsunami height underestimation due to the latter is indicated in parenthesis (\%).}
    \label{table:2D tsunamigenesis amplitudes}
\end{table}

\begin{table}[!t]
    \centering
    \small
    \begin{tabular}{|c |c |c| c|} 
    \hline
    \multicolumn{4}{|c|}{\bfseries Tsunami arrival time}\\
         \hline
       $c_\text{eqk}$ & Dynamic source & Static source at $t=0$ & Static source at $t=t_\text{eqk}$   \\ [0.5ex]
         \hline
         $2 \sqrt{gH_0}$ & 545.2 s & 538.5 s ($-6.7$ s) & 895.4 s ($+350.2$ s) \\ 
         \hline
         $5 \sqrt{gH_0}$ & 540.2 s & 538.5 s ($-1.7$ s) & 681.3 s ($+141.1$ s)  \\ 
         \hline
         $10 \sqrt{gH_0}$ & 538.5 s & 538.5 s ($+0.0$ s) & 609.9 s ($+071.4$ s) \\ 
         \hline
         $20 \sqrt{gH_0}$ & 538.5 s & 538.5 s ($+0.0$ s) & 574.2 s ($+035.7$ s)  \\
         \hline
    \end{tabular}
    \caption{\emph{2D earthquake-induced tsunamigenesis}. Simulated arrival times at coastline $(x,y) = (-95, 0)$ km produced by the FC-based solver employing both the time-dependent behavior of the ground source as well as the classical static (instantaneous) source applied at $t=0$ and at the end of the ground motion at $t=t_\text{eqk}$ (see Remark~\ref{rem:source}). The tsunami arrival time differences for static sourcing with respect to dynamic sourcing are indicated in parenthesis.}
    \label{table:2D tsunamigenesis arrival times}
\end{table}

\begin{table}[!t]
    \centering
    \small
    \begin{tabular}{|c |c |c|} 
         \hline
          $c_\text{eqk}$ & Computation time (1D)& Computation time (2D) \\ [0.5ex]
         \hline
         static & 0.20 s & 41.39 s \\ 
         \hline
         $2 \sqrt{gH_0}$ & 0.20 s & 40.60 s \\ 
         \hline
         $5 \sqrt{gH_0}$ & 0.19 s & 43.78 s \\ 
         \hline
         $10 \sqrt{gH_0}$ & 0.19 s & 43.67 s \\ 
         \hline
         $20 \sqrt{gH_0}$ & 0.19 s & 39.96 s \\ 
         \hline
    \end{tabular}
    \caption{\emph{2D earthquake-induced tsunamigenesis}. Computation times for the FC-based solver in producing the 1D and 2D solutions of this section (averaged over multiple runs). All simulations are conducted in MATLAB (\autoref{remark:FCparams}).}
    \label{table: tsunamigenesis computation times}
\end{table}

\begin{figure}[!t]
	\centering
	\includegraphics[width=\textwidth]{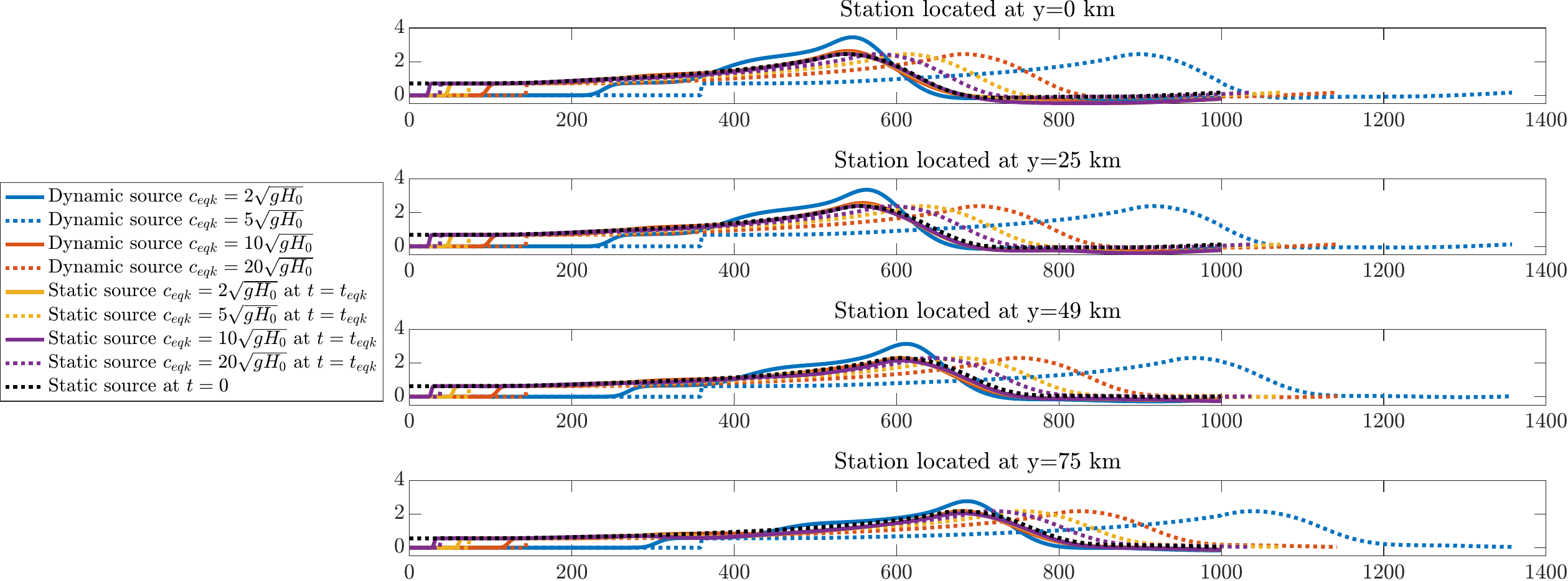}
	\caption{\emph{2D earthquake-induced tsunamigenesis}. Time evolution of the free surface displacement $\eta(x,y,t)$ at various locations (stations) $y$ along the coastline $x=-95$ km for both time-dependent and classical (instantaneous) seafloor pertubations at different sourcing earthquake speeds.}
	\label{fig:surface wave stations}
\end{figure}

For the 2D case presented here, \autoref{fig:surface wave stations} presents the time evolution of the water height at different values (locations) of $y$ on the coastline $x=-95$ km (denoted as simulated "stations" which can correspond in practice to, e.g., DART buoys). The arrival times are expectedly delated in relation to $y=0$ due to the radial nature of the 2D tsunami waves. Even at other locations,  the classical static source model leads to an underestimation of tsunami height and an overestimation of the tsunami arrival time. Such results the efficacy of the FC-based solver in treating problems where complex earthquake ground sources are imposed at the sea floor, particularly through consideration of the dynamic seafloor displacement and velocity. The authors encourage the use of the 1D and 2D formulations of the problem introduced in this section for future benchmark investigation concerning dynamic source models.

\section{Conclusions}\label{Conclusions}
This work introduces a high-order spectral solver based on a Fourier continuation (FC) methodology for solving the nonlinear shallow water equations with dynamic seabed motion. A number of numerical experiments attest to fifth-order spatial accuracy while minimal numerical dispersion errors, making it well-suited for long-distance or high-frequency tsunami propagation modeling. 
% The solver was successfully validated against analytical solutions, experimental data, and classical tsunami benchmarks, confirming its ability to capture critical physical effects relevant to tsunami generation and propagation.
Additionally, several classical and semi-classical benchmarks are provided to validate the present implementation of the solver against analytical solutions as well as experimental data, demonstrating its capability in capturing critical physical effects relevant to tsunami generation and propagation with minimal numerical errors. A number of results are presented in comparison to other solvers, including high-order finite volume and finite difference methods, where the FC-based approach is highly competitive, accurately preserving waveforms over extended spatial and temporal scales at a much lower computational cost than other considered solvers. Furthermore, a new problem configuration, inspired by realistic ground motion settings, is proposed and highlights the relevance of considering dynamic ground motion in tsunami modeling and the efficacy of employing an FC-based approach in such contexts. This includes a preliminary parametric study (a first, to the best of the author's knowledge) of the effects of earthquake speed on tsunamigenesis and the subsequent wave dynamics. 

% a first explicit demonstration of the relevance of considering dynamic ground motion for tsunamigenesis and the efficacy of employing an FC-based approach in such contexts.

% The computational results provided include a number of comparisons to other solvers, including high-order finite volume and high-order finite difference methods, that demonstrate the competitiveness of the new FC solver. 

% Compared to conventional finite difference and finite element methods, the FC-based approach accurately preserves waveforms over extended spatial and temporal scales without requiring excessive computational resources. 

The FC-based approach can be further improved by the addition of explicit shock-capturing schemes (e.g., ENO or WENO~\cite{shu2006essentially,shahbazi2011multi}), although the present implementation is demonstrated in this paper to be reasonably effective without any special treatment. Future work also entails the full coupling of the solver with dynamic rupture models and global-scale bathymetric datasets. Additionally, the applicability of the solver in operational tsunami warning systems can be enhanced by incorporating curvilinear mesh generation strategies, which facilitate adaptive meshing as well as treatment of more complex and realistic geometries~\cite{amlani2016}. The latter may involve more large-scale configurations that require optimization of the solver for high-performance computing (HPC) architectures, such as those with GPU acceleration (the computational complexity here being dominated by the FFT).

Recent contributions~\cite{Elbanna2021,amlanibhat2022} have demonstrated the importance of considering dynamic seafloor motion (in terms of time-dependent displacement and velocity) in the generation of tsunamis due to strike-slip earthquakes. The results of \autoref{sec:Application} extend this idea and imply the importance in capturing the interplay of earthquake dynamics and tsunami propagation even for subduction earthquakes (which are conventionally modeled as instantaneous seafloor displacements). The parameter studies presented herein suggest that such considerations are particularly relevant when the speed of the ground surface waves are comparable to the characteristic wavespeed of the resulting tsunami waves, revealing an unexplored potential for improving tsunami hazard assessment.
% through providing accurate and high-fidelity simulations for further parametric studies. 
The robust and efficient nature of the proposed FC solver, in terms of considering such motions with minimal tuning of numerical parameters, can be highly beneficial for early warning systems and for real-time forecasting, including for the training of physics-based machine learning surrogate models~\cite{Andraud2022} with physically-faithful simulated data.

% Such a solverRecent interest has also been growing in Physics-Informed Neural Networks (PINNs) for accelerating tsunami simulations by approximating solutions at the cost of reduced accuracy . While PINNs offer rapid inference capabilities, they may not yet achieve the high fidelity required for precise tsunami hazard assessments. 

% The solver was successfully validated against analytical solutions, experimental data, and classical tsunami benchmarks, confirming its ability to capture critical physical effects relevant to tsunami generation and propagation. The final 

\section*{Acknowledgments}
T.M., H.S.B., F.A.: this project has received financial support from the CNRS through the MITI interdisciplinary programs. H.S.B.: acknowledges additional support from the European Research Council (ERC) PERSISMO, Grant 865411. F.A.: acknowledges additional support from the Agence Nationale de la Recherche (ANR), Grant  ANR-23-CE46-0008. The authors are grateful to Prof A.E. Elbanna and Dr M. Abdelmeguid for sharing some experimental data, and to Prof M. Derakhti for valuable discussions concerning the benchmark considered in \autoref{piston benchmark}.

% \section*{Data availability}

% All figure and simulation data, including for the benchmarks, is available at: \href{}{DOI:XXXXX}.

% \bibliographystyle{elsarticle-num}

\bibliography{bibliography} %-->reference list is on the template.bib file

\end{document}